\documentclass[11pt,a4paper]{article}

\usepackage[english]{babel}
\usepackage{amsmath}
\usepackage{amssymb}
\usepackage{graphicx}        

\usepackage{cases}
\usepackage{amsfonts}
\usepackage{color}
\usepackage{verbatim}
\usepackage{caption}
\usepackage{subcaption}
\usepackage{fullpage}

\usepackage{array}

\usepackage{epsfig}
\usepackage{float}
\usepackage{url}

\usepackage{theorem}

\usepackage{fullpage}

\usepackage{hyperref}

\usepackage{lmodern} 
\usepackage{amsmath}
\usepackage{amssymb}
\usepackage{graphicx,scalerel}
\usepackage{mathtools} 

\newcommand\nwhat[1]{\hstretch{2}{\hat{\hstretch{.5}{#1}}}}

\def\overbigdot#1{\overset{\hbox{\tiny$\bullet$}}{#1}}
\def\twobigdot#1{\overset{\hbox{\tiny$\bullet \bullet$}}{#1}}

\definecolor{viola}{rgb}{0.3,0,0.7}
\definecolor{ciclamino}{rgb}{0.5,0,0.5}
\definecolor{rosso}{rgb}{0.85,0,0}

\def\pier #1{{\color{black} #1}}
\def\abramo #1{{\color{black} #1}}
\def\abramon #1{{\color{black} #1}}
\def\abramonn #1{{\color{black} #1}}

\def\abramonew #1{{\color{black} #1}}
\def\pco #1{{\color{black} #1}}
\def\abramohh #1{{\color{black} #1}}
\def\abramohhb #1{{\color{black} #1}}
\def\pierhhb #1{{\color{black} #1}}
\def\michhhb #1{{\color{black} #1}}
\def\ap #1{{\color{black} #1}}

\usepackage{stackengine}

\newtheorem{rem}{Remark}[section]
\newtheorem{lem}{Lemma}[section]
\newtheorem{thm}{Theorem}[section]

\newenvironment{pf}{\noindent{\bf Proof. \/}\noindent%
}{\hfill\EndProofMarker}
\newcommand{\EndProofMarker}{$\Box$}
\newcommand{\diver}{\ensuremath{\operatorname{div}}}
\newcommand{\curl}{\ensuremath{\operatorname{curl}}}

\newcommand{\grad}{\ensuremath{\operatorname{grad}}}

\def\XXint#1#2#3{{\setbox0=\hbox{$#1{#2#3}{\int}$ }
\vcenter{\hbox{$#2#3$ }}\kern-.6\wd0}}

\begin{document}
\newcommand{\mic}{\color{blue}}
\newcommand{\mich}{\color{black}}
\abramohhb{ 
\title{ \michhhb{Large deformations in terms of stretch and rotation \\ \ap{and local solution to} the non-stationary problem}}
\date{}
\maketitle 
\begin{center}
\vskip-1.2cm
{\large\sc Abramo Agosti$^{(1)}$}\\
{\normalsize e-mail: {\tt abramo.agosti@unipv.it}}\\[0.25cm]
{\large\sc Michel Fr\'emond$^{(2)}$}\\
{\normalsize e-mail: {\tt michel.fremond@uniroma2.it}}\\[.5cm]
$^{(1)}$
{\small Dipartimento di Matematica ``F. Casorati'', Universit\`a di Pavia}\\
{\small via Ferrata 5, I-27100 Pavia, Italy}\\[.3cm] 
$^{(2)}$
{\small Lagrange Laboratory}\\
{\small \michhhb{ Dipartimento di Ingegneria Civile e Ingegneria Informatica}, \\Universit\`a di Roma ``Tor Vergata"}\\
{\small Via del Politecnico, 1, 00163 Roma, Italy}
\end{center}
%
%
%
%
%
\begin{abstract}
In this paper we consider and generalize a model, recently proposed and analytically investigated in its quasi-stationary approximation by the authors, for \michhhb{visco-elasticity} with large deformations and conditional compatibility, where the independent variables are the stretch and the rotation tensors. The model takes the form of a system of integro-differential coupled equations. Here, its derivation is generalized to consider mixed boundary conditions, which may represent a wider range of physical applications then the case with Dirichlet boundary conditions considered in our previous contribution. This also introduces nontrivial technical difficulties in the theoretical framework, related to the definition and the regularity of the solutions of elliptic operators with mixed boundary conditions.
As a novel contribution, we develop the analysis of the fully non-stationary version of the system where we consider inertia.  In this context, we prove the existence of a local in time weak solution in three space dimensions, employing techniques from PDEs and convex analysis.
\end{abstract}
}

\section{Introduction}
In this paper we consider and generalize a model for \michhhb{visco-elasticity} with large deformations, formulated in terms of the stretch and rotation tensors, recently introduced in \cite{agosti4}. In this model, the stretch tensor satisfies a conditional compatibility condition, with occurrence of defects depending on the magnitude $k$ of an internal force. When defect occurs, the kinematic relation between the stretch tensor $\mathbf{W}$, the rotation tensor $\mathbf{R}$, the displacement vector $\vec{u}$ and the defect tensor $\mathbf{Z}$ is the following:
\begin{gather}
\label{eqn:1intro}
\mathbf{RW-I=}\grad \vec{u} \mathbf{+}\curl \mathbf{Z,} \quad 
\diver \mathbf{Z}=0,
\end{gather}
where $\mathbf{I}$ is the identity tensor.
\newline
The model was derived in \cite{agosti4} in the situation of {\mich a visco-elastic} body fixed at its boundary, employing dissipative principles and starting from a non-standard form of the principle of virtual powers, where the virtual velocities satisfy an internal constraint similar to \eqref{eqn:1intro} depending on the solutions of the problem. In particular, thanks to a Helmholtz--Hodge decomposition for second order tensors, the virtual velocities associated to the displacement vector and the defect tensors were defined to satisfy elliptic problems with homogeneous Dirichlet boundary conditions and with data given in terms of the stretch and rotation tensors and of their associated virtual velocities. Using the technique of Green functions, we inverted the internal constraints and reduced the set of independent virtual velocities to the ones associated to the stretch and rotation tensors, obtaining a system of integro-differential coupled equations in the stretch and rotation variables only. Once the model is solved, the displacement vector and the defect tensor are constructed by solving back the elliptic problems associated to the Helmholtz decomposition of $\mathbf{RW-I}$ in \eqref{eqn:1intro}.
In \cite{agosti4} we developed the analysis for the quasi-stationary approximation of the full system, i.e. neglecting inertia, obtaining the existence of a global in time strong solution in three space dimensions, with uniqueness and continuous dependence from data in the limit of full incompatibility in the system.
\newline
In the present contribution, we extend the model proposed in \cite{agosti4} by considering the more physically appropriate situation of {\mich a visco-elastic} body which is fixed only on a part of its boundary, and which is free to move on the other part. This generalization is not straightforward, since it requires specific mathematical developments to fit the theoretical framework to the present case, in particular the Helmholtz--Hodge decomposition and the definition of Green functions, associated to mixed boundary conditions. Also, we develop the analysis for the actual basic mechanical problem, which is the full non-stationary problem, i.e. considering inertia. As in \cite{agosti4}, the inertia of the system is expressed by a virtual power of acceleration forces containing second-order interaction terms in space, which allows us to obtain sufficient regularity of weak solutions to be able to represent the physical situation of a contact at a point with inertia, for instance, punctual force and torque applied to the head of a nail hammered in a wall. We also impose the positive definiteness of the stretch matrix as an internal constraint in the free energy of the system, which implies that the material is not flattening or crushing.
In the case with inertia, we will find existence of a weak solution in three space dimensions only locally in time, as long as the solution remains continuously in the interior of the proper domain of the indicator function associated to the positivity constraint, preceding the possible realization of external and internal collisions. Hence, this result complies with the mechanical situations of a {\mich visco-elastic} solid possibly undergoing collisions.
\newline
The main technical difficulties in this study are that we need to require non standard mixed boundary conditions for $\mathbf{RW}$ in \eqref{eqn:1intro} if the displacement and the defect variables are solutions of elliptic problems with mixed boundary conditions, and further that no elliptic regularity tools are available in general for elliptic problems with mixed boundary conditions. The latter tools were heavily employed in the analysis of the quasi-stationary problem with homogeneous Dirichlet boundary conditions studied in \cite{agosti4}. Also, the presence of inertia terms complicate the analysis, since they contain higher order time derivatives of the variables nonlinearily coupled with lower order terms, hence hardening the proof of the existence of a solution even at the discrete level. Finally, the inertia terms prevent the uniform control of the subdifferential of the indicator function associated to the positivity constraint in some $L^p$ space, which lead to global in time existence in the quasi-stationary case, while in the present situation we only obtain local in time existence.
\newline 
The limit with full incompatibility in the system, i.e., the case the threshold $k$ for defects to appear is null, will be addressed in a forthcoming paper, together with the possible presence of collisions in the dynamics.
\newline

{
In mechanical parlance, our parti pris is to describe the motion with
stretch matrix $\mathbf{W}$ and rotation matrix $\mathbf{R}$. If they are
known displacement $\vec{u} $ should be given by polar decomposition
\begin{gather}
\mathbf{RW-I=}\grad \vec{u}. \notag
\end{gather}
But there is no fundamental reason for $\mathbf{RW}$ to be a
gradient. Difficulty is overcome by mechanics which experiments defects, and by mathematics which proves with convenient
boundary conditions, that for any matrix $\mathbf{RW}$ there exist matrix $\mathbf{Z}$ accounting for the defects and displacement $\vec{u}$ which satisfy \ref{eqn:1intro}. These quantities are given by Green functions. For their part, matrices $\mathbf{W}$ and $\mathbf{R}$ are given by linear and angular equations of motion
resulting from the principle of virtual power and constitutive laws, one of
them reporting the experimental result that there are no defects if some
stress is not too large.
}

\pier{The paper is organized as follows. In Section \ref{sec:notation} we introduce the necessary notation and some preliminary results, which let us extend the theoretical framework introduced in \cite{agosti4} to the case with mixed boundary conditions. In Section \ref{sec:model} we derive the full model with inertia and mixed boundary conditions. In Section \ref{sec:analysis} we study the existence problem for the full non-stationary problem. We conclude with some observations and future perspectives in Section \ref{sec:conclusions}.}

\section{Notations and preliminaries}
\label{sec:notation}
In this section we introduce the notation and the preliminary results about the functional setting which will be necessary for the model derivation.

\subsection{Geometrical and functional setting}
Let $\mathcal{D}_a\subset \mathbb{R}^3$ be an open bounded and simply connected domain with Lipschitz boundary $\Gamma_a:=\partial \mathcal{D}_a$, with associated unit normal $\vec{n}$, and let $[0,T]$ be a finite time interval, with $T>0$. We introduce the notation $\mathcal{D}_{aT}:= \mathcal{D}_a\times [0,T]$. We remark that the assumption of a simply connected domain is made only to maintain the theoretical framework as simple as possible, in particular for what concerns the characterization of the kernels of the $\curl$ and the $\diver$ operators, which will be used later. The theoretical framework could be extended in a standard manner to consider a connected but not simply connected domain, as it will be further detailed in the forthcoming sections (see Remark \ref{rem:hhsc}).

 We indicate as $M(\mathbb{R}^{3\times 3})$ the linear space of square matrices, endowed with the Frobenius inner product
\[
\mathbf{A}: \mathbf{B}=\sum_{i,j=1}^3\mathbf{A}_{ij}\mathbf{B}_{ij},
\]
for any $\mathbf{A},\mathbf{B}\in M(\mathbb{R}^{3\times 3})$. Tensors are indicated with capital boldface letters, while vectors are indicated by lowercase letters with an arrow superscript.
We also indicate with the notation $: :$ the Frobenius inner product in $M(\mathbb{R}^{3\times 3\times 3})$, and with the notation $: : :$ the Frobenius inner product in $M(\mathbb{R}^{3\times 3\times 3\times 3})$. \pier{The orthogonal subspaces of symmetric and antisymmetric matrices are denoted by} \pier{$Sym(\mathbb{R}^{3\times 3})\subset M(\mathbb{R}^{3\times 3})$ and $Skew(\mathbb{R}^{3\times 3})\subset M(\mathbb{R}^{3\times 3})$}, 
respectively. We indicate the set of special orthogonal matrices as $SO(\mathbb{R}^{3\times 3})$ and the set of positive definite symmetric matrices as  $Sym^+(\mathbb{R}^{3\times 3})$. \abramo{We recall that for any $\mathbf{R}\in SO(\mathbb{R}^{3\times 3})$ there exists a unique $\mathbf{A}\in Skew(\mathbb{R}^{3\times 3})$ such that $\mathbf{R}=\mathrm{e}^{\mathbf{A}}$, where the exponential of a matrix must be intended as $\mathrm{e}^{\mathbf{A}}=\sum_{n=0}^{\infty}\frac{\mathbf{A}^n}{n!}$.}
\pier{For a generic subset $K \subset M(\mathbb{R}^{3\times 3})$, let $I_K : M(\mathbb{R}^{3\times 3}) \to \{0, +\infty\}$  denote the indicator function of $K$, which is defined, for any $\mathbf{A} \in M(\mathbb{R}^{3\times 3})$, by  $ I_k(\mathbf{A})=0 $ if $\mathbf{A} \in K$, $ I_k(\mathbf{A})= +\infty $ if $\mathbf{A} \not\in K$.}

We introduce the space of vector fields $\mathcal{V}:=(\mathbb{R}^3)^{\mathcal{D}_{aT}}$, whose elements are functions from $\mathcal{D}_{aT}$ to $\mathbb{R}^3$. We further introduce the spaces of tensor fields $\mathcal{M}:=(M(\mathbb{R}^{3\times 3}))^{\mathcal{D}_{aT}}$, $\mathcal{SO}:=(SO(\mathbb{R}^{3\times 3}))^{\mathcal{D}_{aT}}$, $\mathcal{S}:=(Sym(\mathbb{R}^{3\times 3}))^{\mathcal{D}_{aT}}$ and $\mathcal{A}:=(Skew(\mathbb{R}^{3\times 3}))^{\mathcal{D}_{aT}}$, with $\mathcal{M}=\mathcal{S}\oplus \mathcal{A}$. Given a tensor $\mathbf{A}\in \mathcal{M}$, we denote by $\textrm{Sym}(\mathbf{A}):=\frac{\mathbf{A}+\mathbf{A}^T}{2}$ its symmetric part and by $\textrm{Skew}(\mathbf{A}):=\frac{\mathbf{A}-\mathbf{A}^T}{2}$ its antisymmetric part. We also need to introduce the space of tensor fields $\mathcal{M}_{\diver}:=\{\mathbf{A}\in \mathcal{M}|\diver \mathbf{A}=\mathbf{0}\}$, where the \pco{\textit{divergence}} of a second order tensor is defined row wise. In the following, we will operate also with the \textit{curl} of second order tensors, which is defined row wise. 

We denote by $L^p(\mathcal{D}_a,K)$ and $W^{r,p}(\mathcal{D}_a,K)$ the standard Lebesgue and Sobolev spaces of functions defined on $\mathcal{D}_a$ with values in a set $K$, {where $K$ may be $\mathbb{R}$ or a 
multiple power of $\mathbb{R}$,} and by  $L^p(0,t;V)$ the Bochner space of functions defined on $(0,t)$ with values in the functional space $V$, with $1\leq p \leq \infty$. If $K\equiv \mathbb{R}$, we simply write  $L^p(\mathcal{D}_a)$ and $W^{r,p}(\mathcal{D}_a)$. 
For a normed space $X$, the associated norm is denoted by $\pier{\|}\cdot\pier{\|}_X$. In the case $p=2$, we use the notations $H^r:=W^{r,2}$, and we denote by $(\cdot,\cdot)$ and $\pier{\|}\cdot\pier{\|}$ the $L^2$ scalar product and induced norm between functions with scalar, vectorial or tensorial values. Moreover, we denote by $C^k(\mathcal{D}_a;K)$ the spaces of continuously differentiable functions up to order $k$ defined on $\mathcal{D}_a$ with values in a set $K$\pco{; by}  
{$C^k([0,t];V)$, $k\geq 0$, the spaces of {continuously differentiable functions up to} order $k$ from $[0,t]$ to the space $V$. The dual space of a Banach space $Y$ is denoted by $Y'$, and their dual product is indicated as $\prescript{}{Y'}<\cdot,\cdot>_Y$. We denote by $W_{\Sigma}^{r,p}
(\mathcal{D}_a;K)$, $r>\frac{1}{2},p\geq 1$, the space of functions in $W^{r,p}
(\mathcal{D}_a;K)$ with zero trace on a Lipschitz continuous subset $\Sigma$ of $\Gamma_a$ with positive measure $|\Sigma|>0$.
As before, when $p=2$ we will indicate 
the latter functional space as $H_{\Sigma}^r(\mathcal{D}_a;K)$. {We will also need the space $ \bar{H}^2_{\Sigma}(\mathcal{D}_a,K):=\{\mathbf{f}\in H^2(\mathcal{D}_a,K)|\mathbf{f}=\mathbf{0},\; \grad \mathbf{f}=\mathbf{0} \;\text{on}\;\Sigma\}$}.}
The traces of functions $H^r(\mathcal{D}_a;K)$ on $\Sigma$ belong to the space $H^{r-1/2}(\Sigma;K)$, for $r>\frac{1}{2}$ (see e.g. \cite{LM}).
We moreover introduce the space $H_{00}^{r}(\Sigma,K):=\{v\in L^2(\Sigma,K): \, \tilde{v}\in H^{r}(\Gamma_a,K)\}$, for $r\geq 1/2$, where $\tilde{v}$ is the extension by zero of $v$ {to the set $\Gamma_a\setminus \Sigma$}.  

We also need to introduce the spaces
\begin{align*}
& L_{\diver}^2(\mathcal{D}_a,K):=\overline{\{\vec{u}\in C^{\infty}(\mathcal{D}_a,K): \, \diver\vec{u}=0 \; \text{in} \; \mathcal{D}_a\}}^{\pier{\|}\cdot\pier{\|}_{L^2(\mathcal{D}_a;K)}},\\
&H(\diver;\mathcal{D}_a,K):=\{\vec{u}\in L^2(\mathcal{D}_a,K): \, \diver\vec{u}\in L^2(\mathcal{D}_a,\hat{K})\},\\
&H(\curl;\mathcal{D}_a,K):=\{\vec{u}\in L^2(\mathcal{D}_a,K): \, \curl \vec{u}\in L^2(\mathcal{D}_a,K)\},\\
& H_{\Sigma,\diver}^1(\mathcal{D}_a,K):={\{\vec{u}\in H_{\Sigma}^1(\mathcal{D}_a,K): \, \diver\vec{u}=0 \; \text{in} \; \mathcal{D}_a\}},
\end{align*}
where in the second definition $K=\mathbb{R}^d$ or $K=\mathbb{R}^{d\times d}$ and $\hat{K}=\mathbb{R}$ or $\hat{K}=\mathbb{R}^{d}$ respectively. The normal traces of functions $H(\diver;\mathcal{D}_a,K)$ on $\Sigma$ belong in general to the space $\left(H_{00}^{1/2}(\Sigma)\right)'$, with $H^{-1/2}(\Sigma)$ continuously embedded in $\left(H_{00}^{1/2}(\Sigma)\right)'$ (see e.g. \cite{LM}). 
The duality pairing between $H_{\Sigma,\diver}^1(\mathcal{D}_a;K)$ and $(H_{\Sigma,\diver}^1\left(\mathcal{D}_a;K)\right)^{\prime}$ is still denoted by $<\cdot,\cdot>$. 

In the following, $C$ denotes a generic positive constant independent of the unknown variables, the discretization and the \michhhb{physical parameters}, the value of which might change from line to line; $C_1, C_2, \dots$ indicate generic positive constants whose particular value must be tracked through the calculations; $C(a,b,\dots)$ denotes a constant depending on the nonnegative parameters $a,b,\dots$.

\subsection{Green functions with mixed boundary conditions}
We endow the space 
$H_{\Sigma}^1(\mathcal{D}_a;K)$ with the inner product $(A,B)_{H_{\Sigma}^1(\mathcal{D}_a;K)}:=(\grad A,\grad B)$, for all $A,B \in H_{\Sigma}^1(\mathcal{D}_a;K)$, and we introduce the Riesz isomorphism $\mathcal{R}:H_{\Sigma}^1(\mathcal{D}_a;K)\rightarrow 
\left(H_{\Sigma}^1(\mathcal{D}_a;K)\right)'$, defined by
\[
<\mathcal{R}A,B>=(A,B)_{H_{\Sigma}^1(\mathcal{D}_a;K)}, \quad \forall A,B\in H_{\Sigma}^1(\mathcal{D}_a;K).
\]
The operator $\mathcal{R}=-\Delta$ is the negative weak Laplace operator with homogeneous Dirichlet boundary conditions on $\Sigma$ and homogeneous Neumann boundary conditions on $\Gamma_a\setminus \Sigma$, which is positive definite and self adjoint. As a consequence of the Lax--Milgram theorem and the Poincar\'e inequality, the inverse operator $(-\Delta)^{-1}:\left(H_{\Sigma}^1(\mathcal{D}_a;K)\right)'\rightarrow H_{\Sigma}^1(\mathcal{D}_a;K)$ is well defined, and we set $A:=(-\Delta)^{-1}F=\mathcal{G}_L\ast F$, for $F\in \left(H_{\Sigma}^1(\mathcal{D}_a;K)\right)'$, where $\mathcal{G}_L$ is the Green function associated to the Laplace operator with 
mixed boundary conditions and $\ast$ denotes the convolution operation, if $-\Delta A=F$ in $\mathcal{D}_a$ in the weak sense, and $A=0$ on $\Sigma$ in the sense of traces. The following Lax--Milgram estimate is valid
\begin{equation}
\label{eqn:1}
\pier{\|}A\pier{\|}_{H_{\Sigma}^1(\mathcal{D}_a;K)}\leq C\pier{\|}F\pier{\|}_{\left(H_{\Sigma}^1(\mathcal{D}_a;K)\right)'}.
\end{equation}
\begin{rem}
\label{ellipticreg} 
In the case with a smooth boundary $\Gamma_a$ and regular $F$, e.g. $\Gamma_a$ of class $C^{m+2}$ and $F\in W^{m,p}(\mathcal{D}_a;K)$, $1<p<\infty$, $m\in \mathbb{N}$, the solution $A\in H_{\Gamma_a}^1(\mathcal{D}_a;K)$ of the elliptic problem $-\Delta A=F$ with homogeneous Dirichlet boundary conditions over $\Gamma_a$ satisfies the elliptic regularity property that $A\in W^{m+2,p}(\mathcal{D}_a;K)$ and $-\Delta A=F$ a.e. in $\mathcal{D}_a$ \cite[Chapter 9]{LM}. In the case with mixed Dirichlet and Neumann boundary conditions, elliptic regularity is locally preserved in the neighborhoods of interior points of $\Sigma$ and $\Gamma_a\setminus \Sigma$, but globally $A\notin W^{m+2,p}(\mathcal{D}_a;K)$ \cite{MS}. For instance, in the case of mixed homogeneous boundary conditions $A\in W^{s,p}(\mathcal{D}_a;K)$, with $s<1/2+2/p$, even in presence of smooth data (see e.g. \cite[Chapter 3, Remark 3.3]{bg}, \cite{savare} and the references therein). Elliptic regularity with mixed boundary conditions is globally valid only in particular cases, for instance when $\bar{\Sigma}\cap (\overline{\Gamma_a\setminus \Sigma})= \emptyset$, e.g. when the solid has the form of a three dimensional annulus with mixed boundary conditions applied separately to the internal and external connected components of the boundary, or in some situations when the solid has a convex shape or is a polyedron \cite[Chapter 8]{grisvard}. Since we want to describe the most general case, in our treatment we will not have at our disposal elliptic regularity instruments. This is a severe technical limitation with respect to the study developed in \cite{agosti4}, where we treated the quasi-stationary problem with homogeneous Dirichlet boundary conditions over the whole domain.
\end{rem} 
Similarly, we can introduce the Riesz isomorphism $\mathcal{R}_{\diver}:H_{\Sigma,\diver}^1(\mathcal{D}_a;K)\rightarrow \left(H_{\Sigma,\diver}^1(\mathcal{D}_a;K)\right)^{\prime}$, defined by
\[
<\mathcal{R}_{\diver}A,B>=(\grad A,\grad B), \quad \forall A,B\in H_{\Sigma,\diver}^1(\mathcal{D}_a;K).
\]
The operator $\mathcal{R}_{\diver}=-P_L\Delta$, where $P_L:L^2(\mathcal{D}_a;K)\to L_{\diver}^2(\mathcal{D}_a;K)$ \pier{denotes the Leray projector, is the} negative projected Laplace operator with homogeneous Dirichlet boundary conditions on $\Sigma$ and homogeneous Neumann boundary conditions on $\Gamma_a\setminus \Sigma$, which is positive definite and self adjoint. As a consequence of the Lax--Milgram theorem and the Poincar\'e inequality, the inverse operator $(-P_L\Delta)^{-1}:\left(H_{\Sigma,\diver}^1(\mathcal{D}_a;K)\right)^{\prime}\rightarrow H_{\Sigma,\diver}^1(\mathcal{D}_a;K)$ is well defined, and we set $A:=(-P_L\Delta)^{-1}F=\mathcal{G}_{L,\diver}\ast F$, for $F\in \left(H_{\Sigma,\diver}^1(\mathcal{D}_a;K)\right)^{\prime}$, where $\mathcal{G}_{L,\diver}$ is the Green function associated to the projected 
Laplace operator with mixed boundary conditions, if $-P_L\Delta A=F$ in $\mathcal{D}_a$ in the weak sense, and $A=0$ on $\Sigma$ in the sense of traces. We again note that, if $A\in H_{\Sigma,\diver}^1(\mathcal{D}_a;K)$ solves $-P_L\Delta A=F$ for some $F\in W^{m,p}(\mathcal{D}_a;K)\cap L_{\diver}^2(\mathcal{D}_a,K)$, $1<p<\infty$, $m\in \mathbb{N}$, and $\Gamma_a$ is of class $C^{m+2}$, in general $A\notin W^{m+2,p}(\mathcal{D}_a;K)\cap L_{\diver}^2(\mathcal{D}_a,K)$.

\subsection{Helmholtz--Hodge decomposition for vector fields with mixed boundary conditions.}
We now give specific forms of the Helmholtz--Hodge decomposition with mixed boundary conditions which will be useful in the forthcoming sections. Similar decompositions were introduced in \cite{gilardimixed} and proved by topological arguments. Here, we report the proof of such decompositions obtained through a constructive procedure by means of the solution of elliptic problems, following similar arguments introduced in \cite[Chapter IX]{dautraylions} for the case of standard (i.e. Dirichlet or Neumann) boundary conditions. The aforementioned elliptic problems will be crucial in the model derivation to define the \textit{kinematic constraints} between the model variables.  
\begin{thm}
\label{thm:hhd}
Let $\mathcal{D}_a\subset \mathbb{R}^3$ be an open bounded and simply connected domain with Lipschitz boundary $\Gamma_a:=\partial \mathcal{D}_a$. Let us assume that $\Gamma_a=\Gamma_D\cup \Gamma_N$, where $\Gamma_D,\Gamma_N$ are connected and Lipschitz continuous subsets of $\Gamma_a$ with positive measures and such that $|\Gamma_D\cap \Gamma_N|=0$.
Let us introduce the spaces
\begin{align*}
& \grad H^1_{c,\Gamma_D}:=\{\vec{w}\in L^2(\mathcal{D}_a,\mathbb{R}^3): \exists p \in H^1(\Omega) \; \text{such that}\; \vec{w}=\grad p, \;\; p|_{\Gamma_D}=c\}, \\
& \curl H^1_{\Gamma_N,div}:=\{\vec{w}\in L^2(\mathcal{D}_a,\mathbb{R}^3): \exists \vec{v} \in H^1(\mathcal{D}_a,\mathbb{R}^3) \; \text{such that}\; \vec{w}=\curl \vec{v},\;\; \diver \vec{v}=0,\;\; \vec{v}|_{\Gamma_N}=\vec{0}\}, \\
& H_{\Gamma_N,div}:=\{\vec{v}\in L^2(\mathcal{D}_a,\mathbb{R}^3): \textrm{div}\vec{v}=0, \; \vec{v}\cdot \vec{n}|_{\Gamma_N}=0\},\\
& H_{\Gamma_D,curl}:=\{\vec{v}\in L^2(\mathcal{D}_a,\mathbb{R}^3): \textrm{curl}\vec{v}=\vec{0}, \; \vec{v}\wedge \vec{n}|_{\Gamma_D}=\vec{0}\},
\end{align*}
where $c\in \mathbb{R}$ is a constant.
For any $\vec{\xi}\in L^2(\mathcal{D}_a,\mathbb{R}^3)$, there exist unique $p \in H^1_{\Gamma_D}$ and $\vec{r}\in H_{\Gamma_N,div}$ such that
\begin{equation}
\label{hh5}
\vec{\xi}=\grad p + \vec{r},
\end{equation}
i.e. the following decomposition is valid 
\begin{equation}
\label{hh5bis}
L^2(\mathcal{D}_a,\mathbb{R}^3)=\grad H^1_{0,\Gamma_D}\oplus H_{\Gamma_N,div}.
\end{equation}
\end{thm}
Moreover, there exist a $\vec{v}\in \grad H^1_{c,\Gamma_D}$, defined for any $c$, with $\vec{v}=\grad u$, and a $\mathbf{q}\in \curl H^1_{\Gamma_N,div}$, with $\mathbf{q}=\curl \vec{\alpha}$, such that 
\begin{equation}
\label{hh6}
\vec{\xi}=\grad u + \curl \vec{\alpha},
\end{equation}
i.e. the following decomposition is valid 
\begin{equation}
\label{hh6bis}
L^2(\mathcal{D}_a,\mathbb{R}^3)= \grad H^1_{c,\Gamma_D}\oplus \curl H^1_{\Gamma_N,div},
\end{equation}
where $u$ and $\vec{\alpha}$ satisfy the elliptic problems with mixed boundary conditions
\begin{equation}
\label{phh6}
\begin{cases}
\Delta u = \textrm{div}\vec{\xi},\\
u|_{\Gamma_D}=c,\\
\grad u \cdot \vec{n}|_{\Gamma_N}=\vec{\xi}\cdot \vec{n}|_{\Gamma_N},
\end{cases}
\begin{cases}
-\Delta \vec{\alpha} = \text{curl}\vec{\xi},\;\; \diver \vec{\alpha}=0,\\
\text{curl} \,\vec{\alpha}\wedge \vec{n}|_{\Gamma_D}=\vec{\xi}\wedge \vec{n}|_{\Gamma_D},\\
\vec{\alpha}|_{\Gamma_N}=0.
\end{cases}
\end{equation}
\begin{rem}
\label{rem:hhsc}
The hypotheses that  $\mathcal{D}_a$ is simply connected and that $\Gamma_D,\Gamma_N$ are connected are made to simplify the presentation of the results. The theorem could be extended in a standard manner to a (not simply) connected domain $\mathcal{D}_a$ with boundary subsets $\Gamma_D,\Gamma_N$ constituted by a finite number of connected components $\Gamma_{D,i},\Gamma_{N,j}$, $i=1,\dots,N$, $j=1,\dots,M$, introducing a finite number of cuts $\Sigma_l$, $l=1,\dots,L$, in the domain constituted by connected orientable Lipschitz submanifolds, glued topologically to the boundary of connected elements of $\Gamma_D,\Gamma_N$, and reducing the analysis to the simply connected subdomain $\mathcal{D}_a\setminus \bigcup_{l=1}^L\Sigma_l$ (see e.g. \cite{dautraylions,gilardimixed}). Anyhow, the situation in which the initial form of the body is topologically simply connected until it develops cuts or holes is mechanically meaningfull.
\end{rem}
\begin{pf}
We start by introducing the following preliminary Green formulas:
\begin{multline}
 \label{green1}
 \int_{\mathcal{D}_a}(\diver \vec{u})\,v =-\int_{\mathcal{D}_a}\vec{u}\cdot \grad v + \prescript{}{\left(H_{00}^{1/2}(\Gamma_N)\right)'}{<}\vec{u}\cdot \vec{n},v>_{H_{00}^{1/2}(\Gamma_N)}, \\ \forall \vec{u}\in H(\diver;\mathcal{D}_a,\mathbb{R}^3), \; v\in H_{\Gamma_D}^1(\mathcal{D}_a),
\end{multline}
\begin{multline}
 \label{green2}
 \int_{\mathcal{D}_a}\curl \vec{u}\cdot \vec{v} =\int_{\mathcal{D}_a}\vec{u}\cdot \curl \vec{v} - \prescript{}{\left(H_{00}^{1/2}(\Gamma_D,\mathbb{R}^3)\right)'}{<}\vec{u}\wedge \vec{n},\vec{v}>_{H_{00}^{1/2}(\Gamma_D,\mathbb{R}^3)}, \\ 
 \forall \vec{u}\in H(\curl;\mathcal{D}_a,\mathbb{R}^3), \; \vec{v}\in H_{\Gamma_N}^1(\mathcal{D}_a,\mathbb{R}^3).
\end{multline}
We then prove that $\grad H^1_{0,\Gamma_D}\perp H_{\Gamma_N,div}$ in $L^2(\mathcal{D}_a,\mathbb{R}^3)$. Indeed, for any $\vec{w}\in \grad H^1_{0,\Gamma_D}$, $\vec{v}\in H_{\Gamma_N,div}$, we have from \eqref{green1} that
\[
\int_{\mathcal{D}_a}\vec{w}\cdot \vec{v} =\int_{\mathcal{D}_a}\grad p\cdot \vec{v} =\prescript{}{\left(H_{00}^{1/2}(\Gamma_N)\right)'}{<}\vec{v}\cdot \vec{n},p>_{H_{00}^{1/2}(\Gamma_N)}-\int_{\mathcal{D}_a}p\text{div}\vec{v} =0.
\]
We now introduce the following elliptic problems,
\begin{equation}
\label{phh1}
\begin{cases}
\Delta p_1 = \text{div}\vec{\xi},\\
p_1|_{\Gamma_a}=0,
\end{cases}
\begin{cases}
\Delta p_2 = 0,\\
p_2|_{\Gamma_D}=0,\\
\grad p_2 \cdot \vec{n}|_{\Gamma_N}=(\vec{\xi}-\grad p_1)\cdot \vec{n}|_{\Gamma_N},
\end{cases}
\end{equation}
which may be intended in $H^{-1}(\mathcal{D}_a)$ and $\left(H_{\Gamma_D}^1(\mathcal{D}_a)\right)'$ respectively.
There exists a unique solution $p_1\in H_{\Gamma_a}^1(\Omega)$ of the first problem in \eqref{phh1}. Since $\vec{\xi}-\grad p_1$ is an element of $L^2(\mathcal{D}_a,\mathbb{R}^3)$ and has zero divergence, it has a normal trace in the space $\left(H_{00}^{1/2}(\Gamma_N)\right)'$. Then, there exists a unique solution $p_2\in H_{\Gamma_D}^1(\Omega)$, which is an harmonic function. See e.g. \cite{LM} for the latter existence results.  Defining $p:=p_1+p_2\in H_{\Gamma_D}^1(\Omega)$ and $\vec{r}:=\vec{\xi}-\grad p$, we have that $\vec{r}\in H_{\Gamma_N,div}$ and \eqref{hh5} is valid.
We observe that the decomposition \eqref{hh5} is unique. Indeed, given two decompositions $\vec{\xi}=\grad p_a+\vec{r}_a=\grad p_b+\vec{r}_b$, taking the scalar product in $L^2(\mathcal{D}_a,\mathbb{R}^3)$ between their difference and $\vec{r}_a-\vec{r}_b$ and integrating over the domain, we obtain that
\[
||\vec{r}_a-\vec{r}_b||^2+\int_{\mathcal{D}_a}(\grad p_a-\grad p_b)\cdot (\vec{r}_a-\vec{r}_b) =0,
\]
and also, using \eqref{green1} and the facts that $p_a,p_b\in H_{\Gamma_D}^1$ and $\vec{r}_a,\vec{r}_b\in H_{\Gamma_N,div}$, that
\begin{align*}
& \int_{\mathcal{D}_a}(\grad p_a-\grad p_b)\cdot (\vec{r}_a-\vec{r}_b) =\prescript{}{\left(H_{00}^{1/2}(\Gamma_N)\right)'}{<}(\vec{r}_a-\vec{r}_b)\cdot \vec{n},(p_a-p_b)>_{H_{00}^{1/2}(\Gamma_N)}\\
&\quad - \int_{\mathcal{D}_a}(p_a-p_b)\text{div}(\vec{r}_a-\vec{r}_b) =0.
\end{align*}
Hence, $\vec{r}_a\equiv \vec{r}_b$ in $H_{\Gamma_N,div}$ and the decomposition is unique.

We now rewrite \eqref{hh5} in a form which will be usefull in the sequel, since it will let us associate an elliptic problem also to the component in $H_{\Gamma_N,div}$ in the decomposition, after having properly characterized the kernels of the $\curl$ and $\diver$ operators in presence of mixed boundary conditions. First of all we note that
\[
 \grad H_{0,\Gamma_D}^1\subset H_{\Gamma_D,curl},
\]
with both sets being closed in $L^2(\mathcal{D}_a,\mathbb{R}^3)$. This is due to the facts that $\curl \grad p=\vec{0}$ for any $p\in H^1(\mathcal{D}_a)$ and that, if $p|_{\Gamma_D}=0$ in $H^{1/2}(\Gamma_D)$, the tangential derivatives of the trace are null, i.e. $\grad p\wedge \vec{n}|_{\Gamma_D}=\vec{0}$ in $\left(H_{00}^{1/2}(\Gamma_D)\right)'$ and $\grad p\in H_{\Gamma_D,curl}$. The closedness of the sets is proved in \cite[Propositions 6.1-6.2]{gilardimixed}. Then, thanks to \eqref{hh5bis}, we have that
\begin{equation}
 \label{hh7}
 H_{\Gamma_D,curl}=\grad H_{0,\Gamma_D}^1 \oplus \underbrace{H_{\Gamma_N,div}\cap H_{\Gamma_D,curl}}_{:=\mathbb{H}}.
\end{equation}
Note that, if $\vec{u}\in \mathbb{H}$, then $\vec{u}=\grad p$, where $p$ is an harmonic function satisfying the system
\[
 \begin{cases}
  \Delta p=0,\\
  p|_{\Gamma_D}=c,\\
  \grad p\cdot \vec{n}|_{\Gamma_N}=0,
 \end{cases}
\]
i.e. $p$ is a constant. Hence,
\[
 H_{\Gamma_D,curl}=\grad H_{0,\Gamma_D}^1 \oplus \mathbb{H}=\grad H_{c,\Gamma_D}^1.
\]
Since $\mathbb{H}$ is finite-dimensional and hence closed, we may write $H_{\Gamma_N,div}=\mathbb{H}\oplus \mathbb{H}^{\perp}$, whence we rewrite \eqref{hh5bis} as
\begin{equation}
\label{hh5tris}
L^2(\mathcal{D}_a,\mathbb{R}^3)=\grad H^1_{0,\Gamma_D}\oplus \mathbb{H}\oplus \mathbb{H}^{\perp}=H_{\Gamma_D,curl}\oplus \mathbb{H}^{\perp}= \grad H^1_{c,\Gamma_D}\oplus \mathbb{H}^{\perp}.
\end{equation}
Finally, since $H_{\Gamma_D,curl}$ is closed in $L^2(\mathcal{D}_a,\mathbb{R}^3)$, we may write the decomposition
\[
 L^2(\mathcal{D}_a,\mathbb{R}^3)=H_{\Gamma_D,curl}\oplus \left(H_{\Gamma_D,curl}\right)^{\perp}.
\]
The orthogonal complement of $H_{\Gamma_D,curl}$ in $L^2(\mathcal{D}_a)$ is the set $\curl H_{\Gamma_N}^1$. Indeed, taking $\vec{u}\in H_{\Gamma_D,curl}$ and $\vec{v}\in H_{\Gamma_N}^1$, from \eqref{green2} we have that
\[
 \int_{\mathcal{D}_a}\vec{u}\cdot \curl \,\vec{v} =\int_{\mathcal{D}_a}\curl \,\vec{u}\cdot \vec{v} +\prescript{}{\left(H_{00}^{1/2}(\Gamma_D,\mathbb{R}^3)\right)'}{<}\vec{u}\wedge \vec{n},\vec{v}>_{H_{00}^{1/2}(\Gamma_D,\mathbb{R}^3)}=0.
\]
Hence we identify $\mathbb{H}^{\perp}\equiv \curl H_{\Gamma_N}^1$, and rewrite \eqref{hh5tris} as
\begin{equation}
\label{hh8}
L^2(\mathcal{D}_a,\mathbb{R}^3)= \grad H^1_{c,\Gamma_D}\oplus \curl H_{\Gamma_N}^1,
\end{equation}
i.e. for any $\vec{\xi}\in L^2(\mathcal{D}_a,\mathbb{R}^3)$ there exist a  $\vec{v}\in \grad H^1_{c,\Gamma_D}$, with $\vec{v}=\grad u$, and a $\vec{\alpha}\in H^1_{\Gamma_N}$, such that 
\begin{equation}
\label{hh6tris}
\vec{\xi}=\grad u + \curl \vec{\alpha}.
\end{equation}
Using \eqref{hh5} to express $\vec{\alpha}$ in \eqref{hh6tris} and the fact that $\curl \grad p=\vec{0}$ for any $p\in H^1(\mathcal{D}_a)$, we equivalently may express \eqref{hh6tris} with the requirement that $\vec{\alpha}\in H^1_{\Gamma_N,div}$, hence \eqref{hh6} and \eqref{hh6bis} are verified.

As a consequence of \eqref{phh1} and \eqref{hh5tris} we may construct the component $u$ in \eqref{hh6} as a solution of the elliptic problem with mixed boundary conditions
\begin{equation}
\label{phh6a}
\begin{cases}
\Delta u = \textrm{div}\vec{\xi},\\
u|_{\Gamma_D}=c,\\
\grad u \cdot \vec{n}|_{\Gamma_N}=\vec{\xi}\cdot \vec{n}|_{\Gamma_N},
\end{cases}
\end{equation}
Taking the curl of \eqref{hh6} and considering that $\text{div}\vec{\alpha}=0$, we may obtain $\vec{\alpha}$ as the solution of the elliptic problem
\begin{equation}
\label{phh6b}
\begin{cases}
-\Delta \vec{\alpha} = \text{curl}\vec{\xi},\;\; \diver \vec{\alpha}=0,\\
\text{curl} \,\vec{\alpha}\wedge \vec{n}|_{\Gamma_D}=\vec{\xi}\wedge \vec{n}|_{\Gamma_D},\\
\vec{\alpha}|_{\Gamma_N}=0.
\end{cases}
\end{equation}
We observe that setting $\vec{d}=\textrm{curl}\vec{\alpha}-\vec{\xi}$, we get from \eqref{phh6b} that $\textrm{curl}\vec{d}=\vec{0}$ and $\vec{d}\wedge \vec{n}|_{\Gamma_D}=\vec{0}$ in $\left(H_{00}^{1/2}(\Gamma_D)\right)'$, hence $\vec{d}\in H_{\Gamma_D,curl}\equiv \grad H^1_{c,\Gamma_D}$ and \eqref{hh6} is again obtained.
The mixed boundary conditions in \eqref{phh6a} and \eqref{phh6b} are then complementary with respect to the decomposition \eqref{hh6}.  
\end{pf}\\
\subsection{Functional inequalities}
We recall the Gagliardo-Nirenberg inequality {(see e.g. \cite{gagliardo,nirenberg,leoni})}.
\begin{lem}
\label{lem:gagliardoniremberg}
Let $\mathcal{D} \subset \mathbb{R}^3$ be a bounded domain with Lipschitz boundary and $f\in W^{m,r}\cap L^q$, $q\geq 1$, $r\leq \infty$, where $f$ can be a function with scalar, vectorial or tensorial values. For any integer $j$ with $0 \leq j < m$, suppose there is $\alpha \in \mathbb{R}$ such that
\[
j-\frac{3}{p}=\left(m-\frac{3}{r}\right)\alpha+(1-\alpha)\left(-\frac{3}{q}\right), \quad \frac{j}{m}\leq \alpha \leq 1.
\]
Then, there exists a positive constant $C$ depending on $\mathcal{D}$, m, j, q, r, and $\alpha$ such that
\begin{equation}
\label{eqn:2}
\pier{\|}D^jf\pier{\|}_{L^p}\leq C\pier{\|}f\pier{\|}_{W^{m,r}}^{\alpha}\pier{\|}f\pier{\|}_{L^q}^{1-\alpha}.
\end{equation}
\end{lem}
We also state the Agmon type inequality in three space dimensions {(see e.g. \cite{agmon})}.
\begin{lem}
\label{lem:agmon}
Let $\mathcal{D} \subset \mathbb{R}^3$ be a bounded domain with Lipschitz boundary and $f \in H^2(\Omega)$, where $f$ can be a function with scalar, vectorial or tensorial values. 
Then, there exists a positive constant $C$ depending on $\mathcal{D}$ such that
\begin{equation}
\label{eqn:agmon}
||f||_{L^{\infty}(\Omega)}\leq C||f||_{H^1(\Omega)}^{\frac{1}{2}}||f|||_{H^2(\Omega)}^{\frac{1}{2}}.
\end{equation}
\end{lem}
\section{Model derivation and main result}
\label{sec:model}
In this section we report the main steps of the model derivation which we introduced in \cite{agosti4}, generalizing the theoretical framework to consider the case of a {\mich visco-elastic} body which is fixed only on a part of its boundary, and is free to move on the other part.

Let $\Gamma_D,\Gamma_N \subseteq \Gamma_a$ be smooth subsets with positive measure of the domain boundary such that $\Gamma_a=\Gamma_D\cup\Gamma_N$, $|\Gamma_D\cap\Gamma_N|=0$.
We consider the motion of a deformable elastic solid in $\mathcal{D}_a$ which is fixed on $\Gamma_D$, which is immobile, while  no traction is applied to the remaining part of the boundary $\Gamma_N$ .
{\mich On this part $\Gamma_N$ the defect matrix $\mathbf{Z}$ remains constant and  keeps its initial value because the external forces being null do not modify the defects. On part $\Gamma_D$, an external action  can modify the defects. We assume there is no interaction with the exterior  related to the defects which are free to evolve.
}
In the time interval $(0,T)$, the motion is described by the displacement map
\[
(\vec{a},t)\rightarrow \vec{a}+\vec{u}(\vec{a},t)\in \mathbb{R}^3, \quad (\vec{a},t)\in \mathcal{D}_{aT},
\]
with initial condition
\[
\vec{u}(\vec{a},0)=\vec{0}\;\; \pier{\text{for}}\;\vec{a}\in \mathcal{D}_a
\]
and \abramon{Dirichlet boundary condition
\begin{align*}
&\vec{u}(\vec{a},t)=\vec{0}, \;\; \text{for} \; (\vec{a},t)\in \Gamma_D\times (0,T).
\end{align*}
}
We assume that the motion is not compatible, \pier{i.e.,} there exists a defect tensor $\mathbf{Z}\in \mathcal{M}_{\diver}$ with $\mathbf{Z}(\vec{a},0)=\mathbf{0}$ \pier{for} $\vec{a}\in \mathcal{D}_a$ and {$\mathbf{Z}(\vec{a},t)=\mathbf{0}$ for $(\vec{a},t)\in \Gamma_N\times (0,T)$}, such that
\begin{equation}
\label{eqn:3}
\grad{\vec{u}}=(\mathbf{R}\mathbf{W}-\mathbf{I})- \curl \mathbf{Z},
\end{equation}
where $\mathbf{R}\in \mathcal{SO}$ is the rotation tensor and $\mathbf{W}\in \mathcal{S}$ is the stretch tensor associated to the deformation gradient tensor, {with $\mathbf{R}\mathbf{W}(\vec{a},0)=\mathbf{I}$ for $\vec{a}\in \mathcal{D}_a$, \abramon{$\mathbf{R}(\vec{a},t)=\mathbf{W}(\vec{a},t)=\mathbf{I}$ 
{{}for $(\vec{a},t)\in\Gamma_D\times (0,T)$}.}
Since the $\grad, \diver, \curl$ operators are applied to second order tensors row-wise, we observe that the existence of the decomposition \eqref{eqn:3} is a consequence of the application of Theorem \ref{thm:hhd}, in particular of formula \eqref{hh6}, to the row vectors of the involved tensors. Given $\mathbf{R}$ and $\mathbf{W}$, the components $\vec{u}$ and $\mathbf{Z}$ in the decomposition \eqref{eqn:3} may be obtained as in \eqref{phh6}, i.e. solving elliptic problems with mixed boundary conditions derived by applying the divergence and the curl operators to \eqref{eqn:3}, ending with the kinematic relations:
\abramon{
\begin{equation}
\label{eqn:4}
\Delta \pier{\vec{u}}=\diver\left(\mathbf{R}\mathbf{W}\right),
\end{equation}
endowed with the boundary conditions $\pier{\vec{u}}(\vec{a},t)=\vec{0}$ {for $(\vec{a},t)\in\Gamma_D\times (0,T)$}, ${\grad \vec{u}}(\vec{a},t)\vec{n}=(\mathbf{R}\mathbf{W}-\mathbf{I})\vec{n}$ {{}for $(\vec{a},t)\in\Gamma_N\times (0,T)$},} and
\begin{equation}
\label{eqn:5}
- \Delta {\mathbf{Z}}=\curl \left(\mathbf{R}\mathbf{W}\right), \quad \diver{\mathbf{Z}}=\mathbf{0},
\end{equation}
endowed with the boundary condition ${\mathbf{Z}}(\vec{a},t)=\mathbf{0}$ for $(\vec{a},t)\in\Gamma_N\times (0,T)$, $\textrm{curl}\mathbf{Z}\wedge \vec{n}=\mathbf{0}$ \pco{{}for $(\vec{a},t)\in\Gamma_D\times (0,T)$}. Here, the notation $\textrm{curl}\mathbf{Z}\wedge \vec{n}$ represents the second order tensor whose rows are given by the vector product of the curl of a row vector of $\mathbf{Z}$ with $\vec{n}$, i.e. for a row index $i$, $(\textrm{curl}\mathbf{Z}\wedge \vec{n})_i=2Skew(\grad \mathbf{Z})_i\vec{n}$.
\abramon{\begin{rem}
The elliptic problem \eqref{eqn:4} has mixed Dirichlet--Neumann boundary conditions which are not standard, representing the situation of zero displacement on the Dirichlet boundary and {\mich displacement normal derivative in agreement with ${\mathbf{Z}}(\vec{a},t)=\mathbf{0}$ } on the Neumann boundary. A more standard mixed homogeneous Dirichlet--Neumann boundary condition, with $\pier{\vec{u}}(\vec{a},t)=\vec{0}$ {for $(\vec{a},t)\in\Gamma_D\times (0,T)$}, ${\grad \vec{u}}(\vec{a},t)\vec{n}=\vec{0}$ {{}for $(\vec{a},t)\in\Gamma_N\times (0,T)$}, would correspond to the situation in which $\mathbf{R}(\vec{a},t)\vec{n}=\vec{n}$ and $\mathbf{W}(\vec{a},t)\vec{n}=\vec{n}$ {{}for $(\vec{a},t)\in\Gamma_N\times (0,T)$}. The latter slip boundary conditions for the tensors $\mathbf{R}, \mathbf{W}$ over the Neumann boundary would represent the physical situation in which the solid is in contact with a device which forbids normal deformation, for
instance a rigid plate on which the structure slides. We highlight the fact that our theoretical framework could be easily adapted to consider the case with mixed homogeneous Dirichlet--Neumann boundary conditions for the displacement vector.
\end{rem}}

The model equations are derived from the principle of virtual powers. We make constitutive assumptions for the \michhhb{internal force} tensors, in terms of the \pierhhb{kinematic} variables, in order for \pier{the system} to satisfy the Clausius--Duhem dissipative equality. Given $\mathbf{R}\in \mathcal{SO}$, $\mathbf{W}\in \mathcal{S}$, we define, for any $t\in (0,T)$, the set $\mathcal{C}$ of virtual velocities as
\abramon{
\begin{align}
\label{eqn:8}
& \displaystyle \mathcal{C}:=\biggl\{\left(\vec{v},\nwhat{\mathbf{W}},\nwhat{\boldsymbol{\Omega}},\nwhat{\mathbf{Z}}\right)\in (\mathcal{V},\mathcal{S},\mathcal{A},\mathcal{M}_{\diver})\ \biggl| \;\; \nwhat{\mathbf{W}}|_{\Gamma_D}=\nwhat{\boldsymbol{\Omega}}|_{\Gamma_D}=\boldsymbol{0},  \;\; \grad\nwhat{\mathbf{W}}|_{\Gamma_D}=\grad\nwhat{\boldsymbol{\Omega}}|_{\Gamma_D}=\mathbf{0},\notag\\
& \displaystyle \qquad\quad
\begin{cases}
\Delta \vec{v}=\diver\left(\mathbf{R}\nwhat{\mathbf{W}}+\nwhat{\boldsymbol{\Omega}}\mathbf{R}\mathbf{W}\right),\\
\vec{v}=\vec{0} \quad \text{on} \; \Gamma_D\; ,\\
(\grad \vec{v})\vec{n}=\left(\mathbf{R}\nwhat{\mathbf{W}}+\nwhat{\boldsymbol{\Omega}}\mathbf{R}\mathbf{W}\right)\vec{n} \quad \text{on} \; \Gamma_N\; ,
\end{cases}
\begin{cases}
-P_L\Delta \nwhat{\mathbf{Z}}=\curl \left(\mathbf{R}\nwhat{\mathbf{W}}+\nwhat{\boldsymbol{\Omega}}\mathbf{R}\mathbf{W}\right),\\
\nwhat{\mathbf{Z}}=\mathbf{0} \quad \text{on} \; \Gamma_N \; ,\\
\textrm{curl} \nwhat{\mathbf{Z}}\wedge \vec{n}=\mathbf{0} \quad \text{on} \; \Gamma_D.
\end{cases}
\!\!\!\biggr\}
\end{align}
}
\abramonn{
\begin{rem}
In the following model derivation, we assume that the virtual velocities are sufficiently smooth to give a meaning to the formal variational formulation expressed by the principle of virtual powers. In the existence Theorem \ref{thm:1} the arguments will be made rigorous and we will explicitly state the regularity classes associated to the test functions of the weak formulation of the problem.
\end{rem}}
The virtual velocities then satisfy the following constraint, which, \abramohh{similarly to \eqref{eqn:3}, is a consequence of \eqref{hh6} applied row-wise:} 
\begin{equation}
\label{eqn:9}
\grad{\vec{v}}=\mathbf{R}\nwhat{\mathbf{W}}+\nwhat{\boldsymbol{\Omega}}\mathbf{R}\mathbf{W}-\curl \nwhat{\mathbf{Z}}.
\end{equation}
We observe that the set $\mathcal{C}$ of virtual velocities is defined in terms of the variables $\mathbf{R}$ and $\mathbf{W}$, and hence depend on the solutions of the \michhhb{equations of motion}. We can formally write
\abramon{
\begin{equation}
\label{eqn:10}
\vec{v}=-\mathcal{G}_L\ast\diver\left(\mathbf{R}\nwhat{\mathbf{W}}+\nwhat{\boldsymbol{\Omega}}\mathbf{R}\mathbf{W}\right)+\mathcal{G}_L\ast_{|_{\Gamma_N}}\left(\mathbf{R}\nwhat{\mathbf{W}}+\nwhat{\boldsymbol{\Omega}}\mathbf{R}\mathbf{W}\right)\vec{n},
\end{equation}
where we used the notation $\mathcal{G}_L\ast_{|_{\Gamma_N}}\vec{f}:=\int_{\Gamma_N}\mathcal{G}_L(\vec{x}-\vec{y})\vec{f}(\vec{y})dS(\vec{y})$,
}
and
\begin{equation}
\label{eqn:11}
\nwhat{\mathbf{Z}}=\mathcal{G}_{L,\diver}\ast\curl\left(\mathbf{R}\nwhat{\mathbf{W}}+\nwhat{\boldsymbol{\Omega}}\mathbf{R}\mathbf{W}\right).
\end{equation}

The principle of virtual powers takes the form
\begin{equation}
\label{eqn:12}
p_{\text{acc}}(\mathcal{D}_a,C)=p_{\text{int}}(\mathcal{D}_a,C)+p_{\text{ext}}(\mathcal{D}_a,C) \quad \forall C\in \mathcal{C},
\end{equation}
where $p_{\text{int}}(\mathcal{D}_a,C)$ is the virtual power of internal forces, $p_{\text{ext}}(\mathcal{D}_a,C)$ is the virtual power of external forces and $p_{\text{acc}}(\mathcal{D}_a,C)$ is the virtual power of acceleration forces, defined in terms of $\mathcal{D}_a$ and of an element $C\in \mathcal{C}$.
\noindent
The virtual power of internal forces is defined as
\begin{align*}
& \displaystyle p_{\text{int}}(\mathcal{D}_a,C):=-\int_{\mathcal{D}_a}\left(\boldsymbol{\Pi}:\grad \vec{v}+\mathbf{X}:: \grad \nwhat{\mathbf{W}}+\mathbf{Y}::: \grad\grad \nwhat{\mathbf{W}}\right) \\
& \displaystyle +\frac{1}{2}\int_{\mathcal{D}_a}\left(\mathbf{M}:\nwhat{\boldsymbol{\Omega}}-\boldsymbol{\Lambda}::\grad \nwhat{\boldsymbol{\Omega}}-\mathbf{C}::: \grad\grad \nwhat{\boldsymbol{\Omega}}\right)+\int_{\mathcal{D}_a} \ap{\boldsymbol{\Gamma}:\curl \nwhat{\mathbf{Z}}},
\end{align*}
where $\boldsymbol{\Pi}$ is the Piola--Kirchhoff--Boussinesq stress tensor, $\mathbf{M}$ represents the momentum, $\boldsymbol{\Lambda}$ the momentum flux and $\mathbf{C}$ the flux of the momentum flux. The quantities $\mathbf{X}, \mathbf{Y}, \boldsymbol{\Gamma}$ are new \michhhb{internal force} tensors associated to the kinematic variables $\mathbf{W}$ and $\mathbf{Z}$. \abramohh{In particular, $\boldsymbol{\Gamma}$ is an internal force accounting for the evolution
of the defects.} We impose the following boundary conditions for the internal forces:
\begin{equation}
\label{pigamma}
(\boldsymbol{\Pi}+\boldsymbol{\Gamma})\wedge \vec{n}=\mathbf{0} \quad \text{for}\;\; (\vec{a},t)\in\Gamma_D\times (0,T).
\end{equation}
Using \eqref{eqn:9}, we rewrite the virtual power of internal forces as
\begin{align}
\label{eqn:13a}
& \displaystyle p_{\text{int}}(\mathcal{D}_a,C):=-\int_{\mathcal{D}_a}\left(\mathbf{R}^T\boldsymbol{\Pi}:\nwhat{\mathbf{W}}+\mathbf{X}:: \grad \nwhat{\mathbf{W}}+\mathbf{Y}:::  \grad \grad \nwhat{\mathbf{W}}\right) \notag \\
& \displaystyle \notag +\frac{1}{2}\int_{\mathcal{D}_a}\left((\mathbf{M}-2\boldsymbol{\Pi}\mathbf{W}\mathbf{R}^T):\nwhat{\boldsymbol{\Omega}}-\boldsymbol{\Lambda}::\grad \nwhat{\boldsymbol{\Omega}}-\mathbf{C}::: \grad \grad \nwhat{\boldsymbol{\Omega}}\right) \notag
\\ 
&+\displaystyle\int_{\mathcal{D}_a}\left(\boldsymbol{\Gamma}+\boldsymbol{\Pi}\right):\curl  \nwhat{\mathbf{Z}},
\end{align}
We integrate by parts the last term in the previous equation, using \eqref{eqn:11} and employing the boundary conditions on $\nwhat{\mathbf{Z}},\nwhat{\mathbf{W}},\nwhat{\boldsymbol{\Omega}}$ and \eqref{pigamma}, obtaining that
\begin{align*}
& \displaystyle \int_{\mathcal{D}_a}\left(\boldsymbol{\Gamma}+\boldsymbol{\Pi}\right):\curl  \nwhat{\mathbf{Z}}= \int_{\mathcal{D}_a}\curl \left(\boldsymbol{\Gamma}+\boldsymbol{\Pi}\right):\nwhat{\mathbf{Z}}-\underbrace{\int_{\Gamma_N}\left(\boldsymbol{\Gamma}+\boldsymbol{\Pi}\right):\nwhat{\mathbf{Z}}\wedge \vec{n}}_{=0}\\
& \displaystyle +\underbrace{\int_{\Gamma_D}\left(\boldsymbol{\Gamma}+\boldsymbol{\Pi}\right)\wedge \vec{n}:\nwhat{\mathbf{Z}}}_{=0}=\int_{\mathcal{D}_a}\curl \left(\boldsymbol{\Gamma}+\boldsymbol{\Pi}\right):\mathcal{G}_{L,\diver}\ast\curl\left(\mathbf{R}\nwhat{\mathbf{W}}+\nwhat{\boldsymbol{\Omega}}\mathbf{R}\mathbf{W}\right)\\
& \displaystyle =\int_{\mathcal{D}_a}\curl\left(\mathcal{G}_{L,\diver}\ast \curl(\boldsymbol{\Gamma}+\boldsymbol{\Pi})\right):\left(\mathbf{R}\nwhat{\mathbf{W}}+\nwhat{\boldsymbol{\Omega}}\mathbf{R}\mathbf{W}\right)\\
& \displaystyle - \underbrace{\int_{\Gamma_D}\mathcal{G}_{L,\diver}\ast \curl(\boldsymbol{\Gamma}+\boldsymbol{\Pi}):\left(\mathbf{R}\nwhat{\mathbf{W}}+\nwhat{\boldsymbol{\Omega}}\mathbf{R}\mathbf{W}\right)\wedge \vec{n}}_{=0}\\
& \displaystyle +\underbrace{\int_{\Gamma_N}\left(\mathcal{G}_{L,\diver}\ast \curl(\boldsymbol{\Gamma}+\boldsymbol{\Pi})\right)\wedge \vec{n}:\left(\mathbf{R}\nwhat{\mathbf{W}}+\nwhat{\boldsymbol{\Omega}}\mathbf{R}\mathbf{W}\right)}_{=0}.
\end{align*} 
Hence, equation \eqref{eqn:13a} becomes
\begin{align}
\label{eqn:13}
& \displaystyle p_{\text{int}}(\mathcal{D}_a,C):=-\int_{\mathcal{D}_a}\left(\mathbf{R}^T\boldsymbol{\Pi}:\nwhat{\mathbf{W}}+\mathbf{X}:: \grad \nwhat{\mathbf{W}}+\mathbf{Y}:::  \grad\grad \nwhat{\mathbf{W}}\right) \notag \\
& \displaystyle \notag +\frac{1}{2}\int_{\mathcal{D}_a}\left((\mathbf{M}-2\boldsymbol{\Pi}\mathbf{W}\mathbf{R}^T):\nwhat{\boldsymbol{\Omega}}-\boldsymbol{\Lambda}::\grad \nwhat{\boldsymbol{\Omega}}-\mathbf{C}::: \grad\grad \nwhat{\boldsymbol{\Omega}}\right) \notag
\\ 
&+\displaystyle\int_{\mathcal{D}_a}\curl\left(\mathcal{G}_{L,\diver}\ast \curl(\boldsymbol{\Gamma}+\boldsymbol{\Pi})\right):\left(\mathbf{R}\nwhat{\mathbf{W}}+\nwhat{\boldsymbol{\Omega}}\mathbf{R}\mathbf{W}\right),
\end{align}
 The virtual power of external forces is defined as
\begin{equation}
\label{eqn:14}
p_{\text{ext}}(\mathcal{D}_a,C):=\int_{\mathcal{D}_a}\vec{\mathcal{F}}_{\text{ext}}\cdot \vec{v}+\int_{\mathcal{D}_a}\boldsymbol{\mathcal{W}}_{\text{ext}}:\nwhat{\mathbf{W}}+ \int_{\mathcal{D}_a}\boldsymbol{\Omega}_{\text{ext}}:\nwhat{\boldsymbol{\Omega}},
\end{equation}
where $\boldsymbol{\mathcal{W}}_{\text{ext}}$ and $\boldsymbol{\Omega}_{\text{ext}}$ are external forces, \abramohhb{possibly depending on $\mathbf{W}$ and $\boldsymbol{R}$}, which perform work by stretching and rotating the system, respectively, while $\vec{\mathcal{F}}_{\text{ext}}$ is an external volume force depending only on time and position, which may account e.g. for gravitation. 
Finally, the virtual power of acceleration forces is defined as 
\begin{align}
\label{eqn:15}
& \displaystyle \notag
p_{\text{acc}}(\mathcal{D}_a,C):=\int_{\mathcal{D}_a}\twobigdot{\vec{u}}\cdot \vec{v}+\int_{\mathcal{D}_a} 
\left( \twobigdot{\mathbf{W}}:
\nwhat{\mathbf{W}}+\grad \twobigdot{\mathbf{W}}::
\grad \nwhat{\mathbf{W}}+\grad\grad\twobigdot{\mathbf{W}}:::\grad\grad\nwhat{\mathbf{W}}\right)\\
& \displaystyle +\int_{\mathcal{D}_a} \left( \overbigdot{\boldsymbol{\Omega}}:\nwhat{\boldsymbol{\Omega}}+\grad \overbigdot{\boldsymbol{\Omega}}::\grad \nwhat{\boldsymbol{\Omega}}+
\grad\grad\overbigdot{\boldsymbol{\Omega}}:::\grad\grad\nwhat{\boldsymbol{\Omega}}\right),
\end{align}
\abramonn{where
\begin{equation}
\label{eqn:15r}
\overbigdot{\mathbf{R}}=\boldsymbol{\Omega}\mathbf{R}.
\end{equation}}
\abramohh{As discussed in \cite{agosti4}, higher order terms in the virtual power of acceleration forces are introduced to be able to represent a situation of a contact at a point with inertia, which requires regularity in space and time of the angular velocity and acceleration variables.}
The first term on the right hand side of \eqref{eqn:15}, and similarly for the first term on the right hand side of \eqref{eqn:14}, can be rewritten using \eqref{eqn:10}, integrating by parts and employing the boundary conditions for $\nwhat{\mathbf{W}},\nwhat{\boldsymbol{\Omega}}$, obtaining that
\abramon{
\begin{align}
& \displaystyle  \notag
\int_{\mathcal{D}_a}\twobigdot{\vec{u}}\cdot \vec{v}=-\int_{\mathcal{D}_a}\twobigdot{\vec{u}}\cdot \mathcal{G}_L\ast\diver\left(\mathbf{R}\nwhat{\mathbf{W}}+\nwhat{\boldsymbol{\Omega}}\mathbf{R}\mathbf{W}\right)+\int_{\mathcal{D}_a}\twobigdot{\vec{u}}\cdot\mathcal{G}_L\ast_{|_{\Gamma_N}}\left(\mathbf{R}\nwhat{\mathbf{W}}+\nwhat{\boldsymbol{\Omega}}\mathbf{R}\mathbf{W}\right)\vec{n}\\
& \displaystyle \notag \quad =\int_{\mathcal{D}_a}\grad \left(\mathcal{G}_L\ast \twobigdot{\vec{u}}\right): \left(\mathbf{R}\nwhat{\mathbf{W}}+\nwhat{\boldsymbol{\Omega}}\mathbf{R}\mathbf{W}\right)-\int_{\Gamma_a}\mathcal{G}_L\ast \twobigdot{\vec{u}}: \left(\mathbf{R}\nwhat{\mathbf{W}}+\nwhat{\boldsymbol{\Omega}}\mathbf{R}\mathbf{W}\right)\vec{n}\\
& \displaystyle \notag \qquad +\int_{\mathcal{D}_a}\twobigdot{\vec{u}}\cdot\mathcal{G}_L\ast_{|_{\Gamma_N}}\left(\mathbf{R}\nwhat{\mathbf{W}}+\nwhat{\boldsymbol{\Omega}}\mathbf{R}\mathbf{W}\right)\vec{n}=\int_{\mathcal{D}_a}\grad \left(\mathcal{G}_L\ast \twobigdot{\vec{u}}\right): \left(\mathbf{R}\nwhat{\mathbf{W}}+\nwhat{\boldsymbol{\Omega}}\mathbf{R}\mathbf{W}\right)\\
& \displaystyle \notag \qquad -\underbrace{\int_{\Gamma_D}\mathcal{G}_L\ast \twobigdot{\vec{u}}: \left(\mathbf{R}\nwhat{\mathbf{W}}+\nwhat{\boldsymbol{\Omega}}\mathbf{R}\mathbf{W}\right)\vec{n}}_{=0}-\int_{\Gamma_N}\mathcal{G}_L\ast \twobigdot{\vec{u}}: \left(\mathbf{R}\nwhat{\mathbf{W}}+\nwhat{\boldsymbol{\Omega}}\mathbf{R}\mathbf{W}\right)\vec{n}\\
& \displaystyle \notag \qquad +\int_{\mathcal{D}_a}\twobigdot{\vec{u}}\cdot\mathcal{G}_L\ast_{|_{\Gamma_N}}\left(\mathbf{R}\nwhat{\mathbf{W}}+\nwhat{\boldsymbol{\Omega}}\mathbf{R}\mathbf{W}\right)\vec{n}=\int_{\mathcal{D}_a}\grad \left(\mathcal{G}_L\ast \twobigdot{\vec{u}}\right): \left(\mathbf{R}\nwhat{\mathbf{W}}+\nwhat{\boldsymbol{\Omega}}\mathbf{R}\mathbf{W}\right)\\
& \displaystyle \notag \qquad -\int_{\mathcal{D}_a}\twobigdot{\vec{u}}\cdot\mathcal{G}_L\ast_{|_{\Gamma_N}}\left(\mathbf{R}\nwhat{\mathbf{W}}+\nwhat{\boldsymbol{\Omega}}\mathbf{R}\mathbf{W}\right)\vec{n}+\int_{\mathcal{D}_a}\twobigdot{\vec{u}}\cdot\mathcal{G}_L\ast_{|_{\Gamma_N}}\left(\mathbf{R}\nwhat{\mathbf{W}}+\nwhat{\boldsymbol{\Omega}}\mathbf{R}\mathbf{W}\right)\vec{n}\\
& \displaystyle \quad \label{eqn:15bis} =
\int_{\mathcal{D}_a}\grad \left(\mathcal{G}_L\ast \twobigdot{\vec{u}}\right): \left(\mathbf{R}\nwhat{\mathbf{W}}+\nwhat{\boldsymbol{\Omega}}\mathbf{R}\mathbf{W}\right),\\
& \displaystyle \notag
\int_{\mathcal{D}_a}\vec{\mathcal{F}}_{\text{ext}}\cdot \vec{v}=-\int_{\mathcal{D}_a}\vec{\mathcal{F}}_{\text{ext}}\cdot \mathcal{G}_L\ast\diver\left(\mathbf{R}\nwhat{\mathbf{W}}+\nwhat{\boldsymbol{\Omega}}\mathbf{R}\mathbf{W}\right)+\int_{\mathcal{D}_a}\vec{\mathcal{F}}_{\text{ext}}\cdot\mathcal{G}_L\ast_{|_{\Gamma_N}}\left(\mathbf{R}\nwhat{\mathbf{W}}+\nwhat{\boldsymbol{\Omega}}\mathbf{R}\mathbf{W}\right)\vec{n}\\
& \displaystyle \quad \label{eqn:15tris} =\int_{\mathcal{D}_a}\grad \left(\mathcal{G}_L\ast \vec{\mathcal{F}}_{\text{ext}}\right): \left(\mathbf{R}\nwhat{\mathbf{W}}+\nwhat{\boldsymbol{\Omega}}\mathbf{R}\mathbf{W}\right).
\end{align}
}
\begin{rem}
 \label{rem:stress}
 We observe that \eqref{eqn:15bis} and \eqref{eqn:15tris} \abramon{are such that the virtual power of
the acceleration forces and the virtual power of the external forces depends only
on virtual velocities $\nwhat{\mathbf{W}}$ and $\nwhat{\boldsymbol{\Omega}}$. It results the principle of virtual power \eqref{eqn:12} is
no longer given in terms of $\mathcal{C}$ but only in terms of $(\nwhat{\mathbf{W}},\nwhat{\boldsymbol{\Omega}})$.}
\end{rem}
\abramonn{
Finally, the principle of virtual powers \eqref{eqn:12} takes the following form:
\begin{equation}
\label{eqn:20w1}
\begin{cases}
\displaystyle \int_{\mathcal{D}_a}\left(\mathbf{R}^T\grad \left(\mathcal{G}_L\ast\left( \twobigdot{\vec{u}}-\vec{\mathcal{F}}_{\textrm ext}\right)\right)+\twobigdot{\mathbf{W}}\right):\nwhat{\mathbf{W}}+\int_{\mathcal{D}_a}\grad \twobigdot{\mathbf{W}}::\grad \nwhat{\mathbf{W}}\\
\displaystyle \;\; +\int_{\mathcal{D}_a}\grad \grad  \twobigdot{\mathbf{W}}::: \grad \grad  \nwhat{\mathbf{W}}-\int_{\mathcal{D}_a}\left(\mathbf{R}^T\curl\left(\mathcal{G}_{L,\diver}\ast \curl(\boldsymbol{\Gamma}+\boldsymbol{\Pi})\right)+\mathbf{R}^T\boldsymbol{\Pi}\right):\nwhat{\mathbf{W}}\\
\displaystyle \;\; +\int_{\mathcal{D}_a}\mathbf{X}::\grad \nwhat{\mathbf{W}}+\int_{\mathcal{D}_a} \mathbf{Y}:::\grad \grad \nwhat{\mathbf{W}}=\int_{\mathcal{D}_a}\mathbf{W}_{\text{ext}}:\nwhat{\mathbf{W}},\\ \\
\displaystyle \int_{\mathcal{D}_a}\left(\grad \left(\mathcal{G}_L\ast \left(\twobigdot{\vec{u}}-\vec{\mathcal{F}}_{\textrm ext}\right)\right)\mathbf{W}\mathbf{R}^T+\overbigdot{\boldsymbol{\Omega}}\right):\nwhat{\boldsymbol{\Omega}}+\int_{\mathcal{D}_a}\grad \overbigdot{\boldsymbol{\Omega}}:: \grad \nwhat{\boldsymbol{\Omega}}\\
\displaystyle \;\; +\int_{\mathcal{D}_a}\grad \grad \overbigdot{\boldsymbol{\Omega}}::: \grad \grad \nwhat{\boldsymbol{\Omega}} -\int_{\mathcal{D}_a}\curl\left(\mathcal{G}_{L,\diver}\ast \curl(\boldsymbol{\Gamma}+\boldsymbol{\Pi})\right)\mathbf{W}\mathbf{R}^T: \nwhat{\boldsymbol{\Omega}}\\
\displaystyle \;\; -\int_{\mathcal{D}_a} \frac{1}{2}(\mathbf{M}-2\boldsymbol{\Pi}\mathbf{W}\mathbf{R}^T):\nwhat{\boldsymbol{\Omega}}+\frac{1}{2}\int_{\mathcal{D}_a} \boldsymbol{\Lambda}:: \grad \nwhat{\boldsymbol{\Omega}}+\frac{1}{2}\int_{\mathcal{D}_a} \mathbf{C}::: \grad \grad \nwhat{\boldsymbol{\Omega}}=\int_{\mathcal{D}_a}\boldsymbol{\Omega}_{\text{ext}}:\nwhat{\boldsymbol{\Omega}},\\ \\
\end{cases}
\end{equation}
valid for all choices of $\nwhat{\mathbf{W}}\in \mathcal{S}$, $\nwhat{\boldsymbol{\Omega}} \in \mathcal{A}$, with $\nwhat{\mathbf{W}}|_{\Gamma_D}=\nwhat{\boldsymbol{\Omega}}|_{\Gamma_D}=\boldsymbol{0}$, $\grad\nwhat{\mathbf{W}}|_{\Gamma_D}=\grad\nwhat{\boldsymbol{\Omega}}|_{\Gamma_D}=\mathbf{0}$, with boundary conditions
\begin{align}
\label{eqn:20w1bc}
& \notag \displaystyle \mathbf{W}= \mathbf{I},\; \grad {\mathbf{W}}=\mathbf{0}\quad \text{on} \; \Gamma_D\times (0,T),\\
& \displaystyle \boldsymbol{\Omega}= \mathbf{0},\; \grad {\boldsymbol{\Omega}}=\mathbf{0}\quad \text{on} \; \Gamma_D\times (0,T).
\end{align}
System \eqref{eqn:20w1} must be coupled with the kinematic relations \eqref{eqn:4} and \eqref{eqn:5}, expressed in a variational form as
\begin{equation}
\label{eqn:20w2}
\begin{cases}
\displaystyle \int_{\mathcal{D}_a}\grad{}{{\vec{u}}}\colon \grad{}\vec{v}=\int_{\mathcal{D}_a}\left(\mathbf{R}\mathbf{W}-\mathbf{I}\right)\colon \grad{}\vec{v},\\ \\
\displaystyle \int_{\mathcal{D}_a}\grad {\mathbf{Z}}:: \grad \nwhat{\mathbf{Z}}=\int_{\mathcal{D}_a}\left(\mathbf{R}\mathbf{W}-\mathbf{I}\right)\colon \curl \nwhat{\mathbf{Z}},
\end{cases}
\end{equation}
valid for {all choices of} $\vec{v}\in \mathcal{V}$ and $\nwhat{\mathbf{Z}}\in \mathcal{M}_{div}$, with $\vec{v}|_{\Gamma_D}=\vec{0}$, $\nwhat{\mathbf{Z}}|_{\Gamma_N}=\mathbf{0}$, with boundary conditions
\begin{align}
\label{eqn:20w2bc}
& \notag \displaystyle \vec{u}=\vec{0} \quad \text{on} \; \Gamma_D\times (0,T),\\
& \displaystyle \mathbf{Z}= \mathbf{0} \quad \text{on} \; \Gamma_N\times (0,T).
\end{align}
Systems \eqref{eqn:20w1}-\eqref{eqn:20w2} is endowed with the initial conditions
\begin{equation}
\label{eqn:20wic}
\mathbf{W}(\vec{a},0)=\mathbf{I}, \; \overbigdot{\mathbf{W}}(\vec{a},0)=\mathbf{0}, \; \mathbf{R}(\vec{a},0)=\mathbf{I}, \; \boldsymbol{\Omega}(\vec{a},0)=\mathbf{0}, \; \mathbf{Z}(\vec{a},0)=\mathbf{0}
\quad {\text{for} \; \vec{a} \in \mathcal{D}_a}.
\end{equation} 
\begin{rem}
We observe that in \eqref{eqn:20w2}$_1$ the Neumann boundary conditions $$(\grad \vec{u})\vec{n}=(\mathbf{RW-I})\vec{n}\; \; \text{on}\; \; \Gamma_N\times (0,T)$$ are enforced weakly within the variational formulation, obtained formally from testing \eqref{eqn:4} with a test function $\vec{v}$ and integrating by parts. 
\end{rem}
Starting from the variational formulations \eqref{eqn:20w1} and \eqref{eqn:20w2} and integrating by parts in the first and second gradient terms, it is possible, upon assigning proper Neumann boundary conditions for the {internal forces} and the acceleration forces, to derive a strong form of the principle of virtual powers (see e.g. \cite{germain2}). In order to proceed in this sense, we report the following integration by parts formula, involving generic third order tensor $\mathbf{A}$ and fourth order tensor $\mathbf{B}$:
\begin{align}
\label{bcparts}
& \displaystyle \notag \int_{\mathcal{D}_a}\left(\mathbf{A}:: \grad \nwhat{\mathbf{W}}+\mathbf{B}:::  \grad\grad \nwhat{\mathbf{W}}\right)=\int_{\Gamma_a}\left[\mathbf{A}-(\diver \mathbf{B})\right]\vec{n}:\nwhat{\mathbf{W}}+\int_{\Gamma_a}\mathbf{B}\vec{n}::\grad \nwhat{\mathbf{W}}\\
& \quad \displaystyle +\int_{\mathcal{D}_a}\left(\diver\diver \mathbf{B}-\diver \mathbf{A}\right):\nwhat{\mathbf{W}}.
\end{align}
}\abramonn{
The quantities $\nwhat{\mathbf{W}}$ and $\grad \nwhat{\mathbf{W}}$ in the right hand side of \eqref{bcparts} are not independent on the surface. To get a
relationship with independent virtual velocities, recalling the following identity (see e.g. \cite[Appendix]{germain2}), valid for a given tensor field $\mathbf{T}$ of any order and for a sufficiently smooth surface $S$,
\[
\int_{S}\diver \mathbf{T}=\int_{S}\left(2\kappa \mathbf{T}\vec{n}+[(\grad \mathbf{T})\vec{n}]\vec{n}\right),
\]
where $\kappa$ is the mean curvature of the surface, we may write
\begin{align}
\label{byparts1}
& \displaystyle \notag \int_{\Gamma_a}\left[\mathbf{A}-(\diver \mathbf{B})\right]\vec{n}:\nwhat{\mathbf{W}}+\int_{\Gamma_a}\mathbf{B}\vec{n}::\grad \nwhat{\mathbf{W}}=\int_{\Gamma_a}\left[\mathbf{A}-(\diver \mathbf{B})\right]\vec{n}:\nwhat{\mathbf{W}}\\
& \displaystyle  \quad +\int_{\Gamma_a}\diver\left(\mathbf{B}\vec{n}\odot \nwhat{\mathbf{W}}\right)-\int_{\Gamma_a}\diver(\mathbf{B}\vec{n}):\nwhat{\mathbf{W}} = \int_{\Gamma_N}W_{N}(\mathbf{A},\mathbf{B}):\nwhat{\mathbf{W}}+\int_{\Gamma_N}K_{N}(\mathbf{B}):\left(\grad \nwhat{\mathbf{W}}\right)\vec{n},
\end{align}
where in the last step we employed the boundary conditions on $\nwhat{\mathbf{W}}$, with $\left(\mathbf{B}\vec{n}\odot \nwhat{\mathbf{W}}\right)_k:=\sum_{i,j,l=1}^3B_{ijkl}n_lW_{ij}$, for $k=1,\dots,3$, and
\begin{equation}
\label{byparts2}
\displaystyle W_{N}(\mathbf{A},\mathbf{B}):=2\kappa (\mathbf{B}\vec{n})\vec{n}+[[(\grad(\mathbf{B}\vec{n}))]\vec{n}]\vec{n}-\diver(\mathbf{B}\vec{n})+\left[\mathbf{A}-(\diver \mathbf{B})\right]\vec{n}, \quad K_{N}(\mathbf{B}):=(\mathbf{B}\vec{n})\vec{n}.
\end{equation}
\begin{rem}
\label{parts2}
We observe that, in the case with regularity $\nwhat{\mathbf{W}}\in \bar{H}_{\Gamma_D}^2\left(\mathcal{D}_a,Sym\left(\mathbb{R}^{3\times 3}\right)\right)$, $\mathbf{A}\in L^2(\mathcal{D}_a,\mathbb{R}^{3\times 3\times 3})$, $\mathbf{B}\in L^2(\mathcal{D}_a,\mathbb{R}^{3\times 3\times 3\times 3})$, the boundary terms in \eqref{byparts2} should be interpreted as proper dual products and the corresponding Neumann boundary conditions would be valid in the trace spaces $K_N(\mathbf{B})\in \left(H^{1/2}_{00}(\Gamma_N,\mathbb{R}^{3\times 3\times 3})\right)'$, $W_{N}(\mathbf{A},\mathbf{B})\in \left(H^{3/2}_{00}(\Gamma_N,\mathbb{R}^{3\times 3})\right)'$.
\end{rem}
Given the integration by parts formula \eqref{byparts1}, the principle of virtual powers \eqref{eqn:20w1} results in a linear function of $\nwhat{\mathbf{W}}$ and $\left(\grad \nwhat{\mathbf{W}}\right)\vec{n}$ on $\Gamma_N$ to be null for any $\nwhat{\mathbf{W}}$ and $\left(\grad \nwhat{\mathbf{W}}\right)\vec{n}$ on $\Gamma_N$, which gives two surface equations to be satisfied as Neumann boundary conditions over $\Gamma_N$.
\newline
Then, the principle of virtual powers \eqref{eqn:12} implies the following {equations}, valid in $\mathcal{D}_{aT}$, which are coupled to the kinematic relations \eqref{eqn:4} and \eqref{eqn:5}:
\begin{equation}
\label{eqn:20}
\begin{cases}
\displaystyle \text{Sym}\left(\mathbf{R}^T\grad \left(\mathcal{G}_L\ast\left( \twobigdot{\vec{u}}-\vec{\mathcal{F}}_{\textrm ext}\right)\right)\right)+\twobigdot{\mathbf{W}}-\Delta \twobigdot{\mathbf{W}}+{{}\diver \Delta \grad{}} \twobigdot{\mathbf{W}}\\
\displaystyle -\text{Sym}\left(\mathbf{R}^T\curl\left(\mathcal{G}_{L,\diver}\ast \curl(\boldsymbol{\Gamma}+\boldsymbol{\Pi})\right)+\mathbf{R}^T\boldsymbol{\Pi}\right)-\diver \mathbf{X}+\diver \diver \mathbf{Y}=\mathbf{W}_{\text{ext}},\\ \\
\displaystyle \mathbf{W}= \mathbf{I},\; \grad {\mathbf{W}}=\mathbf{0}\quad \text{on} \; \Gamma_D\times (0,T),\\ \\
\displaystyle W_N(\mathbf{X},\mathbf{Y})+W_N(\grad \twobigdot{\mathbf{W}},\grad\grad \twobigdot{\mathbf{W}})=\mathbf{0}, \;\; K_N(\mathbf{Y})+K_N(\grad\grad \twobigdot{\mathbf{W}})=\mathbf{0} \;\; \text{on} \; \Gamma_N\times (0,T),\\ \\
\displaystyle \text{Skew}\left(\grad \left(\mathcal{G}_L\ast \left(\twobigdot{\vec{u}}-\vec{\mathcal{F}}_{\textrm ext}\right)\right)\mathbf{W}\mathbf{R}^T\right)+\overbigdot{\boldsymbol{\Omega}}-\Delta \overbigdot{\boldsymbol{\Omega}}+{{}{\diver \Delta \grad}{}}\overbigdot{\boldsymbol{\Omega}}\\
\displaystyle -\text{Skew}\left(\curl\left(\mathcal{G}_{L,\diver}\ast \curl(\boldsymbol{\Gamma}+\boldsymbol{\Pi})\right)\mathbf{W}\mathbf{R}^T -\frac{1}{2}(\mathbf{M}-2\boldsymbol{\Pi}\mathbf{W}\mathbf{R}^T)\right)-\frac{1}{2}\diver \boldsymbol{\Lambda}+\frac{1}{2}{\diver\diver} \mathbf{C}=\boldsymbol{\Omega}_{\text{ext}},\\ \\
\displaystyle \boldsymbol{\Omega}= \mathbf{0},\; \grad {\boldsymbol{\Omega}}=\mathbf{0}\quad \text{on} \; \Gamma_D\times (0,T),\\ \\
\displaystyle W_N\left(\frac{\boldsymbol{\Lambda}}{2},\frac{\mathbf{C}}{2}\right)+W_N(\grad \overbigdot{\boldsymbol{\Omega}},\grad\grad \overbigdot{\boldsymbol{\Omega}})=\mathbf{0}, \;\; K_N\left(\frac{\mathbf{C}}{2}\right)+K_N(\grad\grad \overbigdot{\boldsymbol{\Omega}})=\mathbf{0} \;\; \text{on} \; \Gamma_N\times (0,T),\\ \\
\overbigdot{\mathbf{R}}=\boldsymbol{\Omega}\mathbf{R},\\ \\
\Delta {\vec{u}}=\diver\left(\mathbf{R}\mathbf{W}\right), \quad  {\vec{u}=\mathbf{0} \quad \text{on} \; \Gamma_D\times (0,T)}, \quad {(\grad \vec{u})\vec{n}=(\mathbf{RW-I})\vec{n}}\quad \text{on} \; \Gamma_N\times (0,T),
\\ \\
- P_L\Delta {\mathbf{Z}}=\curl \left(\mathbf{R}\mathbf{W}\right),\quad {\mathbf{Z}=\mathbf{0} \quad \text{on} \; \Gamma_N\times (0,T)}, \quad \curl \mathbf{Z}\wedge \vec{n}=\mathbf{0} \; \text{on} \; \Gamma_D\times (0,T).
\end{cases}
\end{equation}
\begin{rem}
We observe that the PDE system \eqref{eqn:20}, which represents the strong form of the variational equations \eqref{eqn:20w1}-\eqref{eqn:20w2} associated to the principle of virtual powers \eqref{eqn:12}, is complemented by an involved set of mixed boundary conditions. The boundary conditions required in the variational formulation of the problem, i.e. \eqref{eqn:20w1bc}-\eqref{eqn:20w2bc}, are more simple and specified in the definition of the virtual velocities \eqref{eqn:8}. In reality the dynamics is described by the principle of virtual powers in its variational formulation and by dissipative {\mich and non dissipative} laws, where {\mich the forces, even intricate and sophistigated,  may be experimented with their powers.} For this reason, in the following we will work with the variational forms \eqref{eqn:20w1}-\eqref{eqn:20w2}, which will be directly linked to the weak formulation of the problem expressed in Theorem \ref{thm:1}.
\end{rem}
}
We now assign general constitutive assumptions for $\boldsymbol{\Pi},\mathbf{M},\mathbf{X},\mathbf{Y},\boldsymbol{\Lambda},\mathbf{C},\boldsymbol{\Gamma}$ in order for \eqref{eqn:20} to satisfy the Clausius--Duhem dissipative equality in isothermal situations, which has the form
\abramo{
\begin{equation}
\label{eqn:21}
\frac{d\psi}{dt}+\left(\frac{d D}{d\overbigdot{C}}(\overbigdot{C}),\overbigdot{C}\right)=-p_{\text{int}}(\mathcal{D}_a,\overbigdot{C}),
\end{equation}
}%
where $\overbigdot{C}:=(\overbigdot{\mathbf{W}},\boldsymbol{\Omega})$ is the actual velocity, $\psi$ is the free energy of the system and $D$ is the dissipation potential. We assume the following form for the free energy of the system:
\begin{align}
\label{eqn:22}
\displaystyle \notag &\psi(\mathbf{W},\mathbf{R},\mathbf{Z}):=\frac{1}{2}\|\mathbf{W}-\mathbf{I}\|^2+\nwhat{\psi}(\mathbf{W})+\frac{1}{2}\pier{\|}\grad \mathbf{W}\pier{\|}^2+\frac{1}{2}\pier{\|}\grad \mathbf{R}\pier{\|}^2\\
& \displaystyle +\pier{\int_{\mathcal{D}_a}k |\curl \mathbf{Z}|}+\frac{1}{2}\pier{\|}\curl \mathbf{Z}\pier{\|}^2+\frac{1}{2}\pier{\|} \grad \grad \mathbf{W}\pier{\|}^2,
\end{align}  
where $k> 0$ is a material parameter, representing a threshold for the norm of the reaction term $\boldsymbol{\Gamma+\Pi}$ below which there is compatibility in the system, and
\begin{align}
\label{eqn:23}
\nwhat{\psi}(\mathbf{W}):=\int_{\mathcal{D}_a}I_{\text{SPD}_{\alpha}}(\mathbf{W}),
\end{align}
where $I_{\text{SPD}_{\alpha}}$ is the indicator function of the set
\begin{align}
\label{eqn:24}
\text{SPD}_{\alpha}:=\{\mathbf{W}\in \pier{Sym(\mathbb{R}^{3\times 3})}: \text{det}\mathbf{W}\geq \alpha^3,\; \text{tr}(\text{cof}\mathbf{W})\geq 3\alpha^2,\; \text{tr}\mathbf{W}\geq 3\alpha\},
\end{align}
which is closed and convex for any $\alpha \geq 0$ \cite{MM}.
If $\alpha> 0$ the elements of  $\text{SPD}_{\alpha}$ are positive definite matrices, with all their eigenvalues being \abramohhb{not smaller than $\alpha$ at the same time}. 
The functional \eqref{eqn:23} may be written also as
\begin{align}
\label{eqn:25}
\nwhat{\psi}(\mathbf{W})=\int_{\mathcal{D}_a}I_{S}(\mathbf{W})+\int_{\mathcal{D}_a}I_{C_{\alpha}}(\mathbf{W}),
\end{align}
where \pier{$I_{S}$} is the indicator function of the set of symmetric matrices and \pier{$I_{C_{\alpha}}$} is the indicator function of the set
\begin{align}
\label{eqn:26}
\text{C}_{\alpha}:=\{\mathbf{W}\in M(\mathbb{R}^{3\times 3}): \text{det}\mathbf{W}\geq \alpha^3,\; \text{tr}(\text{cof}\mathbf{W})\geq 3\alpha^2,\; \text{tr}\mathbf{W}\geq 3\alpha\}
\end{align}
for $\alpha> 0$. 
We observe that the function $I_{S}(\mathbf{W})+I_{C_{\alpha}}(\mathbf{W})$ is a convex and l.s.c. function \cite{MM}.
\newline
\abramonn{
In the case $\mathbf{W} \in \mathcal{S}$, we have that 
$$\nwhat{\psi}(\mathbf{W})\equiv  \psi_{C_\alpha} ( {\mathbf{W}} ):=\int_{\mathcal{D}_a}I_{C_\alpha}({\mathbf{W}}). $$
}
\pier{Let us introduce for future convenience the notation}
\begin{equation}
\label{eqn:22bis}
 \pier{\psi_D(\mathbf{A}):= \pier{\int_{\mathcal{D}_a}k |\mathbf{A}}| +\frac{1}{2}\pier{\|}\mathbf{A}\pier{\|}^2, \quad \hbox{for all } \,\mathbf{A}\in L^2(\mathcal{D}_a;\mathbb{R}^{3\times 3}).}
\end{equation}
\newline
\pier{Moreover, we} assume the following form for the dissipation potential of the system, \abramohhb{containing viscous contributions}:
 \begin{align}
\label{eqn:27}
& \displaystyle \notag D(\overbigdot{\mathbf{W}},\boldsymbol{\Omega}):=\frac{1}{2}\pier{\|}\overbigdot{\mathbf{W}}\pier{\|}^2+\frac{1}{2}\pier{\|}\grad \boldsymbol{\Omega}\pier{\|}^2 + \frac{1}{2}\pier{\|}\grad\overbigdot{\mathbf{W}}\pier{\|}^2\\
& \displaystyle   \quad +\frac{1}{2}\pier{\|} \grad \grad \overbigdot{\mathbf{W}}\pier{\|}^2+\frac{1}{2}\pier{\|} \grad \grad  \boldsymbol{\Omega}\pier{\|}^2,
\end{align}
Using \eqref{eqn:13}, \eqref{eqn:22} and \eqref{eqn:27} in \eqref{eqn:21}, and observing from \eqref{eqn:5} that 
\[
\overbigdot{\mathbf{Z}}=\mathcal{G}_{L,\diver}\ast\curl\left(\mathbf{R}\overbigdot{\mathbf{W}}+{\boldsymbol{\Omega}}\mathbf{R}\mathbf{W}\right),
\]
we obtain the following general constitutive assumptions
\begin{equation}
\label{eqn:29}
\displaystyle \mathbf{R}^T\boldsymbol{\Pi}=\pier{\mathbf{W}-\mathbf{I} 
+\boldsymbol{\chi}_{\alpha}+\overbigdot{\mathbf{W}}},
\end{equation}
where 
\pier{$\boldsymbol{\chi}_{\alpha} \in \partial \pier{\nwhat{\psi}(\mathbf{W})}$;}
\begin{equation}
\label{eqn:30}
\displaystyle \mathbf{M}=2\boldsymbol{\Pi}\mathbf{W}\mathbf{R}^T;
\end{equation}
\begin{equation}
\label{eqn:31}
\displaystyle \boldsymbol{\Sigma}:=-(\boldsymbol{\Gamma}+\boldsymbol{\Pi})\in \partial \psi_D(\curl \mathbf{Z})
= \pier{%
\begin{cases}
\displaystyle k\frac{\curl \mathbf{Z}}{|\curl\mathbf{Z}|}+\curl\mathbf{Z}\quad \hbox{if }\, |\curl\mathbf{Z}| \not= 0,\\
\hbox{any } \, \mathbf{M}_D, \, \hbox{ with } |\mathbf{M}_D|\leq k,\, \mathbf{M}_D\wedge \vec{n}|_{\Gamma_D}=\mathbf{0},  \, \hbox{ if }\, |\curl\mathbf{Z}| = 0;
\end{cases}
}
\end{equation}
\begin{equation}
\label{eqn:32}
\begin{cases}
 \displaystyle \mathbf{X}=\grad \mathbf{W}+\grad \overbigdot{\mathbf{W}};\\
 \displaystyle  \mathbf{Y}=\grad\grad \mathbf{W}+\grad\grad \overbigdot{\mathbf{W}};\\
 \displaystyle \boldsymbol{\Lambda}=(\grad \mathbf{R})\mathbf{R}^T+\grad \boldsymbol{\Omega};\\
\displaystyle \mathbf{C}=\grad \left((\grad \mathbf{R})\mathbf{R}^T\right)+\grad\grad \boldsymbol{\Omega}.
\end{cases}
\end{equation}
We observe that the constitutive assumption \eqref{eqn:31} with boundary conditions of \eqref{eqn:20} satisfies \eqref{pigamma}.
We remark that the constitutive assumptions \eqref{eqn:29}--\eqref{eqn:32} \pier{comply with} the principle of objectivity (see \cite{agosti4} for details). 
\noindent
\abramonn{
We now introduce the variable $\boldsymbol{\Theta}(\vec{a},t):=\int_0^t\boldsymbol{\Omega}(\vec{a},s)ds$, 
{$(\vec{a}, t) \in \mathcal{D}_a\times (0,T)$}, and observe that, given $\boldsymbol{\Omega}\in \mathcal{A}$, the differential equation $\overbigdot{\mathbf{R}}=\boldsymbol{\Omega}\mathbf{R}$, with the initial condition $\mathbf{R}(\vec{a},0)=\mathbf{I}$ in \eqref{eqn:20wic}, uniquely {define} a rotation tensor 
\[
\mathbf{R}(\vec{a},t)=\mathrm{e}^{\boldsymbol{\Theta}(\vec{a},t)}\quad {\text{for} \; (\vec{a}, t) \in \mathcal{D}_a\times (0,T)}.
\]
Since $\mathbf{R}:\mathbf{R}=\mathrm{e}^{\boldsymbol{\Theta}(\vec{a},t)}\mathrm{e}^{-\boldsymbol{\Theta}(\vec{a},t)}:\mathbf{I}=3$, we have that
\begin{equation}
\label{eqn:39}
\mathbf{R}
\in L^{\infty}(\mathcal{D}_{aT},\mathbb{R}^{3\times 3}).
\end{equation}
With the latter change of variables, inserting \eqref{eqn:29}--\eqref{eqn:32} in \eqref{eqn:20w1} we finally obtain the variational formulation
\begin{equation}
\label{eqn:48}
\begin{cases}
\displaystyle \int_{\mathcal{D}_a}\mathrm{e}^{-\boldsymbol{\Theta}}\grad \left(\mathcal{G}_L\ast \left(\twobigdot{\vec{u}}-\vec{\mathcal{F}}_{\textrm ext}\right)\right):\nwhat{\mathbf{W}}+\int_{\mathcal{D}_a}\left(\twobigdot{\mathbf{W}}+\overbigdot{\mathbf{W}}+\mathbf{W}-\mathbf{I}+{{}\boldsymbol{\chi}_{\alpha}}\right):\nwhat{\mathbf{W}}\\
\displaystyle \;\; + \int_{\mathcal{D}_a}\mathrm{e}^{-\boldsymbol{\Theta}}\curl\left(\mathcal{G}_{L,\diver}\ast\left(\curl 
\boldsymbol{\Sigma}\right)\right):\nwhat{\mathbf{W}}+\int_{\mathcal{D}_a}\grad\left(\twobigdot{\mathbf{W}}+\overbigdot{\mathbf{W}}+\mathbf{W}\right)::\grad \nwhat{\mathbf{W}}\\
\displaystyle + \int_{\mathcal{D}_a}\grad \grad \left(\twobigdot{\mathbf{W}}+\overbigdot{\mathbf{W}}+\mathbf{W}\right)::: \grad \grad  \nwhat{\mathbf{W}}=\int_{\mathcal{D}_a}\mathbf{W}_{\text{ext}}:\nwhat{\mathbf{W}},\\ \\
\displaystyle \int_{\mathcal{D}_a}\left( \grad \left(\mathcal{G}_L\ast \left(\twobigdot{\vec{u}}-\vec{\mathcal{F}}_{\textrm ext}\right)\right)\mathbf{W}\mathrm{e}^{-\boldsymbol{\Theta}}+\twobigdot{\boldsymbol{\Theta}}\right):\nwhat{\boldsymbol{\Omega}}+\int_{\mathcal{D}_a}\curl\left(\mathcal{G}_{L,\diver}\ast\left(\curl \boldsymbol{\Sigma}\right)\right)\mathbf{W}\mathrm{e}^{-\boldsymbol{\Theta}}:\nwhat{\boldsymbol{\Omega}}\\
\displaystyle \;\; +\int_{\mathcal{D}_a}\grad\left(\twobigdot{\boldsymbol{\Theta}}+\frac{1}{2}\overbigdot{\boldsymbol{\Theta}}+\frac{1}{2}{\boldsymbol{\Theta}}\right)::\grad \nwhat{\boldsymbol{\Omega}}+\int_{\mathcal{D}_a}\grad\grad\left(\twobigdot{\boldsymbol{\Theta}}+\frac{1}{2}\overbigdot{\boldsymbol{\Theta}}+\frac{1}{2}{\boldsymbol{\Theta}}\right)::: \grad \grad \nwhat{\boldsymbol{\Omega}}\\
\displaystyle \;\;=\int_{\mathcal{D}_a}\boldsymbol{\Omega}_{\text{ext}}:\nwhat{\boldsymbol{\Omega}},\\ \\
\displaystyle {\boldsymbol{\chi}_{\alpha} \in \partial {\nwhat{\psi}(\mathbf{W})}, \quad  
\boldsymbol{\Sigma}\in \partial \psi_D(\curl \mathbf{Z}),}\\ \\
\displaystyle \int_{\mathcal{D}_a}\grad{}{{\vec{u}}}\colon \grad{}\vec{v}=\int_{\mathcal{D}_a}\left(\mathrm{e}^{\boldsymbol{\Theta}}\mathbf{W}-\mathbf{I}\right)\colon \grad{}\vec{v},\\ \\
\displaystyle \int_{\mathcal{D}_a}\grad {\mathbf{Z}}:: \grad \nwhat{\mathbf{Z}}=\int_{\mathcal{D}_a}\left(\mathrm{e}^{\boldsymbol{\Theta}}\mathbf{W}-\mathbf{I}\right)\colon \curl \nwhat{\mathbf{Z}},
\end{cases}
\end{equation}
valid for all choices of $\nwhat{\mathbf{W}}\in \mathcal{S}$, $\nwhat{\boldsymbol{\Omega}} \in \mathcal{A}$, $\vec{v}\in \mathcal{V}$ and $\nwhat{\mathbf{Z}}\in \mathcal{M}_{div}$, with $\nwhat{\mathbf{W}}|_{\Gamma_D}=\nwhat{\boldsymbol{\Omega}}|_{\Gamma_D}=\boldsymbol{0}$, $\grad\nwhat{\mathbf{W}}|_{\Gamma_D}=\grad\nwhat{\boldsymbol{\Omega}}|_{\Gamma_D}=\mathbf{0}$, $\vec{v}|_{\Gamma_D}=\vec{0}$, $\nwhat{\mathbf{Z}}|_{\Gamma_N}=\mathbf{0}$, with boundary conditions
{
\begin{align}
\label{eqn:34}
& \displaystyle \notag \mathbf{W}=\mathbf{I},\;\boldsymbol{\Theta}=\mathbf{0}, \grad \mathbf{W}=\grad \boldsymbol{\Theta}=\mathbf{0} \;\; \text{on} \; \Gamma_D\times (0,T), \\
& \displaystyle \notag \vec{u} =\mathbf{0}, \; \quad \text{on} \; \Gamma_D\times (0,T), 
\\ 
& \displaystyle \mathbf{Z}=\mathbf{0}, \; \quad \text{on} \; \Gamma_N\times (0,T),
\end{align}
}
 and initial conditions
\begin{equation}
\label{eqn:35}
\mathbf{W}(\vec{a},0)=\mathbf{I}, \; \overbigdot{\mathbf{W}}(\vec{a},0)=\mathbf{0}, \; \boldsymbol{\Theta}=\overbigdot{\boldsymbol{\Theta}}=\boldsymbol{0},\quad {\text{for} \; \vec{a} \in \mathcal{D}_a}.
\end{equation}
We observe that an initial condition for $\mathbf{Z}$ could be also defined by assuming that \eqref{eqn:48}$_5$ is valid for $t=0$, {i.e.,}
\begin{equation}
\label{eqn:50z}
\mathbf{Z}(\vec{a},0)= (\mathcal{G}_{L,\diver}\ast \curl \mathbf{W}(\vec{a},0)) { \quad \text{for} \; \vec{a}  \in \mathcal{D}_a}.
\end{equation}
In the case {$\mathbf{W}(\vec{a},0)=\mathbf{I}$, then we have} $\mathbf{Z}(\vec{a},0)=\boldsymbol{0}$ ${ \ \text{for} \; \vec{a}  \in \mathcal{D}_a}$. A similar argument can be applied to obtain the initial conditions for $\vec{u}$ and $\overbigdot{\vec{u}}$.
}
\newline
We state now the main theorem of the present paper. We start by introducing the following assumptions on the data:
\begin{itemize}
  \item[\textbf{A1}:] $\mathcal{D}_a\subset \mathbb{R}^3$ is an open bounded and simply connected domain with Lipschitz boundary $\Gamma_a:=\partial \mathcal{D}_a$. Moreover, $\Gamma_a=\Gamma_D\cup \Gamma_N$, where $\Gamma_D,\Gamma_N$ are connected Lipschitz subsets of $\Gamma_a$ with positive measures and such that $|\Gamma_D\cap \Gamma_N|=0$;
  \item[\textbf{A2}:] The initial data are $\mathbf{W}_0=\mathbf{I}$, $\overbigdot{\mathbf{W}}_0=\mathbf{0}$, $\boldsymbol{\Theta}_0=\overbigdot{\boldsymbol{\Theta}}_0=\mathbf{0}$, $\vec{u}_0=\overbigdot{\vec{u}}_0=\vec{0}$, $\mathbf{Z}_0=\mathbf{0}$. Note that the initial data for $\vec{u}$ and $\mathbf{Z}$ are given as the solutions of the elliptic problems \eqref{eqn:48}$_4$, \eqref{eqn:48}$_5$ with right hand sides written at time $t=0$;
  \item[\textbf{A3}:] The forcing term \pier{$\pier{\mathbf{W}_{\text{ext}}}: H^2(\mathcal{D}_a;Sym(\mathbb{R}^{3\times 3}))\times (0,T)\rightarrow L^2(\mathcal{D}_a;Sym(\mathbb{R}^{3\times 3}))$ is measurable in $t\in (0,T)$ and continuous in $\mathbf{W} \in H^2(\mathcal{D}_a;Sym(\mathbb{R}^{3\times 3}))$, \pier{and it satisfies} 
  \[
   \pier{\|}\pier{\mathbf{W}_{\text{ext}}}(\mathbf{W},t)\pier{\|}_{L^2(\mathcal{D}_a;Sym(\mathbb{R}^{3\times 3}))}\leq L\left(\pier{\|}\mathbf{W}\pier{\|}_{H^2(\mathcal{D}_a;Sym(\mathbb{R}^{3\times 3}))}+1\right),
  \]
  for a.e. $t\in (0,T)$, for all $\mathbf{W}\in H^2(\mathcal{D}_a;Sym(\mathbb{R}^{3\times 3}))$ and for some $L\in \mathbb{R}$.} 
  \abramo{Similarly, the forcing term ${\boldsymbol{\Omega}_{\text{ext}}}: H^2(\mathcal{D}_a;Skew(\mathbb{R}^{3\times 3}))\times (0,T)\rightarrow L^2(\mathcal{D}_a;Skew(\mathbb{R}^{3\times 3}))$ is measurable in $t\in (0,T)$ and continuous in $\boldsymbol{\Theta} \in H^2(\mathcal{D}_a;Skew(\mathbb{R}^{3\times 3}))$, {and it satisfies} 
  \[
   \pier{\|}\pier{\boldsymbol{\Omega}_{\text{ext}}}(\boldsymbol{\Theta},t)\pier{\|}_{L^2(\mathcal{D}_a;Skew(\mathbb{R}^{3\times 3}))}\leq G\left(\pier{\|}\boldsymbol{\Theta}\pier{\|}_{H^2(\mathcal{D}_a;Skew(\mathbb{R}^{3\times 3}))}+1\right),
  \]
  for a.e. $t\in (0,T)$, for all $\boldsymbol{\Theta}\in H^2(\mathcal{D}_a;Skew(\mathbb{R}^{3\times 3}))$ and for some $G\in \mathbb{R}$.} 
  \newline 
  Finally, $\vec{\mathcal{F}}_{\textrm ext}\in L^{\infty}\left(0,T;\left(H^1_{\Gamma_D}(\mathcal{D}_a,\mathbb{R}^3)\right)'\right)$.
  \end{itemize}
\begin{thm}
 \label{thm:1}
 Let assumptions \textbf{A1}-\textbf{A3} be satisfied.  Then, there exist a $\hat{T}$, with $0<\hat{T}\leq T$, and a \pco{quintuplet}
 $(\mathbf{W},\boldsymbol{\Theta},\boldsymbol{\Sigma},\pier{\vec{u}},\mathbf{Z})$, with
\begin{align}
\label{eqn:57thm}
\mathbf{W}\in W^{1,\infty}\cap H^2(0,\hat{T};\bar{H}_{\Gamma_D}^2(\mathcal{D}_a,Sym(\mathbb{R}^{3\times 3}))),
\end{align}
and \pier{\rm $\mathbf{W}(\vec{a},t)\in \text{SPD}_{\alpha}$} for all $(\vec{a},t)\in \mathcal{D}_{a\hat{T}}$,  
\begin{equation}
\label{eqn:58thm}
\boldsymbol{\Theta}\in W^{1,\infty}\cap H^2(0,\hat{T};\bar{H}_{\Gamma_D}^2(\mathcal{D}_a,Skew(\mathbb{R}^{3\times 3}))),
\end{equation}
\begin{equation}
 \label{eqn:59thm}
 \boldsymbol{\Sigma}\in L^{\infty}(0,\hat{T};L^2(\mathcal{D}_a;\mathbb{R}^{3\times 3})),
\end{equation}
\begin{equation}
 \label{eqn:60thm}
 \vec{u}\in W^{1,\infty}\cap H^2(0,\hat{T};H_{\Gamma_D}^1(\mathcal{D}_a,\mathbb{R}^{3})),
\end{equation}
\begin{equation}
 \label{eqn:61thm}
 \mathbf{Z}\in W^{1,\infty}\cap H^2(0,\hat{T};{H}_{\Gamma_N}^1(\mathcal{D}_a,\mathbb{R}^{3\times 3})\cap \mathcal{M}_{div}),
\end{equation}
which solve the following weak formulation associated to \abramonew{\eqref{eqn:48}}:
\begin{equation}
\label{eqn:thm1}
\begin{cases}
\displaystyle \int_{\mathcal{D}_a}\textrm{Sym}\left(\mathrm{e}^{-\boldsymbol{\Theta}}\grad \left(\mathcal{G}_L\ast \twobigdot{\vec{u}}\right)\right)\colon \nwhat{\mathbf{W}}+ \int_{\mathcal{D}_a}\twobigdot{\mathbf{W}}\colon \nwhat{\mathbf{W}}\\
\displaystyle{}+\int_{\mathcal{D}_a}\pier{\grad{}} \twobigdot{\mathbf{W}}::\pier{\grad{}}  \nwhat{\mathbf{W}}+\int_{\mathcal{D}_a}\pier{\grad\grad{}} \twobigdot{\mathbf{W}}\colon \pier{\grad\grad{}}  \nwhat{\mathbf{W}}\\
\displaystyle{}+\int_{\mathcal{D}_a}\textrm{Sym}\left(\mathrm{e}^{-\boldsymbol{\Theta}}\curl\left(
\mathcal{G}_{L,\diver}\ast \curl \boldsymbol{\Sigma}\right)\right)\colon  \nwhat{\mathbf{W}}\\
\displaystyle{}+ \int_{\mathcal{D}_a}\left(\mathbf{W}-\mathbf{I}+\frac{\pier{{}d\tilde{\psi}_{C_\alpha}}}{d \mathbf{W}}(\mathbf{W})+\overbigdot{\mathbf{W}}\right)\colon \nwhat{\mathbf{W}}+\int_{\mathcal{D}_a}\pier{\grad{}} (\mathbf{W}+\overbigdot{\mathbf{W}})::\pier{\grad{}} \nwhat{\mathbf{W}}
\\
\displaystyle{}+\int_{\mathcal{D}_a}\pier{\grad\grad{}} (\mathbf{W}+\overbigdot{\mathbf{W}})\colon \pier{\grad\grad{}}\nwhat{\mathbf{W}}=\int_{\mathcal{D}_a}\mathbf{W}_{\text{ext}}(\mathbf{W},t)\colon \nwhat{\mathbf{W}}\\
\displaystyle{}+\int_{\mathcal{D}_a}\textrm{Sym}\left(\mathrm{e}^{-\boldsymbol{\Theta}}\grad \left(\mathcal{G}_L\ast \vec{\mathcal{F}}_{\textrm ext}\right)\right)\colon \nwhat{\mathbf{W}},\\ \\
\displaystyle \int_{\mathcal{D}_a}\textrm{Skew}\left(\grad \left(\mathcal{G}_L\ast \twobigdot{\vec{u}}\right)\mathbf{W}\mathrm{e}^{-\boldsymbol{\Theta}}\right)\colon \nwhat{\boldsymbol{\Omega}}+\int_{\mathcal{D}_a}\twobigdot{\boldsymbol{\Theta}}\colon \nwhat{\boldsymbol{\Omega}}\\
\displaystyle{}+\int_{\mathcal{D}_a}\grad{} \twobigdot{\boldsymbol{\Theta}}::\pier{\grad{}} \nwhat{\boldsymbol{\Omega}}+\int_{\mathcal{D}_a}\pier{\grad\grad{}} \twobigdot{\boldsymbol{\Theta}}\colon \pier{\grad\grad{}} \nwhat{\boldsymbol{\Omega}}\\
\displaystyle{}+\int_{\mathcal{D}_a}\textrm{Skew}\left(\left(\curl\left(
\mathcal{G}_{L,\diver}\ast \curl \boldsymbol{\Sigma}\right)\mathbf{W}\mathrm{e}^{-\boldsymbol{\Theta}}\right)\right)\colon \nwhat{\boldsymbol{\Omega}} \\
\displaystyle{}+\frac{1}{2}\int_{\mathcal{D}_a}\pier{\grad{}} (\boldsymbol{\Theta}+\overbigdot{\boldsymbol{\Theta}})::\pier{\grad{}} \nwhat{\boldsymbol{\Omega}}+\frac{1}{2}\int_{\mathcal{D}_a}\pier{\grad\grad{}}(\boldsymbol{\Theta}+ \overbigdot{\boldsymbol{\Theta}})\colon\pier{\grad\grad{}} \nwhat{\boldsymbol{\Omega}}\\
\displaystyle=\int_{\mathcal{D}_a}\boldsymbol{\Omega}_{\text{ext}}\abramo{(\boldsymbol{\Theta},t)}\colon \nwhat{\boldsymbol{\Omega}}+\int_{\mathcal{D}_a}\textrm{Skew}\left(\grad \left(\mathcal{G}_L\ast \vec{\mathcal{F}}_{\textrm ext}\right)\mathbf{W}\mathrm{e}^{-\boldsymbol{\Theta}}\right)\colon \nwhat{\boldsymbol{\Omega}},\\ \\
\displaystyle \int_{\mathcal{D}_a}(\curl \nwhat{\mathbf{Z}}-\curl \mathbf{Z}):\boldsymbol{\Sigma}+\int_{\mathcal{D}_a}\psi_D(\curl \mathbf{Z})\leq \int_{\mathcal{D}_a}\psi_D(\curl \nwhat{\mathbf{Z}}),\\ \\
\displaystyle \int_{\mathcal{D}_a}\grad{}{{\vec{u}}}\colon \grad{}\vec{v}=\int_{\mathcal{D}_a}\left(\mathrm{e}^{\boldsymbol{\Theta}}\mathbf{W}-\mathbf{I}\right)\colon \grad{}\vec{v},\\ \\
\displaystyle \int_{\mathcal{D}_a}\grad {\mathbf{Z}}:: \grad \nwhat{\mathbf{Z}}=\int_{\mathcal{D}_a}\left( \mathrm{e}^{\boldsymbol{\Theta}}\mathbf{W}-\mathbf{I}\right)\colon \curl \nwhat{\mathbf{Z}},
\end{cases}
\end{equation}
for a.e. $t \in [0,\hat{T}]$, \abramonn{where $\tilde{\psi}_{C_\alpha}$ is a proper regularization of ${\psi}_{C_\alpha}$ which will be introduced later}, for \pier{all choices of} $\nwhat{\mathbf{W}} \in \bar{H}_{\Gamma_D}^2(\mathcal{D}_a,Sym(\mathbb{R}^{3\times 3}))$, $\nwhat{\boldsymbol{\Omega}} \in \bar{H}^2_{\Gamma_D}(\mathcal{D}_a,Skew(\mathbb{R}^{3\times 3}))$, $\vec{v}\in H^1_{\Gamma_D}(\mathcal{D}_a,\mathbb{R}^3)$ and $\nwhat{\mathbf{Z}}\in H_{\Gamma_N}^1(\mathcal{D}_a,\mathbb{R}^{3\times 3})\cap \mathcal{M}_{div}$,
and with initial conditions
\begin{equation}
    \label{eqn:thmic}
   \pier{ \mathbf{W}(\cdot,0)=\mathbf{I}, \; \overbigdot{\mathbf{W}}(\cdot,0)=\mathbf{0}, \quad \boldsymbol{\Theta} (\cdot,0)=\overbigdot{\boldsymbol{\Theta}} (\cdot,0)=\boldsymbol{0} \quad  \hbox{in }
   \, \mathcal{D}_a, \quad \vec{u}(\cdot,0)=\overbigdot{\vec{u}}(\cdot,0)=\vec{0}\quad  \hbox{a.e. in}
   \, \mathcal{D}_a.}
\end{equation}
\end{thm}
\section{Proof of the main result}
\label{sec:analysis}
In this section we prove Theorem \ref{thm:1}. The strategy of the proof is the following: since the higher order time derivatives are nonlinearly coupled with lower order terms in the inertia, we need to introduce a time regularization in the system in order to be able to prove the existence of a solution. Also,  we introduce proper regularizations to deal with subdifferentials. Further, a Faedo--Galerkin discretization of the regularized system is introduced, which leads to the proof of existence of a solution and to the derivation of a-priori estimates, uniform in the regularization and discretization parameters, which let us identify a solution in the continuum and unregularized limit.
\subsection{Regularization}
In order to proceed, we introduce the three following level of regularizations for specific components of system  \eqref{eqn:48}.
\begin{itemize}
\item We introduce a decreasing non negative smooth approximation $\tilde{I}_{C_{\alpha}}$ of the indicator function $I_{C_{\alpha}}$ from the interior of its effective domain. For instance, we may introduce the function
\begin{equation}
\label{eqn:itilde}
\tilde{I}_{+}^{\beta}(x):=
\begin{cases}
\displaystyle \frac{\widetilde{(\beta-x)}_+}{x}, \quad \text{if} \; x>0,\\
\displaystyle +\infty, \quad \text{if} \; x\leq 0,
\end{cases}
\end{equation}
for a given $\beta>0$, where $\widetilde{(\beta-x)}_+$ is a smooth regularization of the positive part function $(\cdot)_+:=\max(0,\cdot)$, and define 
\begin{equation}
\label{eqn:ica}
\tilde{I}_{C_{\alpha}}(\mathbf{W}):=\tilde{I}_{+}^{1-\alpha^3}(\text{det}\mathbf{W}-\alpha^3)+\tilde{I}_{+}^{3-3\alpha^2}(\text{tr}(\text{cof}\mathbf{W})-3\alpha^2)+\tilde{I}_{+}^{3-3\alpha}(\text{tr}\mathbf{W}-3\alpha).
\end{equation}
We observe that
\begin{equation}
\label{eqn:ica2}
\tilde{I}_{C_{\alpha}}(\mathbf{I})=0.
\end{equation}
Then, we define  $\tilde{\psi}_{C_\alpha} ( \nwhat{\mathbf{W}} ):=\int_{\mathcal{D}_a}\tilde{I}_{C_\alpha}(\nwhat{\mathbf{W}})$, and $\tilde{\chi}_{\alpha}({\mathbf{W}})=\frac{d \tilde{\psi}_{C_\alpha}}{d {\mathbf{W}}}( {\mathbf{W}} )$.
\item We replace the convex function $\psi_D$ and its subdifferential $\partial \psi_D$ by their Moreau--Yosida \pier{approximations}  $\psi_D^{\lambda}$ and $\partial \psi_D^{\lambda}$, depending on a regularization parameter $\lambda >0$. 
We refer to, e.g., \cite[pp.~28 and~39]{brezis}) for definitions and properties of these approximations, recalling simply that if  $f:L^2(\mathcal{D}_a;\mathbb{R}^{3\times 3})\to [0,+\infty] $ is a proper convex lower semicontinuous function and $\partial f$ denotes its subdifferential, then  
\[
 \partial f^{\lambda} :=\frac{I-\left(I+\lambda \partial f \right)^{-1}}{\lambda}, \quad \lambda \in (0,1),
\]
where $I$ here denotes the identity operator. In particular, $\partial f^{\lambda}$ is a monotone and $\frac{1}{\lambda}$-Lipschitz continuous function. Moreover, due the special form of $\psi_D$ defined in \eqref{eqn:22bis},
we have that the following bounds are valid uniformly in $\lambda$:
\begin{align}
\label{eqn:dpsi1}
& \displaystyle \frac{1}{2}\pier{\|}\mathbf{A}\pier{\|}^2\leq C+\psi_D^{\lambda}(\mathbf{A}), \quad \hbox{for all } \, \mathbf{A}\in L^2(\mathcal{D}_a;\mathbb{R}^{3\times 3}),\\
\label{eqn:dpsi2}
& \displaystyle \pier{\|}\partial \psi_D^{\lambda}(\mathbf{A})\pier{\|}^2\leq C\bigl( \psi_D^{\lambda}(\mathbf{A})+1\bigr), \quad \hbox{for all } \, \mathbf{A}\in L^2(\mathcal{D}_a;\mathbb{R}^{3\times 3}).
\end{align}
\item We add a time regularization term $-\lambda \overbigdot{\vec{u}}$, with $\lambda$ the same regularization parameter used to define the Moreau--Yosida \pier{approximations}  $\psi_D^{\lambda}$, to the left hand side of \eqref{eqn:48}$_4$.
\end{itemize}
Given $\lambda>0$, we then introduce the following regularized version of problem \eqref{eqn:48}
\abramonn{
\begin{equation}
\label{symskewreg}
\begin{cases}
\displaystyle \int_{\mathcal{D}_a}\mathrm{e}^{-\boldsymbol{\Theta}}\grad \left(\mathcal{G}_L\ast \left(\twobigdot{\vec{u}}-\vec{\mathcal{F}}_{\textrm ext}\right)\right):\nwhat{\mathbf{W}}+\int_{\mathcal{D}_a}\left(\twobigdot{\mathbf{W}}+\overbigdot{\mathbf{W}}+\mathbf{W}-\mathbf{I}+\frac{d \tilde{\psi}_{C_\alpha}}{d {\mathbf{W}}}( {\mathbf{W}} )\right):\nwhat{\mathbf{W}}\\
\displaystyle \;\; + \int_{\mathcal{D}_a}\mathrm{e}^{-\boldsymbol{\Theta}}\curl\left(\mathcal{G}_{L,\diver}\ast\left(\curl 
\boldsymbol{\Sigma}\right)\right):\nwhat{\mathbf{W}}+\int_{\mathcal{D}_a}\grad\left(\twobigdot{\mathbf{W}}+\overbigdot{\mathbf{W}}+\mathbf{W}\right)::\grad \nwhat{\mathbf{W}}\\
\displaystyle + \int_{\mathcal{D}_a}\grad \grad \left(\twobigdot{\mathbf{W}}+\overbigdot{\mathbf{W}}+\mathbf{W}\right)::: \grad \grad  \nwhat{\mathbf{W}}=\int_{\mathcal{D}_a}\mathbf{W}_{\text{ext}}:\nwhat{\mathbf{W}},\\ \\
\displaystyle \int_{\mathcal{D}_a}\left( \grad \left(\mathcal{G}_L\ast \left(\twobigdot{\vec{u}}-\vec{\mathcal{F}}_{\textrm ext}\right)\right)\mathbf{W}\mathrm{e}^{-\boldsymbol{\Theta}}+\twobigdot{\boldsymbol{\Theta}}\right):\nwhat{\boldsymbol{\Omega}}+\int_{\mathcal{D}_a}\curl\left(\mathcal{G}_{L,\diver}\ast\left(\curl \boldsymbol{\Sigma}\right)\right)\mathbf{W}\mathrm{e}^{-\boldsymbol{\Theta}}:\nwhat{\boldsymbol{\Omega}}\\
\displaystyle \;\; +\int_{\mathcal{D}_a}\grad\left(\twobigdot{\boldsymbol{\Theta}}+\frac{1}{2}\overbigdot{\boldsymbol{\Theta}}+\frac{1}{2}{\boldsymbol{\Theta}}\right)::\grad \nwhat{\boldsymbol{\Omega}}+\int_{\mathcal{D}_a}\grad\grad\left(\twobigdot{\boldsymbol{\Theta}}+\frac{1}{2}\overbigdot{\boldsymbol{\Theta}}+\frac{1}{2}{\boldsymbol{\Theta}}\right)::: \grad \grad \nwhat{\boldsymbol{\Omega}}\\
\displaystyle \;\;=\int_{\mathcal{D}_a}\boldsymbol{\Omega}_{\text{ext}}:\nwhat{\boldsymbol{\Omega}},\\ \\
\displaystyle \boldsymbol{\Sigma}\in \partial \psi_D^{\lambda}(\curl \mathbf{Z}),\\ \\
\displaystyle \lambda \int_{\mathcal{D}_a}\overbigdot{\vec{u}}\cdot \vec{v}+\int_{\mathcal{D}_a}\grad{}{{\vec{u}}}\colon \grad{}\vec{v}=\int_{\mathcal{D}_a}\left(\mathrm{e}^{\boldsymbol{\Theta}}\mathbf{W}-\mathbf{I}\right)\colon \grad{}\vec{v},\\ \\
\displaystyle \int_{\mathcal{D}_a}\grad {\mathbf{Z}}:: \grad \nwhat{\mathbf{Z}}=\int_{\mathcal{D}_a}\left(\mathrm{e}^{\boldsymbol{\Theta}}\mathbf{W}-\mathbf{I}\right)\colon \curl \nwhat{\mathbf{Z}},
\end{cases}
\end{equation}
valid for all choices of $\nwhat{\mathbf{W}}\in \mathcal{S}$, $\nwhat{\boldsymbol{\Omega}} \in \mathcal{A}$, $\vec{v}\in \mathcal{V}$ and $\nwhat{\mathbf{Z}}\in \mathcal{M}_{div}$, with $\nwhat{\mathbf{W}}|_{\Gamma_D}=\nwhat{\boldsymbol{\Omega}}|_{\Gamma_D}=\boldsymbol{0}$, $\grad\nwhat{\mathbf{W}}|_{\Gamma_D}=\grad\nwhat{\boldsymbol{\Omega}}|_{\Gamma_D}=\mathbf{0}$, $\vec{v}|_{\Gamma_D}=\vec{0}$, $\nwhat{\mathbf{Z}}|_{\Gamma_N}=\mathbf{0}$,} with the same boundary and initial conditions as \eqref{eqn:34} and \eqref{eqn:35}. For ease of notation, we have not explicitly indicated the dependence of the solutions from the regularization parameters $\lambda$ and $\lambda$.
\subsection{Faedo--Galerkin approximation}
Let us introduce the finite dimensional spaces which will be used to formulate the Galerkin ansatz to approximate the solutions of \pier{the system} \eqref{symskewreg}. As a first step, we introduce the following fourth order elliptic problem with mixed boundary conditions, associated to the operator $\Upsilon:=\diver \Delta \grad-\Delta$,
\begin{equation}
\label{fg1}
\begin{cases}
\Upsilon\vec{v}=\diver \Delta \grad\vec{v}-\Delta\vec{v}=\vec{f} \;\; \text{in}\;\; \mathcal{D}_a,\\
\vec{v}=\vec{0},\; \grad \vec{v}=\mathbf{0} \;\; \text{on}\;\; \Gamma_D,\\
\abramonn{W_N(\grad \vec{v},\grad \grad \vec{v})=\vec{0}, \;K_N(\grad\grad \vec{v})=\vec{0} \;\; \text{on}\;\; \Gamma_N,}
\end{cases}
\end{equation}
for a given $\vec{f}\in L^2(\mathcal{D}_a,\mathbb{R}^3)$. Taking the $L^2$ scalar product of the previous partial differential equation with test functions \abramon{$\vec{w}\in \bar{H}^2_{\Gamma_D}(\mathcal{D}_a,\mathbb{R}^3)$}, using similar integration by parts formula as \eqref{bcparts}-\eqref{byparts1} we obtain the following weak formulation associated to \eqref{fg1}:
\begin{equation}
\label{fg2}
(\grad\grad \vec{v},\grad\grad \vec{w})+(\grad \vec{v},\grad\vec{w})=(\vec{f},\vec{w}), \;\; \text{for any}\;\;\vec{w}\in \bar{H}^2_{\Gamma_D}.
\end{equation}
Thanks to the Lax--Milgram Lemma and the Poincar\'e inequality, there exists a unique weak solution $\vec{v}\in \bar{H}^2_{\Gamma_D}$ to \eqref{fg2}, which also satisfies the Lax--Milgram estimate
\[
||\vec{v}||_{H^2(\mathcal{D}_a,\mathbb{R}^3)}\leq C||\vec{f}||_{L^2(\mathcal{D}_a,\mathbb{R}^3)}.
\]
Hence, the inverse operator $\Upsilon^{-1}:L^2(\mathcal{D}_a,\mathbb{R}^3)\to L^2(\mathcal{D}_a,\mathbb{R}^3)$ is well defined, and since $\vec{u}=\Upsilon^{-1}\vec{f}\in H^2(\mathcal{D}_a,\mathbb{R}^3)\subset\subset L^2(\mathcal{D}_a,\mathbb{R}^3)$, it is compact. It is also self-adjoint. Indeed, given $\vec{f},\vec{g}\in L^2(\mathcal{D}_a,\mathbb{R}^3)$, with $\vec{u}=\Upsilon^{-1}\vec{f}$, $\vec{w}=\Upsilon^{-1}\vec{g}$, we have that
\[
(\Upsilon^{-1}\vec{f},\vec{g})=(\vec{u},\vec{g})=(\grad\grad \vec{u},\grad\grad \vec{w})+(\grad \vec{u},\grad\vec{w})=(\vec{w},\vec{f})=(\vec{f},\Upsilon^{-1}\vec{g}).
\]
Hence, $\Upsilon^{-1}$ admits a countable set of eigenvectors $\{\vec{\xi}_i\}_{i\in \mathbb{N}}$, i.e. $\Upsilon^{-1}\vec{\xi}_i=\mu_i \vec{\xi}_i$, which is an orthonormal basis of $L^2(\mathcal{D}_a,\mathbb{R}^3)$ and an orthogonal basis in $\bar{H}^2_{\Gamma_D}(\mathcal{D}_a,\mathbb{R}^3)$. Setting $\gamma_i=\mu_i^{-1}$, we then define the eigenfunctions $\{\vec{\xi}_i\}_{i\in \mathbb{N}}$ of the elliptic operator with mixed boundary conditions
\[
\abramon{
\begin{cases}
\diver \Delta \grad \vec{\xi}_i-\Delta \vec{\xi}_i=\gamma_i \vec{\xi}_i \quad \text{in} \; \mathcal{D}_a,\\
\vec{\xi}_i=\vec{0},\; \grad \vec{\xi}_i=\mathbf{0} \;\; \text{on}\;\; \Gamma_D,\\
\abramonn{W_N(\grad \vec{\xi}_i,\grad \grad \vec{\xi}_i)=\vec{0}, \;K_N(\grad\grad \vec{\xi}_i)=\vec{0} \;\; \text{on}\;\; \Gamma_N,}
\end{cases}
}
\]
with \pier{$0<\gamma_0\leq \gamma_1 \leq \dots \leq \gamma_m\to \infty$}. We observe that each eigenvalue has geometric multiplicity $3$, i.e. $\gamma_{3k}=\gamma_{3k+1}=\gamma_{3k+2}$ for any $k\in \mathbb{N}$. 
We then introduce the numbers $n_i:=i\bmod 3$, $h_i:=\{0 \; \text{if} \; n_i=0, \; 1 \; \text{if} \; n_i>0\}$, and the functions $\{\mathbf{S}_{i+j+h_i}\}_{i\in \mathbb{N}; i\leq j\leq i+2-n_i}$ defined by
\[
 \mathbf{S}_{i+j+h_i}:=\vec{\xi}_i\otimes \vec{\xi}_j+\vec{\xi}_j\otimes \vec{\xi}_i,
\]
whith $i\in \mathbb{N}; i\leq j\leq i+2-n_i$. We observe that, given $i\in \mathbb{N}$ with $n_i=0$, the elements $(\mathbf{S}_{i},\mathbf{S}_{i+1},\mathbf{S}_{i+2},\mathbf{S}_{i+3},\mathbf{S}_{i+4},\mathbf{S}_{i+5})$ span the $6$-th dimensional linear eigenspace of symmetric tensors associated to the eigenvalue $\gamma_i$. 
We also introduce the projection operator
\[
PS_m: \bar{H}^2_{\Gamma_D}(\mathcal{D}_a,\mathbb{R}^{3\times 3})\to \text{span}\{\mathbf{S}_0,\mathbf{S}_1,\dots,\mathbf{S}_{6m+5}\}.
\]
We then introduce the functions $\{\mathbf{A}_{n_i+j-1}\}_{i\in \mathbb{N}; i< j\leq i+2-n_i}$ defined by
\[
 \mathbf{A}_{n_i+j-1}:=\vec{\xi}_i\otimes \vec{\xi}_j-\vec{\xi}_j\otimes \vec{\xi}_i,
\]
whith $i\in \mathbb{N}; i< j\leq i+2-n_i$. We observe that, given $i\in \mathbb{N}$ with $n_i=0$, the elements $(\mathbf{A}_{i},\mathbf{A}_{i+1},\mathbf{A}_{i+2})$ span the $3$-th dimensional linear eigenspace of anti-symmetric tensors associated to the eigenvalue $\lambda_i$. 
We then introduce the projection operator
\[
PA_m:\pier{\bar{H}^2_{\Gamma_D}(\mathcal{D}_a,\mathbb{R}^{3\times 3})} \to \text{span}\{\mathbf{A}_0,\mathbf{A}_1,\dots,\mathbf{S}_{3m+2}\}.
\]
We finally introduce the eigenfunctions $\{\theta_i\}_{i\in \mathbb{N}}$ of the Laplace operator with mixed boundary conditions, \pier{i.e.,}
\[
-\Delta \theta_i=\rho_i \theta_i \quad \text{in} \; \mathcal{D}_a, \quad \theta_i =0 \quad \text{on} \; \Gamma_D, \quad \grad \theta\cdot \vec{n}=0 \quad \text{on} \; \Gamma_N,
\]
with \pier{$0<\rho_0\leq \rho_1 \leq \dots \leq \rho_m\to \infty$}. The sequence $\{\theta_i\}_{i\in \mathbb{N}}$ can be chosen as an orthonormal basis in $L^2(\mathcal{D}_a)$ and an orthogonal basis in $H_{\Gamma_D}^1(\mathcal{D}_a)$.
We then introduce the functions $\{\vec{v}_{3k+i}\}_{k\in \mathbb{N}; i=0,\dots,2}$ defined by
\[
 \vec{v}_{3k+i}:=\theta_k\mathbf{e}_i,
\]
where $\{\mathbf{e}_i\}_{i=0,\dots,2}$ is the canonical basis of $\mathbb{R}^3$.
We observe that, given $k\in \mathbb{N}$, the elements $\vec{v}_{3k+i}$ span the $3$-th dimensional linear eigenspace of vector fields associated to the eigenvalue $\rho_k$. 
We then introduce the projection operator
\[
PV_m:\pier{H_{\Gamma_D}^1(\mathcal{D}_a;\mathbb{R}^{3\times 3})}\to \text{span}\{\vec{v}_0,\vec{v}_1,\dots,\vec{v}_{3m+2}\}.
\]
We make the Galerkin ansatz
\begin{align}
\label{galerkinansatz}
& \displaystyle \pier{\mathbf{W}_m(\vec{a}, t)=\mathbf{I}+\sum_{i=0}^{6m+5}x_i^m(t)\mathbf{S}_i(\vec{a}), \quad 
\boldsymbol{\Theta}_m(\vec{a}, t)=\sum_{i=0}^{3m+2}y_i^m(t)\mathbf{A}_i(\vec{a})},
\end{align}
for $(\vec{a}, t) \in \mathcal{D}_a \times (0,T)$, with 
\begin{align*}
&\mathbf{S}_i\in \mathcal{S}\cap \bar{H}^2_{\Gamma_D}(\mathcal{D}_a,\mathbb{R}^{3\times 3}), \\
&\mathbf{A}_i\in \mathcal{A}\cap \bar{H}^2_{\Gamma_D}(\mathcal{D}_a,\mathbb{R}^{3\times 3}),
\end{align*}\
to approximate the solutions $\mathbf{W}$, $\boldsymbol{\Theta}$ of \pier{the system \eqref{symskewreg}}. 
\abramon{
Moreover, we approximate the solution $\vec{u}$ in \eqref{symskewreg}$_4$ as
\abramonn{
\begin{equation}
\label{galerkinansatz2}
{\vec{u}}_m(\vec{a}, t)=\sum_{i=0}^{3m+2}z_i^m(t)\vec{v}_i(\vec{a}),
\end{equation}
}
with 
\[
\vec{v}_i\in H_{\Gamma_D}^1(\mathcal{D}_a;\mathbb{R}^{3}).
\]
}
Given \eqref{galerkinansatz}, we define
\begin{align}
 \label{eqn:sigmammapp}
&{}
\boldsymbol{\Sigma}_m=\partial \psi_D^{\lambda}\left(\curl \left(\mathbf{Z}_m\right)\right),
\\
 \label{eqn:phimzmapp}
 &
 \begin{cases}
- P_L\Delta {\mathbf{Z}_m}=\curl \left(\mathrm{e}^{\boldsymbol{\Theta_m}}\mathbf{W}_m\right),\\
\pier{\ \mathbf{Z}_m =\mathbf{0} \quad \text{on} \; \Gamma_N\times (0,T)},\;\; \curl \mathbf{Z}_m \wedge \vec{n}=\mathbf{0}\quad \text{on} \; \Gamma_D\times (0,T).
\end{cases}
\end{align}
Given the elliptic \pier{problem in \eqref{eqn:phimzmapp} with approximated right hand sides, we then have
\begin{align}
 \label{phimzm}
&\mathbf{Z}_m(\vec{a}, t)=\mathcal{G}_{L,\diver}\ast \curl\left(\mathrm{e}^{\boldsymbol{\Theta}_m}\mathbf{W}_m\right)(\vec{a}, t)
\;\; \text{for} \;  (\vec{a}, t) \in \mathcal{D}_a\times (0,T), \;\; \mathbf{Z}_m\in H_{\Gamma_N}^1(\mathcal{D}_a,\mathbb{R}^{3\times 3})\cap \mathcal{M}_{div}.
\end{align}
}
\abramonn{Taking in \eqref{symskewreg} $\nwhat{\mathbf{W}}=\mathbf{S}_i$, $\nwhat{\boldsymbol{\Omega}}=\mathbf{A}_j$, $\vec{v}=\vec{v}_l$, with $i=0, \dots, 6m+5$, $j,l=0, \dots, 3m+2$, and considering the time derivative of \eqref{symskewreg}$_4$, we obtain the following Galerkin approximation of System \eqref{symskewreg}:
\begin{equation}
\label{eqn:48m}
\begin{cases}
{\int_{\mathcal{D}_a}\mathrm{e}^{-\boldsymbol{\Theta}_m}\grad \left(\mathcal{G}_L\ast \twobigdot{\vec{u}}_m\right)\colon \mathbf{S}_i}+ \int_{\mathcal{D}_a}\twobigdot{\mathbf{W}}_m\colon \mathbf{S}_i\\
{}+\int_{\mathcal{D}_a}{\grad{}} \twobigdot{\mathbf{W}}_m::{\grad{}} \mathbf{S}_i+\int_{\mathcal{D}_a}{\grad\grad{}} \twobigdot{\mathbf{W}}_m\colon {\grad\grad{}} \mathbf{S}_i\\
{}+\int_{\mathcal{D}_a}\mathrm{e}^{-\boldsymbol{\Theta}_m}\curl\left(
\mathcal{G}_{L,\diver}\ast \curl \left[\partial \psi_D^{\lambda}\left(\curl \left(\mathcal{G}_{L,\diver}\ast \curl\left(\mathrm{e}^{\boldsymbol{\Theta}_m}\mathbf{W}_m\right)\right)\right)\right]\right)\colon \mathbf{S}_i\\
{}+ \int_{\mathcal{D}_a}\left(\mathbf{W}_m-\mathbf{I}+\frac{{d{}\tilde{\psi}_{C_\alpha}}}{d \mathbf{W}_m}(\mathbf{W}_m)+\overbigdot{\mathbf{W}}_m\right)\colon \mathbf{S}_i+\int_{\mathcal{D}_a}{\grad{}} (\mathbf{W}_m+\overbigdot{\mathbf{W}_m})::{\grad{}} \mathbf{S}_i
\\
{}+\int_{\mathcal{D}_a}{\grad\grad{}} (\mathbf{W}_m+\overbigdot{\mathbf{W}_m})\colon {\grad\grad{}} \mathbf{S}_i=\int_{\mathcal{D}_a}\mathbf{W}_{\text{ext}}(\mathbf{W}_m,t)\colon \mathbf{S}_i\\
{}+\int_{\mathcal{D}_a}\mathrm{e}^{-\boldsymbol{\Theta}_m}\grad \left(\mathcal{G}_L\ast \vec{\mathcal{F}}_{\textrm ext}\right)\colon \mathbf{S}_i,\\ \\
\int_{\mathcal{D}_a}\grad \left(\mathcal{G}_L\ast \twobigdot{\vec{u}}_m\right)\mathbf{W}_m\mathrm{e}^{-\boldsymbol{\Theta}_m}\colon \mathbf{A}_j\\
{}+\int_{\mathcal{D}_a}\twobigdot{\boldsymbol{\Theta}}_m\colon \mathbf{A}_j+\int_{\mathcal{D}_a}\grad{} \twobigdot{\boldsymbol{\Theta}}_m::{\grad{}} \mathbf{A}_j+\int_{\mathcal{D}_a}{\grad\grad{}} \twobigdot{\boldsymbol{\Theta}}_m\colon {\grad\grad{}} \mathbf{A}_j\\
{}+\int_{\mathcal{D}_a}\left(\curl\left(
\mathcal{G}_{L,\diver}\ast \curl \left[\partial \psi_D^{\lambda}\left(\curl \left(\mathcal{G}_{L,\diver}\ast \curl\left(\mathrm{e}^{\boldsymbol{\Theta}_m}\mathbf{W}_m\right)\right)\right)\right]\right)\mathbf{W}_m\mathrm{e}^{-\boldsymbol{\Theta}_m}\right)\colon \mathbf{A}_j \\
{}+\frac{1}{2}\int_{\mathcal{D}_a}{\grad{}} (\boldsymbol{\Theta}_m+\overbigdot{\boldsymbol{\Theta}}_m)::{\grad{}} \mathbf{A}_j+\frac{1}{2}\int_{\mathcal{D}_a}{\grad\grad{}}(\boldsymbol{\Theta}+ \overbigdot{\boldsymbol{\Theta}}_m)\colon{\grad\grad{}} \mathbf{A}_j\\
=\int_{\mathcal{D}_a}\boldsymbol{\Omega}_{\text{ext}}{(\boldsymbol{\Theta}_m,t)}\colon \mathbf{A}_j+\int_{\mathcal{D}_a}\grad \left(\mathcal{G}_L\ast \vec{\mathcal{F}}_{\textrm ext}\right)\mathbf{W}_m\mathrm{e}^{-\boldsymbol{\Theta}_m}\colon \mathbf{A}_j,\\ \\
\int_{\mathcal{D}_a}\twobigdot{{\vec{u}}}_m\cdot \vec{v}_l+\frac{1}{\lambda}\int_{\mathcal{D}_a}\grad{}\overbigdot{{\vec{u}}}_m\colon \grad{}\vec{v}_l=\frac{1}{\lambda}\int_{\mathcal{D}_a}\left(\overbigdot{\boldsymbol{\Theta}}_m\mathrm{e}^{\boldsymbol{\Theta}_m}\mathbf{W}_m+\mathrm{e}^{\boldsymbol{\Theta}_m}\overbigdot{\mathbf{W}}_m\right)\colon \grad{}\vec{v}_l,\\ \\
- P_L\Delta {\mathbf{Z}_m}=\curl \left(\mathrm{e}^{\boldsymbol{\Theta_m}}\mathbf{W}_m\right),
\end{cases}
\end{equation}
}{in $[0,t]$}, with $0<t\leq T$, with boundary conditions as in \eqref{eqn:34} and with initial conditions 
\begin{equation}
    \label{eqn:48mic}
  { \mathbf{W}_m(\vec{a},0)=\mathbf{I}, \quad \overbigdot{ \mathbf{W}}_m(\vec{a},0)=\boldsymbol{\Theta}_m(\vec{a},0)=\overbigdot{\boldsymbol{\Theta}}_m(\vec{a},0)=\boldsymbol{0}, \quad \vec{u}_m(\vec{a},0)=\overbigdot{\vec{u}}_m(\vec{a},0)=\vec{0}, \quad \vec{a}\in \mathcal{D}_a.}
\end{equation}
\abramonn{
We observe from \eqref{eqn:48m}$_4$ that
\begin{equation}
\label{phimdd}
\twobigdot{{\vec{u}}}_m(t)=-\frac{1}{\lambda}\sum_l\rho_l\overbigdot{z}_l^m\vec{v}_l+\frac{1}{\lambda}\sum_l\left(\int_{\mathcal{D}_a}\left(\overbigdot{\boldsymbol{\Theta}}_m\mathrm{e}^{\boldsymbol{\Theta}_m}\mathbf{W}_m+\mathrm{e}^{\boldsymbol{\Theta}_m}\overbigdot{\mathbf{W}}_m\right):\grad \vec{v}_l\right)\vec{v}_l.
\end{equation}
Substituting \eqref{phimdd} in \eqref{eqn:48m}$_1$ and  \eqref{eqn:48m}$_2$, we obtain that system
\eqref{eqn:48m} defines the following collection of initial value problems for a system of coupled second order ODEs,
\begin{equation}
\label{48tempm}
\begin{cases}
\displaystyle (1+\gamma_i)\left(\twobigdot{x_i^m}+\overbigdot{x_i^m}+x_i^m\right)=x_i^m+\int_{\mathcal{D}_a}\biggl(-\frac{{d{}\tilde{\psi}_{C_\alpha}}}{d \mathbf{W}_m}\biggl({{}\mathbf{I} +{}}\sum_lx_l^m\mathbf{S}_l\biggr)+\mathbf{W}_{\text{ext}}\biggl({{}\mathbf{I} +{}} \sum_lx_l^m\mathbf{S}_l,t\biggr)\biggr)\colon \mathbf{S}_i\\
\displaystyle +\int_{\mathcal{D}_a}PS_m\biggl[\mathrm{e}^{-\sum_ly_l^m\mathbf{A}_l}\grad \left(\mathcal{G}_L\ast \vec{\mathcal{F}}_{\textrm ext}\right)\biggr]\colon \mathbf{S}_i- \int_{\mathcal{D}_a}PS_m\biggl[\mathrm{e}^{-\sum_ly_l^m\mathbf{A}_l}\curl\biggl(\mathcal{G}_{L,\diver}\ast \\
\displaystyle 
\hskip2cm \displaystyle \curl \biggl[\partial \psi_D^{\lambda}\biggl(\curl \biggl(
\mathcal{G}_{L,\diver}\ast \curl\biggl(\mathrm{e}^{\sum_ry_r^m\mathbf{A}_r}\biggl({{}\mathbf{I} +{}}\sum_kx_k^m\mathbf{S}_k\biggr)\biggr)\biggr)\biggr)\biggr]\biggr)\biggr]\colon \mathbf{S}_i\\
\displaystyle - \int_{\mathcal{D}_a}PS_m\biggl[\mathrm{e}^{-\sum_ly_l^m\mathbf{A}_l}\grad \biggl(\mathcal{G}_L\ast \biggl[-\frac{1}{\lambda}\sum_l\rho_l\overbigdot{z}_l^m\vec{v}_l+\frac{1}{\lambda}\sum_l\biggl(\int_{\mathcal{D}_a}\biggl(\sum_l\overbigdot{y_l^m}\mathbf{A}_r\mathrm{e}^{\sum_ry_r^m\mathbf{A}_r}\left(\mathbf{I}+\sum_sx_s^m\mathbf{S}_s\right)\\
\displaystyle
\hskip2cm +\mathrm{e}^{\sum_ly_l^m\mathbf{A}_l}\sum_r\overbigdot{x_r^m}\mathbf{S}_r\biggr):\grad \vec{v}_l\biggr)\vec{v}_l\biggr]\biggr)\biggr]\colon \mathbf{S}_i,
\\ \\
\displaystyle (1+\gamma_j)\twobigdot{y_j^m}=-\gamma_j\left(\frac{1}{2}\overbigdot{y_j^m}+\frac{1}{2}y_j^m\right)+\int_{\mathcal{D}_a}\boldsymbol{\Omega}_{\text{ext}}{(\sum_ly_l^m\mathbf{A}_l)}\colon \mathbf{A}_j\\
\displaystyle+\int_{\mathcal{D}_a}PA_m\biggl[\grad \left(\mathcal{G}_L\ast \vec{\mathcal{F}}_{\textrm ext}\right)\left(\mathbf{I}+\sum_lx_l^m\mathbf{S}_l\right)\mathrm{e}^{\sum_ry_r^m\mathbf{A}_r}\biggr]\colon \mathbf{A}_j\\
\displaystyle -\int_{\mathcal{D}_a}PA_m\biggl[\curl \biggl(\mathcal{G}_{L,\diver}\ast \curl \biggl[\partial \psi_D^{\lambda}\biggl(\curl \biggl(
\mathcal{G}_{L,\diver}\ast\\
\hskip2cm \displaystyle  \curl\biggl(\mathrm{e}^{\sum_ry_r^m\mathbf{A}_r}\biggl({{}\mathbf{I} +{}}\sum_kx_k^m\mathbf{S}_k\biggr)
\biggr)\!\biggr)\!\biggr)\!\biggr]\!\biggr) \biggl({{}\mathbf{I} +{}}\sum_kx_k^m\mathbf{S}_k\biggr)\mathrm{e}^{-\sum_py_p^m\mathbf{A}_p}\biggr]\colon \mathbf{A}_j\\
\displaystyle - \int_{\mathcal{D}_a}PA_m\biggl[\grad \biggl(\mathcal{G}_L\ast \biggl[-\frac{1}{\lambda}\sum_l\rho_l\overbigdot{z}_l^m\vec{v}_l+\frac{1}{\lambda}\sum_l\biggl(\int_{\mathcal{D}_a}\biggl(\sum_l\overbigdot{y_l^m}\mathbf{A}_r\mathrm{e}^{\sum_ry_r^m\mathbf{A}_r}\left(\mathbf{I}+\sum_sx_s^m\mathbf{S}_s\right)\\
\displaystyle
\hskip2cm +\mathrm{e}^{\sum_ly_l^m\mathbf{A}_l}\sum_r\overbigdot{x_r^m}\mathbf{S}_r\biggr):\grad \vec{v}_l\biggr)\vec{v}_l
\biggr]\biggr)\left(\mathbf{I}+\sum_lx_l^m\mathbf{S}_l\right)\mathrm{e}^{\sum_ry_r^m\mathbf{A}_r}\biggr]\colon \mathbf{A}_j,\\ \\
\displaystyle \twobigdot{z_l^m}=-\frac{\rho_l}{\lambda}\overbigdot{z_l^m}+\frac{1}{\lambda}\int_{\mathcal{D}_a}PV_m\left[\sum_l\overbigdot{y_k^m}\mathbf{A}_k\mathrm{e}^{\sum_ry_r^m\mathbf{A}_r}\left(\mathbf{I}+\sum_sx_s^m\mathbf{S}_s\right)+\mathrm{e}^{\sum_ky_k^m\mathbf{A}_k}\sum_r\overbigdot{x_r^m}\mathbf{S}_r\right]:\grad \vec{v}_l,\\
\displaystyle x_i^m(0)=\overbigdot{x}_i^m(0)=0, \, \  y_j^m(0)=\overbigdot{y}_j^m(0)=0, \, \  z_l^m(0)=\overbigdot{z}_l^m(0)=0,
\end{cases}
\end{equation}
for $i=0, \dots, 6m+5, \quad j,l=0, \dots, 3m+2$.
We conclude that \eqref{48tempm} is of the form}
\begin{equation}
\begin{cases}
\label{eqn:48ode}
\displaystyle \twobigdot{\mathbf{x}}=f\left(\mathbf{x},\mathbf{y},\mathbf{z},\overbigdot{\mathbf{x}},\overbigdot{\mathbf{y}},\overbigdot{\mathbf{z}},t\right),\\
\displaystyle \twobigdot{\mathbf{y}}=g\left(\mathbf{x},\mathbf{y},\mathbf{z},\overbigdot{\mathbf{x}},\overbigdot{\mathbf{y}},\overbigdot{\mathbf{z}},t\right),\\
\displaystyle \twobigdot{\mathbf{z}}=h\left(\mathbf{x},\mathbf{y},\mathbf{z},\overbigdot{\mathbf{x}},\overbigdot{\mathbf{y}},\overbigdot{\mathbf{z}},t\right),\\
\displaystyle \mathbf{x}(0)=\overbigdot{\mathbf{x}}(0)=\mathbf{y}(0)=\overbigdot{\mathbf{y}}(0)=\mathbf{z}(0)=\overbigdot{\mathbf{z}}(0)=\mathbf{0}.
\end{cases}
\end{equation}
\begin{rem}
We observe that
\[
\mathbf{W}_m\in \mathring{C}_{\alpha} \iff \mathbf{x}^m\in \mathring{C}_{\alpha}^m,
\]
where $\mathring{C}_{\alpha}^m\subset \mathbb{R}^{3m+5}$ is an open neighborhood of $\mathbf{x}^m=\mathbf{0}$.
\end{rem}
Due to Assumptions \textbf{A3}, to the smoothness of $\tilde{\psi}_{C_{\alpha}}$, to the Lipschitz \pier{continuity of $\partial \psi_D^{\lambda}$} and to the regularity in space of the functions $\mathbf{S}_i,  \mathbf{A}_j, \vec{v}_l$, \pier{the system}~\eqref{eqn:48ode} is a coupled system of \pier{second-order} ODEs in the variables $x_i^m, \, y_j^m, \, z_l^m$, with a right hand side which is measurable in time and continuous in the independent variables. In particular, let us observe, thanks to \eqref{phimdd} and the regularity of $\vec{\mathcal{F}}_{\textrm ext}$, to the fact that the elements of the  images of $\mathcal{G}_{L}$ and $\mathcal{G}_{L,\diver}$ have $H^1$ regularity and thanks to the regularity in space of the functions $\mathbf{S}_i,  \mathbf{A}_j, \vec{v}_l$, that all integrands in the nonlinearly coupled terms in the right hand side of \eqref{eqn:48m} are integrable.
Then,  we can apply the Carath\'{e}odory's existence theorem to infer that there exist a sufficiently small $T_m$ with $0<T_m\leq T$ and a local solution $(x_i^m,y_j^m,z_l^m)$ of \eqref{eqn:48ode}, for $i=0, \dots, 6m+5$, $j,l=0, \dots, 3m+2$, such that 
\[
\mathbf{x}^m,\mathbf{y}^m,\mathbf{z}^m\in W^{2,\infty}([0,T_m]), \quad \text{and} \; \mathbf{x}^m\in \mathring{C}_{\alpha}^m.
\]
\pier{Once we have} a solution to \eqref{eqn:48ode}, 
dealing with \pier{the} elliptic problems with regular right-hand sides in \pier{\eqref{eqn:48m} leads to the elements $\mathbf{Z}_m$ solving} \eqref{eqn:48m}$_4$. 

\pier{Next,
thanks to some uniform} estimates, we will extend these solutions by continuity to the interval $[0,\hat{T}]$, with $\hat{T}$ independent on the discretization and regularization parameters, and we will study the limit as $m\to \infty$ and $\lambda\to 0$. 

\subsection{A priori estimates}
We now deduce a priori estimates, uniform in the discretization parameter $m$ and in the regularization parameters $\lambda$, for the solutions of system \eqref{eqn:48m}, which can be rewritten, combining the equations over $i=0, \dots, 6m+5$, $j=0, \dots, 3m+2$ and $l=0, \dots, 3m+2$, as
\begin{equation}
\label{eqn:48m2}
\begin{cases}
\int_{\mathcal{D}_a}\mathrm{e}^{-\boldsymbol{\Theta}_m}\grad \left(\mathcal{G}_L\ast \twobigdot{\vec{u}}_m\right)\colon \nwhat{\mathbf{W}}_m\\
{}+ \int_{\mathcal{D}_a}\twobigdot{\mathbf{W}}_m\colon \nwhat{\mathbf{W}}_m+\int_{\mathcal{D}_a}\pier{\grad{}} \twobigdot{\mathbf{W}}_m::\pier{\grad{}} \nwhat{\mathbf{W}}_m+\int_{\mathcal{D}_a}\pier{\grad\grad{}} \twobigdot{\mathbf{W}}_m\colon \pier{\grad\grad{}} \nwhat{\mathbf{W}}_m\\
{}+\int_{\mathcal{D}_a}\mathrm{e}^{-\boldsymbol{\Theta}_m}\curl\left(
\mathcal{G}_{L,\diver}\ast \curl \boldsymbol{\Sigma}_m\right)\colon \nwhat{\mathbf{W}}_m\\
{}+ \int_{\mathcal{D}_a}\left(\mathbf{W}_m-\mathbf{I}+\frac{\pier{{}\tilde{\psi}_{C_\alpha}}}{d \mathbf{W}_m}(\mathbf{W}_m)+\overbigdot{\mathbf{W}}_m\right)\colon \nwhat{\mathbf{W}}_m+\int_{\mathcal{D}_a}\pier{\grad{}} (\mathbf{W}_m+\overbigdot{\mathbf{W}_m})::\pier{\grad{}} \nwhat{\mathbf{W}}_m
\\
{}+\int_{\mathcal{D}_a}\pier{\grad\grad{}} (\mathbf{W}_m+\overbigdot{\mathbf{W}_m})\colon \pier{\grad\grad{}}\nwhat{\mathbf{W}}_m=\int_{\mathcal{D}_a}\mathbf{W}_{\text{ext}}(\mathbf{W}_m,t)\colon \nwhat{\mathbf{W}}_m\\
{}+\int_{\mathcal{D}_a}\mathrm{e}^{-\boldsymbol{\Theta}_m}\grad \left(\mathcal{G}_L\ast \vec{\mathcal{F}}_{\textrm ext}\right)\colon \nwhat{\mathbf{W}}_m,\\ \\
\int_{\mathcal{D}_a}\grad \left(\mathcal{G}_L\ast \twobigdot{\vec{u}}_m\right)\mathbf{W}_m\mathrm{e}^{-\boldsymbol{\Theta}_m}\colon \nwhat{\boldsymbol{\Omega}}_m\\
{}+\int_{\mathcal{D}_a}\twobigdot{\boldsymbol{\Theta}}_m\colon \nwhat{\boldsymbol{\Omega}}_m+\int_{\mathcal{D}_a}\grad{} \twobigdot{\boldsymbol{\Theta}}_m::\pier{\grad{}} \nwhat{\boldsymbol{\Omega}}_m+\int_{\mathcal{D}_a}\pier{\grad\grad{}} \twobigdot{\boldsymbol{\Theta}}_m\colon \pier{\grad\grad{}} \nwhat{\boldsymbol{\Omega}}_m\\
{}+\int_{\mathcal{D}_a}\left(\curl\left(
\mathcal{G}_{L,\diver}\ast \curl \boldsymbol{\Sigma}_m\right)\mathbf{W}_m\mathrm{e}^{-\boldsymbol{\Theta}_m}\right)\colon \nwhat{\boldsymbol{\Omega}}_m \\
{}+\frac{1}{2}\int_{\mathcal{D}_a}\pier{\grad{}} (\boldsymbol{\Theta}_m+\overbigdot{\boldsymbol{\Theta}}_m)::\pier{\grad{}} \nwhat{\boldsymbol{\Omega}}_m+\frac{1}{2}\int_{\mathcal{D}_a}\pier{\grad\grad{}}(\boldsymbol{\Theta}+ \overbigdot{\boldsymbol{\Theta}}_m)\colon\pier{\grad\grad{}} \nwhat{\boldsymbol{\Omega}}_m\\
=\int_{\mathcal{D}_a}\boldsymbol{\Omega}_{\text{ext}}\abramo{(\boldsymbol{\Theta}_m,t)}\colon \nwhat{\boldsymbol{\Omega}}_m+\int_{\mathcal{D}_a}\grad \left(\mathcal{G}_L\ast \vec{\mathcal{F}}_{\textrm ext}\right)\mathbf{W}_m\mathrm{e}^{-\boldsymbol{\Theta}_m}\colon \nwhat{\boldsymbol{\Omega}}_m,\\ \\
\boldsymbol{\Sigma}_m=\partial \psi_D^{\lambda}(\curl \mathbf{Z}_m),\\ \\
\lambda \int_{\mathcal{D}_a}\overbigdot{{\vec{u}}}_m\cdot \vec{v}_m+\int_{\mathcal{D}_a}\grad{}{{\vec{u}}}_m\colon \grad{}\vec{v}_m=\int_{\mathcal{D}_a}\left(\mathrm{e}^{\boldsymbol{\Theta}_m}\mathbf{W}_m-\mathbf{I}\right)\colon \grad{}\vec{v}_m,\\ \\
- P_L\Delta {\mathbf{Z}_m}=\curl \left(\mathrm{e}^{\boldsymbol{\Theta_m}}\mathbf{W}_m\right),
\end{cases}
\end{equation}
for a.e. $t \in [0,T_m]$ \pier{and all} $\nwhat{\mathbf{W}}_m \in \text{span}\left\{\mathbf{S}_0,\mathbf{S}_1,\dots,\mathbf{S}_{6m+5}\right\}$, $\nwhat{\boldsymbol{\Omega}}_m \in \text{span}\left\{\mathbf{A}_0,\mathbf{A}_1,\dots,\mathbf{A}_{3m+2}\right\}$,  $\vec{v}_m \in \text{span}\left\{\vec{v}_0,\vec{v}_1,\dots,\vec{v}_{3m+2}\right\}$ and with initial conditions defined in \eqref{eqn:48mic}. 

\subsubsection{First a priori estimate}
The first a-priori estimate is obtained by taking $\nwhat{\mathbf{W}}_m=\overbigdot{\mathbf{W}}_m$ in \eqref{eqn:48m2}$_1$ and $\nwhat{\boldsymbol{\Omega}}_m=\overbigdot{\boldsymbol{\Theta}}_m$ in \eqref{eqn:48m2}$_2$. Moreover, we take the time derivative of \eqref{eqn:48m2}$_4$ and $\vec{v}_m=\mathcal{G}_L\ast\twobigdot{\vec{u}}_m$. We further take the time derivative of \eqref{eqn:48m2}$_5$, multiply it by $\mathcal{G}_{L,\diver}\ast\left( \curl \boldsymbol{\Sigma}_m\right)$ and integrate over $\mathcal{D}_a$. Finally, we sum all the contributions.
\newline
We observe from \eqref{eqn:48m2}$_5$ and from the regularity in space of the functions $\mathbf{S}_i, \mathbf{A}_j$ that $\mathbf{Z}_m\in H^1(\mathcal{D}_a;\mathbb{R}^{3\times 3})$, for any $t\in [0,T_m]$. Hence, from \eqref{eqn:48m2}$_3$ and the Lipschitz continuity of $\partial \psi_D^{\lambda}$ we obtain that $\boldsymbol{\Sigma}_m\in L^2(\mathcal{D}_a;\mathbb{R}^{3\times 3})$ for any $t\in [0,T_m]$, and as a consequence the dual  product of equation~\pier{\eqref{eqn:48m2}$_5$} with the element $\mathcal{G}_{L,\diver}\ast\left( \curl \boldsymbol{\Sigma}_m\right)\in H^1(\mathcal{D}_a;\mathbb{R}^{3\times 3})$ is well defined for any $t\in [0,T_m]$.
We observe that
\begin{align*}
&\lambda \int_{\mathcal{D}_a}\twobigdot{{\vec{u}}}_m\cdot \mathcal{G}_L\ast\twobigdot{\vec{u}}_m+\int_{\mathcal{D}_a}\grad{}{{\overbigdot{{\vec{u}}}}}_m\colon \grad{}{\mathcal{G}_L\ast\twobigdot{\vec{u}}}_m=\lambda ||\twobigdot{{\vec{u}}}_m||_{\left(H_{\Gamma_D}^1(\mathcal{D}_a,\mathbb{R}^3)\right)'}+\frac{1}{2}\frac{d}{dt}||\overbigdot{{\vec{u}}}_m||^2.
\end{align*}
Also, the contribution from \pier{\eqref{eqn:48m2}$_5$}, after integration by parts and considering the boundary conditions and \eqref{pigamma}, gives that
\begin{align*}
& \prescript{}{\left(H_{\Gamma_N}^{1}(\mathcal{D}_a,\mathbb{R}^{3\times 3})\right)'}{<}-\Delta \overbigdot{\mathbf{Z}}_m,\mathcal{G}_{L,\diver}\ast \left(\curl\boldsymbol{\Sigma}_m\right)>_{H_{\Gamma_N}^{1}(\mathcal{D}_a,\mathbb{R}^{3\times 3})}\\
& \quad=\prescript{}{\left(H_{\Gamma_N}^{1}(\mathcal{D}_a,\mathbb{R}^{3\times 3})\right)'}{<}\curl\boldsymbol{\Sigma}_m,\overbigdot{\mathbf{Z}}_m>_{H_{\Gamma_N}^{1}(\mathcal{D}_a,\mathbb{R}^{3\times 3})}\\
& \quad = \int_{\mathcal{D}_a} \curl \overbigdot{\mathbf{Z}}_m\colon \partial \psi_D^{\lambda}(\curl \mathbf{Z}_m).
\end{align*}
We then obtain that
\begin{align}
\label{eqn:53t0}
& \displaystyle \notag \frac{d}{dt}\biggl(\frac{1}{2}\pier{\|}\overbigdot{\mathbf{W}}_m\pier{\|}^2+\frac{1}{2}\pier{\|}\grad{} \overbigdot{\mathbf{W}}_m\pier{\|}^2+\frac{1}{2}\pier{\|}\grad\grad \overbigdot{\mathbf{W}}_m\pier{\|}^2+\frac{1}{2}\pier{\|}\mathbf{W}_m-\mathbf{I}\pier{\|}^2
+\pier{\tilde{\psi}_{C_{\alpha}}}(\mathbf{W}_m)\\
& \displaystyle \notag +\frac{1}{2}\pier{\|}\grad \mathbf{W}_m\pier{\|}^2+\frac{1}{2}\pier{\|}\grad\grad \mathbf{W}_m\pier{\|}^2+\frac{1}{2}\pier{\|} \overbigdot{\boldsymbol{\Theta}}_m\pier{\|}^2+\frac{1}{2}\pier{\|}\grad \overbigdot{\boldsymbol{\Theta}}_m\pier{\|}^2+\frac{1}{2}\pier{\|}\grad\grad \overbigdot{\boldsymbol{\Theta}}_m\pier{\|}^2\\
& \displaystyle \notag +\frac{1}{4}\pier{\|} {\boldsymbol{\Theta}}_m\pier{\|}^2+\frac{1}{4}\pier{\|}\grad \boldsymbol{\Theta}_m\pier{\|}^2+\frac{1}{4}\pier{\|}\grad\grad \boldsymbol{\Theta}_m\pier{\|}^2+\frac{1}{2}\pier{\|}\overbigdot{{\vec{u}}}_m\pier{\|}^2+\psi_D^{\lambda}(\curl \mathbf{Z}_m)\biggr)\\
& \displaystyle \notag +\pier{\|}\overbigdot{\mathbf{W}}_m\pier{\|}^2+\pier{\|}\grad \overbigdot{\mathbf{W}}_m\pier{\|}^2+\pier{\|}\grad\grad \overbigdot{\mathbf{W}}_m\pier{\|}^2+\frac{1}{2}\pier{\|} \grad\overbigdot{\boldsymbol{\Theta}}_m\pier{\|}^2+\frac{1}{2}\pier{\|} \grad\grad\overbigdot{\boldsymbol{\Theta}}_m\pier{\|}^2\\
& \displaystyle \notag +\lambda\pier{\|}\twobigdot{\vec{u}}_m\pier{\|}_{\left(H_{\Gamma_D}^1(\mathcal{D}_a,\mathbb{R}^3)\right)'}^2=\int_{\mathcal{D}_a}\mathbf{W}_{\text{ext}}(\mathbf{W}_m,t)\colon \overbigdot{\mathbf{W}}_m+\int_{\mathcal{D}_a}\boldsymbol{\Omega}_{\text{ext}}\abramo{(\boldsymbol{\Theta}_m,t)}\colon \overbigdot{\boldsymbol{\Theta}}_m+\frac{1}{2}\int_{\mathcal{D}_{a}} {\boldsymbol{\Theta}}_m\colon \overbigdot{\boldsymbol{\Theta}}_m\\
& \displaystyle \int_{\mathcal{D}_a}\mathrm{e}^{-\boldsymbol{\Theta}_m}\grad \left(\mathcal{G}_L\ast \vec{\mathcal{F}}_{\textrm ext}\right)\colon \overbigdot{\boldsymbol{\Theta}}_m+\int_{\mathcal{D}_a}\grad \left(\mathcal{G}_L\ast \vec{\mathcal{F}}_{\textrm ext}\right)\mathbf{W}_m\mathrm{e}^{-\boldsymbol{\Theta}_m}\colon \overbigdot{\boldsymbol{\Theta}}_m,
\end{align}
where we added $\frac{1}{4}\frac{d}{dt}\pier{\|} \boldsymbol{\Theta}_m\pier{\|}^2$ to the left and $\frac{1}{2}\int_{\mathcal{D}_{a}}\overbigdot{\boldsymbol{\Theta}}_m\colon {\boldsymbol{\Theta}}_m$ to the right hand side.
We now observe that, thanks to Assumption \textbf{A3}, we have that
\begin{equation}
\label{ap1temp1}
\grad \left(\mathcal{G}_L\ast \vec{\mathcal{F}}_{\textrm ext}\right)\in L^{\infty}\left(0,T_m;L^2(\mathcal{D}_a,\mathbb{R}^{3\times 3})\right),
\end{equation}
and, using a trilinear H\"{o}lder inequality, \eqref{eqn:2}, the Cauchy--Schwarz and the Young inequalities, that
\begin{align}
\label{ap1temp2}
& \displaystyle \notag\int_{\mathcal{D}_{aT_m}}\grad \left(\mathcal{G}_L\ast \vec{\mathcal{F}}_{\textrm ext}\right)\mathbf{W}_m\mathrm{e}^{-\boldsymbol{\Theta}_m}\colon \overbigdot{\boldsymbol{\Theta}}_m\\
& \displaystyle  \notag \quad \leq \int_0^{T_m}\pier{\|}\grad \left(\mathcal{G}_L\ast \vec{\mathcal{F}}_{\textrm ext}\right)\pier{\|}\,\pier{\|}\mathbf{W}_m\pier{\|}_{L^3(\mathcal{D}_a,Sym(\mathbb{R}^{3\times 3}))}\pier{\|}\mathrm{e}^{-\boldsymbol{\Theta}_m}\pier{\|}_{L^{\infty}(\mathcal{D}_a,\mathbb{R}^{3\times 3})}\pier{\|}\overbigdot{\boldsymbol{\Theta}}_m\pier{\|}_{L^6(\mathcal{D}_a,Skew(\mathbb{R}^{3\times 3}))}\\
& \displaystyle \notag \quad \leq C\int_0^{T_m}\pier{\|}\mathbf{W}_m\pier{\|}^{1/2}\pier{\|}\mathbf{W}_m\pier{\|}_{H^1(\mathcal{D}_a,Sym(\mathbb{R}^{3\times 3}))}^{1/2}\pier{\|}\overbigdot{\boldsymbol{\Theta}}_m\pier{\|}_{H^1(\mathcal{D}_a,Skew(\mathbb{R}^{3\times 3}))}\\
& \displaystyle \quad \leq C\int_0^{T_m}\left(\pier{\|}\mathbf{W}_m\pier{\|}^2+\pier{\|}\overbigdot{\boldsymbol{\Theta}}_m\pier{\|}^2+\pier{\|}\grad \mathbf{W}_m\pier{\|}^2+\pier{\|}\grad \overbigdot{\boldsymbol{\Theta}}_m\pier{\|}^2\right).
\end{align}
Using in \eqref{eqn:53t0} the results \eqref{ap1temp1}, \eqref{ap1temp2}, Assumptions \textbf{A2} and \textbf{A3}, integrating in time over the interval $[0,T_m]$, we obtain that
\begin{align}
\label{eqn:53}
& \displaystyle \notag \frac{1}{2}\pier{\|}\overbigdot{\mathbf{W}}_m\pier{\|}^2+\frac{1}{2}\pier{\|}\grad{} \overbigdot{\mathbf{W}}_m\pier{\|}^2+\frac{1}{2}\pier{\|}\grad\grad \overbigdot{\mathbf{W}}_m\pier{\|}^2+\frac{1}{2}\pier{\|}\mathbf{W}_m-\mathbf{I}\pier{\|}^2
+\pier{\tilde{\psi}_{C_{\alpha}}}(\mathbf{W}_m)\\
& \displaystyle \notag +\frac{1}{2}\pier{\|}\grad \mathbf{W}_m\pier{\|}^2+\frac{1}{2}\pier{\|}\grad\grad \mathbf{W}_m\pier{\|}^2+\frac{1}{2}\pier{\|} \overbigdot{\boldsymbol{\Theta}}_m\pier{\|}^2+\frac{1}{2}\pier{\|}\grad \overbigdot{\boldsymbol{\Theta}}_m\pier{\|}^2+\frac{1}{2}\pier{\|}\grad\grad \overbigdot{\boldsymbol{\Theta}}_m\pier{\|}^2\\
& \displaystyle \notag +\frac{1}{4}\pier{\|} {\boldsymbol{\Theta}}_m\pier{\|}^2+\frac{1}{4}\pier{\|}\grad \boldsymbol{\Theta}_m\pier{\|}^2+\frac{1}{4}\pier{\|}\grad\grad \boldsymbol{\Theta}_m\pier{\|}^2+\frac{1}{2}\pier{\|}\overbigdot{{\vec{u}}}_m\pier{\|}^2+\psi_D^{\lambda}(\curl \mathbf{Z}_m)\\
& \displaystyle \notag +\int_0^{T_m}\biggl(\pier{\|}\overbigdot{\mathbf{W}}_m\pier{\|}^2+\pier{\|}\grad \overbigdot{\mathbf{W}}_m\pier{\|}^2+\pier{\|}\grad\grad \overbigdot{\mathbf{W}}_m\pier{\|}^2+\frac{1}{2}\pier{\|} \grad\overbigdot{\boldsymbol{\Theta}}_m\pier{\|}^2+\frac{1}{2}\pier{\|} \grad\grad\overbigdot{\boldsymbol{\Theta}}_m\pier{\|}^2\\
& \displaystyle \notag +\lambda\pier{\|}\twobigdot{\vec{u}}_m\pier{\|}_{\left(H_{\Gamma_D}^1(\mathcal{D}_a,\mathbb{R}^3)\right)'}^2\biggr)\leq C+C\int_0^{T_m}\biggl(\frac{1}{2}\pier{\|}\mathbf{W}_m-\mathbf{I}\pier{\|}^2+\frac{1}{2}\pier{\|}\grad \mathbf{W}_m\pier{\|}^2+\frac{1}{2}\pier{\|}\grad\grad \mathbf{W}_m\pier{\|}^2\\
&\displaystyle +\frac{1}{2}\pier{\|}\overbigdot{\mathbf{W}}_m\pier{\|}^2+\frac{1}{4}\pier{\|} {\boldsymbol{\Theta}}_m\pier{\|}^2+\frac{1}{4}\pier{\|}\grad \boldsymbol{\Theta}_m\pier{\|}^2+\frac{1}{4}\pier{\|}\grad\grad \boldsymbol{\Theta}_m\pier{\|}^2+\frac{1}{2}\pier{\|} \overbigdot{\boldsymbol{\Theta}}_m\pier{\|}^2+\frac{1}{2}\pier{\|}\grad \overbigdot{\boldsymbol{\Theta}}_m\pier{\|}^2\biggr),
\end{align}
Thanks to the Gronwall lemma, we thus have that
\begin{align}
\label{eqn:54}
& \displaystyle \notag \frac{1}{2}\pier{\|}\overbigdot{\mathbf{W}}_m\pier{\|}^2+\frac{1}{2}\pier{\|}\grad{} \overbigdot{\mathbf{W}}_m\pier{\|}^2+\frac{1}{2}\pier{\|}\grad\grad \overbigdot{\mathbf{W}}_m\pier{\|}^2+\frac{1}{2}\pier{\|}\mathbf{W}_m-\mathbf{I}\pier{\|}^2
+\pier{\tilde{\psi}_{C_{\alpha}}}(\mathbf{W}_m)\\
& \displaystyle \notag
+\frac{1}{2}\pier{\|}\grad \mathbf{W}_m\pier{\|}^2+\frac{1}{2}\pier{\|}\grad\grad \mathbf{W}_m\pier{\|}^2 +\frac{1}{2}\pier{\|} \overbigdot{\boldsymbol{\Theta}}_m\pier{\|}^2+\frac{1}{2}\pier{\|}\grad \overbigdot{\boldsymbol{\Theta}}_m\pier{\|}^2\\
& \displaystyle \notag+\frac{1}{2}\pier{\|}\grad\grad \overbigdot{\boldsymbol{\Theta}}_m\pier{\|}^2+\frac{1}{4}\pier{\|}\grad \boldsymbol{\Theta}_m\pier{\|}^2+\frac{1}{4}\pier{\|}\grad\grad \boldsymbol{\Theta}_m\pier{\|}^2+\frac{1}{2}\pier{\|}\overbigdot{{\vec{u}}}_m\pier{\|}^2+\psi_D^{\lambda}(\curl \mathbf{Z}_m)\\
& \displaystyle \notag +\int_0^{T_m}\biggl(\pier{\|}\overbigdot{\mathbf{W}}_m\pier{\|}^2+\pier{\|}\grad \overbigdot{\mathbf{W}}_m\pier{\|}^2+\pier{\|}\grad\grad \overbigdot{\mathbf{W}}_m\pier{\|}^2+\frac{1}{2}\pier{\|} \grad\overbigdot{\boldsymbol{\Theta}}_m\pier{\|}^2+\frac{1}{2}\pier{\|} \grad\grad\overbigdot{\boldsymbol{\Theta}}_m\pier{\|}^2\\
& \displaystyle +\lambda\pier{\|}\twobigdot{\vec{u}}_m\pier{\|}_{\left(H_{\Gamma_D}^1(\mathcal{D}_a,\mathbb{R}^3)\right)'}^2\biggr)\leq C,
\end{align}
where the constant in the right hand side of \eqref{eqn:54} depends only on the initial data, on the domain $\mathcal{D}_a$ and not on the discretization parameter $m$ and on the regularization parameter $\lambda$. As a consequence of \eqref{eqn:54}, we obtain that
\begin{align}
\label{wm}
&\mathbf{W}_m \;\; \text{is u.b. in} \;\; W^{1,\infty}(0,T_m;\bar{H}_{\Gamma_D}^2(\mathcal{D}_a,Sym(\mathbb{R}^{3\times 3})))\subset\subset C^0(\mathcal{D}_{aT_m},Sym(\mathbb{R}^{3\times 3}),\\
\label{wm2}
&\boldsymbol{\Theta}_m \;\; \text{is u.b. in} \;\; W^{1,\infty}(0,T_m;\bar{H}_{\Gamma_D}^2(\mathcal{D}_a,Skew(\mathbb{R}^{3\times 3})))\subset\subset C^0(\mathcal{D}_{aT_m},Skew(\mathbb{R}^{3\times 3}),\\
\label{wm3}
& \vec{u}_m \;\; \text{is u.b. in} \;\; W^{1,\infty}(0,T_m;L^2(\mathcal{D}_a,\mathbb{R}^{3})) \;\; \text{and $\lambda-$bounded in }\;\; H^2\left(0,T_m;{\left(H_{\Gamma_D}^1(\mathcal{D}_a,\mathbb{R}^3)\right)'}\right),
\end{align}
where ''u.b.'' stands form ''uniformly bounded'' and ''$\lambda-$bounded'' means that the bound depends on $\lambda$. 
We now integrate \eqref{eqn:53t0} in time over the interval $[0,t]$, for any $t\in [0,T_m]$, considering Assumption \textbf{A2, A3},  the property \eqref{eqn:ica2} and the fact that $\psi_D^{\lambda}(\curl \mathbf{Z}_0)=0$, obtaining in particular that
\begin{align*}
& \displaystyle \pier{\|}\mathbf{W}_m-\mathbf{I}\pier{\|}^2+\pier{\|}\grad \mathbf{W}_m\pier{\|}^2+\pier{\|}\grad\grad \mathbf{W}_m\pier{\|}^2\leq \int_{\mathcal{D}_{at}}\mathbf{W}_{\text{ext}}(\mathbf{W}_m,t)\colon \overbigdot{\mathbf{W}}_m\\
& \displaystyle +\int_{\mathcal{D}_{at}}\boldsymbol{\Omega}_{\text{ext}}\abramo{(\boldsymbol{\Theta}_m,t)}\colon \overbigdot{\boldsymbol{\Theta}}_m+\frac{1}{2}\int_{\mathcal{D}_{at}} {\boldsymbol{\Theta}}_m\colon \overbigdot{\boldsymbol{\Theta}}_m+\int_{\mathcal{D}_{at}}\mathrm{e}^{-\boldsymbol{\Theta}_m}\grad \left(\mathcal{G}_L\ast \vec{\mathcal{F}}_{\textrm ext}\right)\colon \overbigdot{\boldsymbol{\Theta}}_m\\
& \displaystyle +\int_{\mathcal{D}_{at}}\grad \left(\mathcal{G}_L\ast \vec{\mathcal{F}}_{\textrm ext}\right)\mathbf{W}_m\mathrm{e}^{-\boldsymbol{\Theta}_m}\colon \overbigdot{\boldsymbol{\Theta}}_m\leq Ct,
\end{align*}
where in the last inequality we employed \eqref{wm} and \eqref{wm2}.
Hence, the inequality \eqref{eqn:agmon} implies that
\begin{equation}
\label{wmc0}
|\mathbf{W}_m(t)-\mathbf{I}|_{C^0(\mathcal{D}_{a},Sym(\mathbb{R}^{3\times 3})}\leq C||\mathbf{W}_m(t)||_{H^1(\mathcal{D}_a,Sym(\mathbb{R}^{3\times 3})}^{\frac{1}{2}}||\mathbf{W}_m(t)||_{H^2(\mathcal{D}_a,Sym(\mathbb{R}^{3\times 3})}^{\frac{1}{2}}\leq C\sqrt{t},
\end{equation}
for any $t\in [0,T_m]$. Hence, since $\mathbf{W}_m\in \mathring{C}_{\alpha}$ in the time interval $[0,T_m]$, the derivative of the invariants of $\mathbf{W}_m$ with respect to $\mathbf{W}_m$ are uniformly bounded in $\mathcal{D}_{aT_m}$, and as a consequence we deduce that
\begin{equation}
\label{detwmc0}
|\text{det}(\mathbf{W}_m(t))-\text{det}(\mathbf{I})|_{C^0(\mathcal{D}_{a})}\leq C|\mathbf{W}_m(t)-\mathbf{I}|_{C^0(\mathcal{D}_{a},Sym(\mathbb{R}^{3\times 3})}\leq C\sqrt{t},
\end{equation}
\begin{equation}
\label{trwmc0}
|\text{Tr}(\mathbf{W}_m(t))-\text{Tr}(\mathbf{I})|_{C^0(\mathcal{D}_{a})}\leq C|\mathbf{W}_m(t)-\mathbf{I}|_{C^0(\mathcal{D}_{a},Sym(\mathbb{R}^{3\times 3})}\leq C\sqrt{t},
\end{equation}
\begin{equation}
\label{iimc0}
|\text{II}(\mathbf{W}_m(t))-\text{II}(\mathbf{I})|_{C^0(\mathcal{D}_{a})}\leq C|\mathbf{W}_m(t)-\mathbf{I}|_{C^0(\mathcal{D}_{a},Sym(\mathbb{R}^{3\times 3})}\leq C\sqrt{t},
\end{equation}
which together imply that, for each invariant $\Gamma(\mathbf{W}_m(t))\in \{I(\mathbf{W}_m(t)),II(\mathbf{W}_m(t)),II(\mathbf{W}_m(t))\}$, there exists a positive constant $C$ independent on $m$ such that
\[
\lambda_{I,II,III}-C\sqrt{t}\leq \Gamma(\mathbf{W}_m(t)) \leq \lambda_{I,II,II}+C\sqrt{t},
\]
where $\lambda_{I,II,III}=\{3,6,1\}$. Hence, there exists a $\hat{T}$ independent on $m$ such that
\[
\mathbf{x}^m\in \mathring{C}_{\alpha}^m \;\; \text{for}\;\; t\in[0,\hat{T}).
\]
We observe that the estimate \eqref{eqn:54} may be extended by continuity to the interval $[0,\hat{T})$.
\newline
Using \eqref{eqn:dpsi1} and \eqref{eqn:54} we have that
\begin{equation}
\label{eqn:542}
\sup_{t\in (0,\hat{T})}\pier{\|} (\curl \mathbf{Z}_m ) \pier{(t)} \pier{\|}^2\leq C.
\end{equation}
\pier{Moreover, in view of \eqref{eqn:dpsi2}, from \eqref{eqn:48m2}$_3$ and \eqref{eqn:54} it follows} that
\begin{equation}
 \label{eqn:543}
 \sup_{t\in (0,\hat{T})}\pier{\|}\boldsymbol{\Sigma}_m \pier{(t)} \pier{\|}^2\leq C.
\end{equation}
We now multiply \pier{the equality $\boldsymbol{\Sigma}_m=\partial \psi_D^{\lambda}(\curl \mathbf{Z}_m)$ in} \eqref{eqn:48m2}$_3$ by $\curl \left(\mathcal{G}_{L,\diver}\ast \left(\curl \boldsymbol{\Sigma}_m\right)\right)\in H^1(\mathcal{D}_a;\mathbb{R}^{3\times 3})$ and integrate over $\mathcal{D}_a$. Employing multiple integration by parts, the Cauchy--Schwarz and Young inequalities and \eqref{eqn:dpsi2}, we obtain that
\begin{align*}
\displaystyle & \int_{\mathcal{D}_a}\boldsymbol{\Sigma}_m:\curl \left(\mathcal{G}_{L,\diver}\ast \left(\curl \boldsymbol{\Sigma}_m\right)\right)=\int_{\mathcal{D}_a}\curl \boldsymbol{\Sigma}_m: \mathcal{G}_{L,\diver}\ast \left(\curl \boldsymbol{\Sigma}_m\right)\\
\displaystyle &=  \pier{\|}\curl \left(\mathcal{G}_{L,\diver}\ast \left(\curl \boldsymbol{\Sigma}_m\right)\right):\curl \left(\mathcal{G}_{L,\diver}\ast \left(\curl \boldsymbol{\Sigma}_m\right)\right)\pier{\|}^2\\
\displaystyle &= \int_{\mathcal{D}_a}\partial \psi_D^{\lambda}(\curl \mathbf{Z}_m):\curl \left(\mathcal{G}_{L,\diver}\ast \left(\curl \boldsymbol{\Sigma}_m\right)\right)\\
\displaystyle &\leq  C\psi_D^{\lambda}(\curl \mathbf{Z}_m)+C+\frac{1}{2}\pier{\|}\curl \left(\mathcal{G}_{L,\diver}\ast \left(\curl \boldsymbol{\Sigma}_m\right)\right):\curl \left(\mathcal{G}_{L,\diver}\ast \left(\curl \boldsymbol{\Sigma}_m\right)\right)\pier{\|}^2.
\end{align*}
Hence, given the estimate \eqref{eqn:54}, we have that 
\begin{equation}
 \label{eqn:544}
 \sup_{t\in (0,\hat{T})}\pier{\|}\curl \boldsymbol{\Sigma}_m\pier{(t)}\pier{\|}_{\left(H_{\Gamma_D,\diver}^1\left(\mathcal{D}_a,\mathbf{R}^{3\times 3}\right)\right)^{\prime}}^2\leq C,
\end{equation}
and, from a Lax--Milgram estimate associated to the operator $-P_L\Delta$, 
\begin{equation}
 \label{eqn:545}
\sup_{t\in (0,\hat{T})} \pier{\|}\pier{(\mathcal{G}_{L,\diver}\ast (\curl \boldsymbol{\Sigma}_m))}\pier{(t)}\pier{\|}_{H_{\Gamma_D,\diver}^1\left(\mathcal{D}_a,\mathbf{R}^{3\times 3}\right)}^2\leq C.
\end{equation} 
\subsubsection{Second a priori estimate}
The second a-priori estimate is obtained by taking $\nwhat{\mathbf{W}}_m=\twobigdot{\mathbf{W}}_m$ in \eqref{eqn:48m2}$_1$ and $\nwhat{\boldsymbol{\Omega}}_m=\twobigdot{\boldsymbol{\Theta}}_m$ in \eqref{eqn:48m2}$_2$. Moreover, we take the second time derivative of \eqref{eqn:48m2}$_4$ and $\vec{v}_m=\mathcal{G}_L\ast\twobigdot{\vec{u}}_m$. Finally, we sum all the contributions and integrate over $\mathcal{D}_{a}$. 
We obtain, collecting and rearranging some terms, that
\begin{align}
\label{ap21}
& \displaystyle \notag  \frac{d}{dt}\biggl(\frac{1}{2}\pier{\|}\overbigdot{\mathbf{W}}_m\pier{\|}^2+\frac{1}{2}\pier{\|}\grad \overbigdot{\mathbf{W}}_m\pier{\|}^2+\frac{1}{2}\pier{\|}\grad\grad \overbigdot{\mathbf{W}}_m\pier{\|}^2+\frac{1}{4}\pier{\|} \grad \overbigdot{\boldsymbol{\Theta}}_m\pier{\|}^2+\frac{1}{4}\pier{\|}\grad\grad \overbigdot{\boldsymbol{\Theta}}_m\pier{\|}^2\\
& \displaystyle \notag  +\frac{\lambda}{2}\pier{\|}\twobigdot{\vec{u}}_m\pier{\|}_{\left(H_{\Gamma_D}^1(\mathcal{D}_a,\mathbb{R}^3)\right)'}^2\biggr)+\pier{\|}\twobigdot{\mathbf{W}}_m\pier{\|}^2+\pier{\|}\grad \twobigdot{\mathbf{W}}_m\pier{\|}^2+\pier{\|}\grad\grad \twobigdot{\mathbf{W}}_m\pier{\|}^2+\frac{1}{2}\pier{\|} \twobigdot{\boldsymbol{\Theta}}_m\pier{\|}^2\\
& \displaystyle \notag+\frac{1}{2}\pier{\|} \grad \twobigdot{\boldsymbol{\Theta}}_m\pier{\|}^2+\frac{1}{2}\pier{\|} \grad\grad\twobigdot{\boldsymbol{\Theta}}_m\pier{\|}^2+\pier{\|}\twobigdot{\vec{u}}_m\pier{\|}^2 =-\int_{\mathcal{D}_a}\mathrm{e}^{-\boldsymbol{\Theta}_m}\curl\left(
\mathcal{G}_{L,\diver}\ast \curl \boldsymbol{\Sigma}_m\right)\colon \twobigdot{\mathbf{W}}_m\\
& \displaystyle \notag -\int_{\mathcal{D}_a}\left(\mathbf{W}_m-\mathbf{I}+\frac{\pier{{}\tilde{\psi}_{C_\alpha}}}{d \mathbf{W}_m}(\mathbf{W}_m)\right)\colon \twobigdot{\mathbf{W}}_m-\int_{\mathcal{D}_a}\pier{\grad{}} (\mathbf{W}_m)::\pier{\grad{}} \twobigdot{\mathbf{W}}_m\\
& \displaystyle \notag -\int_{\mathcal{D}_a}\pier{\grad\grad{}} (\mathbf{W}_m)::\pier{\grad\grad{}} \twobigdot{\mathbf{W}}_m+\int_{\mathcal{D}_a}\mathbf{W}_{\text{ext}}(\mathbf{W}_m,t)\colon \twobigdot{\mathbf{W}}_m\\
& \displaystyle \notag +\int_{\mathcal{D}_a}\mathrm{e}^{-\boldsymbol{\Theta}_m}\grad \left(\mathcal{G}_L\ast \vec{\mathcal{F}}_{\textrm ext}\right)\colon \twobigdot{\mathbf{W}}_m-\int_{\mathcal{D}_a}\left(\curl\left(
\mathcal{G}_{L,\diver}\ast \curl \boldsymbol{\Sigma}_m\right)\mathbf{W}_m\mathrm{e}^{-\boldsymbol{\Theta}_m}\right)\colon \twobigdot{\boldsymbol{\Theta}}_m\\
& \displaystyle \notag -\frac{1}{2}\int_{\mathcal{D}_a}\pier{\grad{}} \boldsymbol{\Theta}_m::\pier{\grad{}} \twobigdot{\boldsymbol{\Theta}}_m-\frac{1}{2}\int_{\mathcal{D}_a}\pier{\grad\grad{}}\boldsymbol{\Theta}\colon\pier{\grad\grad{}} \twobigdot{\boldsymbol{\Theta}}_m+\int_{\mathcal{D}_a}\boldsymbol{\Omega}_{\text{ext}}\abramo{(\boldsymbol{\Theta}_m,t)}\colon \twobigdot{\boldsymbol{\Theta}}_m\\
& \displaystyle \notag -\int_{\mathcal{D}_a}\grad \left(\mathcal{G}_L\ast \vec{\mathcal{F}}_{\textrm ext}\right)\mathbf{W}_m\mathrm{e}^{-\boldsymbol{\Theta}_m}\colon \twobigdot{\boldsymbol{\Theta}}_m+\int_{\mathcal{D}_a}\overbigdot{\boldsymbol{\Theta}}_m^2\mathrm{e}^{\boldsymbol{\Theta}_m}\mathbf{W}_m\colon \grad \left(\mathcal{G}_L\ast \twobigdot{\vec{u}}_m\right)\\
& \displaystyle +2\int_{\mathcal{D}_a}\overbigdot{\boldsymbol{\Theta}}_m\mathrm{e}^{\boldsymbol{\Theta}_m}\overbigdot{\mathbf{W}}_m\colon \grad \left(\mathcal{G}_L\ast \twobigdot{\vec{u}}_m\right)
\end{align}
Thanks to the Lax--Milgram estimate \eqref{eqn:1} and to \eqref{wm} and \eqref{wm2}, the last to terms on the right hand side of \eqref{ap21} can be controlled as
\begin{align*}
&\displaystyle \int_{\mathcal{D}_a}\overbigdot{\boldsymbol{\Theta}}_m^2\mathrm{e}^{\boldsymbol{\Theta}_m}\mathbf{W}_m\colon \grad \left(\mathcal{G}_L\ast \twobigdot{\vec{u}}_m\right)+2\int_{\mathcal{D}_a}\overbigdot{\boldsymbol{\Theta}}_m\mathrm{e}^{\boldsymbol{\Theta}_m}\overbigdot{\mathbf{W}}_m\colon \grad \left(\mathcal{G}_L\ast \twobigdot{\vec{u}}_m\right)\\
& \displaystyle \quad \leq C \left(\pier{\|}\overbigdot{\boldsymbol{\Theta}}_m^2\mathrm{e}^{\boldsymbol{\Theta}_m}\mathbf{W}_m\pier{\|}_{L^{\infty}(\mathcal{D}_a,\mathbb{R}^{3\times 3})}+(\pier{\|}\overbigdot{\boldsymbol{\Theta}}_m\mathrm{e}^{\boldsymbol{\Theta}_m}\overbigdot{\mathbf{W}}_m\pier{\|}_{L^{\infty}(\mathcal{D}_a,\mathbb{R}^{3\times 3})}\right)\pier{\|}\grad \left(\mathcal{G}_L\ast \twobigdot{\vec{u}}_m\right)\pier{\|}\\
& \displaystyle \quad \leq C +\frac{1}{2}\pier{\|}\twobigdot{\vec{u}}_m\pier{\|}.
\end{align*}
The remaining terms on the right hand side of \eqref{ap21} can be controlled in a similar way using Assumption \textbf{A3}, \eqref{eqn:54}, \eqref{wm}, \eqref{wm2}, \eqref{eqn:545}. Finally, using the Cauchy--Schwarz and Young inequalities, integrating in time over the interval $[0,\hat{T}]$ and using Assumption \textbf{A2}, we obtain that
\begin{align}
\label{eqn:55}
& \displaystyle \notag \frac{1}{2}\pier{\|}\overbigdot{\mathbf{W}}_m\pier{\|}^2+\frac{1}{2}\pier{\|}\grad \overbigdot{\mathbf{W}}_m\pier{\|}^2+\frac{1}{2}\pier{\|}\grad\grad \overbigdot{\mathbf{W}}_m\pier{\|}^2+\frac{1}{2}\pier{\|} \grad \overbigdot{\boldsymbol{\Theta}}_m\pier{\|}^2+\frac{1}{2}\pier{\|}\grad\grad \overbigdot{\boldsymbol{\Theta}}_m\pier{\|}^2\\
& \displaystyle \notag +\frac{\lambda}{2}\pier{\|}\twobigdot{\vec{u}}_m\pier{\|}_{\left(H_{\Gamma_D}^1(\mathcal{D}_a,\mathbb{R}^3)\right)'}^2+\int_0^{\hat{T}}\biggl(\pier{\|}\twobigdot{\mathbf{W}}_m\pier{\|}^2+\pier{\|}\grad \twobigdot{\mathbf{W}}_m\pier{\|}^2+\pier{\|}\grad\grad \twobigdot{\mathbf{W}}_m\pier{\|}^2+\frac{1}{2}\pier{\|} \twobigdot{\boldsymbol{\Theta}}_m\pier{\|}^2\\
& \displaystyle \notag +\frac{1}{2}\pier{\|} \grad \twobigdot{\boldsymbol{\Theta}}_m\pier{\|}^2+\frac{1}{2}\pier{\|} \grad\grad\twobigdot{\boldsymbol{\Theta}}_m\pier{\|}^2+\pier{\|}\twobigdot{\vec{u}}_m\pier{\|}^2\biggr)\\
&\displaystyle \notag \leq C+\frac{1}{2}\int_0^{\hat{T}}\biggl(\pier{\|}\twobigdot{\mathbf{W}}_m\pier{\|}^2+\pier{\|}\grad \twobigdot{\mathbf{W}}_m\pier{\|}^2+\pier{\|}\grad\grad \twobigdot{\mathbf{W}}_m\pier{\|}^2+\frac{1}{2}\pier{\|} \twobigdot{\boldsymbol{\Theta}}_m\pier{\|}^2\\
& \displaystyle+\frac{1}{2}\pier{\|} \grad \twobigdot{\boldsymbol{\Theta}}_m\pier{\|}^2+\frac{1}{2}\pier{\|} \grad\grad\twobigdot{\boldsymbol{\Theta}}_m\pier{\|}^2+\pier{\|}\twobigdot{\vec{u}}_m\pier{\|}^2\biggr)
\end{align}
 Thanks to \eqref{eqn:55}, we obtain that
\begin{align}
\label{2wm}
&\mathbf{W}_m \;\; \text{is u.b. in} \;\; W^{1,\infty}\cap H^2(0,\hat{T};\bar{H}_{\Gamma_D}^2(\mathcal{D}_a,Sym(\mathbb{R}^{3\times 3}))),\\
\label{2wm2}
&\boldsymbol{\Theta}_m \;\; \text{is u.b. in} \;\; W^{1,\infty}\cap H^2(0,\hat{T};\bar{H}_{\Gamma_D}^2(\mathcal{D}_a,Skew(\mathbb{R}^{3\times 3}))),\\
\label{2wm3}
& \vec{u}_m \;\; \text{is u.b. in} \;\; H^2(0,\hat{T};L^2(\mathcal{D}_a,\mathbb{R}^{3})) \;\; \text{and $\lambda-$bounded in }\;\; W^{2,\infty}\left(0,\hat{T};{\left(H_{\Gamma_D}^1(\mathcal{D}_a,\mathbb{R}^3)\right)'}\right),
\end{align}
\subsubsection{Higher order estimates for the displacement and the defect variables}
As a consequence of  \eqref{2wm} and \eqref{2wm2} we have that
\[
\overbigdot{\boldsymbol{\Theta}}_m\mathrm{e}^{\boldsymbol{\Theta}_m}\mathbf{W}_m+\mathrm{e}^{\boldsymbol{\Theta}_m}\overbigdot{\mathbf{W}}_m \;\; \text{is u.b. in} \;\; L^{\infty}\left(\mathcal{D}_{a\hat{T}},\mathbb{R}^{3\times 3}\right),
\]
and also that
\[
\twobigdot{\boldsymbol{\Theta}}_m\mathrm{e}^{\boldsymbol{\Theta}_m}\mathbf{W}_m+\overbigdot{\boldsymbol{\Theta}}_m^2\mathrm{e}^{\boldsymbol{\Theta}_m}\mathbf{W}_m+2\overbigdot{\boldsymbol{\Theta}}_m\mathrm{e}^{\boldsymbol{\Theta}_m}\mathbf{W}_m+\mathrm{e}^{\boldsymbol{\Theta}_m}\twobigdot{\mathbf{W}}_m \;\; \text{is u.b. in} \;\; L^2\left(0,\hat{T},L^{\infty}\left(\mathcal{D}_{a},\mathbb{R}^{3\times 3}\right)\right).
\]
Hence, taking the time derivative of \eqref{eqn:48m2}$_4$ and $\vec{v}_m=\overbigdot{\vec{u}}_m$ we obtain that
\begin{align*}
& \displaystyle \pier{\|}\grad \overbigdot{\vec{u}}_m\pier{\|}^2\leq C\pier{\|}\grad \overbigdot{\vec{u}}_m\pier{\|}-\lambda\int_{\mathcal{D}_a}\twobigdot{{\vec{u}}}_m\cdot \overbigdot{{\vec{u}}}_m\\
& \displaystyle \quad \leq \frac{C^2}{2\epsilon}+\frac{\epsilon}{2}\pier{\|}\grad \overbigdot{\vec{u}}_m\pier{\|}^2+\frac{\lambda}{2\epsilon}\pier{\|}\twobigdot{{\vec{u}}}_m\pier{\|}_{\left(H_{\Gamma_D}^1(\mathcal{D}_a,\mathbb{R}^3)\right)'}+\frac{\epsilon \lambda}{2}(1+C_P)\pier{\|}\grad \overbigdot{\vec{u}}_m\pier{\|}^2,
\end{align*}
where $C_P$ is the Poincar\'{e} constant and $\epsilon>0$. Choosing $\epsilon$ sufficiently small and using \eqref{eqn:55} we thus conclude that 
\begin{equation}
\label{umhigh1}
\vec{u}_m \;\; \text{is u.b. in} \;\; W^{1,\infty}(0,\hat{T};H_{\Gamma_D}^1(\mathcal{D}_a,\mathbb{R}^{3})).
\end{equation}
Moreover,  taking the second time derivative of \eqref{eqn:48m2}$_4$ and $\vec{v}_m=\twobigdot{\vec{u}}_m$, integrating in time over the interval $[0,\hat{T}]$, we obtain that
\[
\frac{\lambda}{2}\pier{\|}\twobigdot{\vec{u}}_m\pier{\|}+\frac{1}{2}\int_0^{\hat{T}}\pier{\|}\grad \twobigdot{\vec{u}}_m\pier{\|}^2\leq C,
\]
from which we conclude that 
\begin{equation}
\label{umhigh2}
\vec{u}_m \;\; \text{is u.b. in} \;\; H^2(0,\hat{T};H_{\Gamma_D}^1(\mathcal{D}_a,\mathbb{R}^{3})) \;\; \text{and $\lambda-$bounded in }\;\; W^{2,\infty}\left(0,\hat{T};L^2(\mathcal{D}_a,\mathbb{R}^3)\right).
\end{equation}
Finally, taking the $L^2$ scalar product of the first and second time derivative of \eqref{eqn:48m2}$_5$ with $\overbigdot{\mathbf{Z}}_m$ and $\twobigdot{\mathbf{Z}}_m$ respectively, we obtain that
\begin{equation}
\label{zmhigh}
\mathbf{Z}_m \;\; \text{is u.b. in} \;\; W^{1,\infty}(0,\hat{T};H_{\Gamma_N}^1(\mathcal{D}_a,\mathbb{R}^{3\times 3})\cap \mathcal{M}_{div})\cap H^2(0,\hat{T};H_{\Gamma_N}^1(\mathcal{D}_a,\mathbb{R}^{3\times 3})\cap \mathcal{M}_{div}).
\end{equation}

\subsubsection{Passing to the limit}
Collecting the \pier{bounds}~\eqref{wm}, \eqref{wm2}, \eqref{eqn:543}, \eqref{eqn:545}, \eqref{2wm}, \eqref{2wm2}, \eqref{umhigh1}, \eqref{umhigh2} and \eqref{zmhigh} which are uniform in $m$ and $\lambda$, from the Banach--Alaoglu, the Aubin--Lions and 
the \pier{Arzel{\`a}--Ascoli} lemmas, 
we finally obtain the convergence properties, up to subsequences
, which we still label by the index $m$ (without reporting the index $\lambda$), as follows:
\begin{align}
\label{eqn:convwm1} & \mathbf{W}_m \overset{\ast}{\rightharpoonup} \mathbf{W} \quad \text{in} \quad W^{1,\infty}(0,\hat{T};\bar{H}_{\Gamma_D}^2(\mathcal{D}_a;Sym(\mathbb{R}^{3\times 3}))),\\
\label{eqn:convwm2} &  \mathbf{W}_m \rightharpoonup \mathbf{W} \quad \text{in} \quad H^{2}(0,\hat{T};\bar{H}_{\Gamma_D}^2(\mathcal{D}_a;Sym(\mathbb{R}^{3\times 3}))),\\
\label{eqn:convwm3} & \mathbf{W}_m \to \mathbf{W} \quad \text{in} \quad C^{1}(\pier{[0,\hat{T}]};W^{1,p}(\mathcal{D}_a;Sym(\mathbb{R}^{3\times 3}))) , \;\; p\in [1,6), \;\; \text{and} \;\; \text{a.e. in} \; \; \mathcal{D}_{a\hat{T}},\\
\label{eqn:convwm4} & \mathbf{W}_m \to \mathbf{W} \quad \text{in} \quad C^{1}([0,\hat{T}];C^0(\overline{\mathcal{D}}_a;Sym(\mathbb{R}^{3\times 3}))), 
\end{align}
\begin{align}
\label{eqn:convtm1} &  \boldsymbol{\Theta}_m \overset{\ast}{\rightharpoonup}\boldsymbol{\Theta} \quad \text{in} \quad W^{1,\infty}(0,\hat{T};\bar{H}_{\Gamma_D}^2(\mathcal{D}_a;Skew(\mathbb{R}^{3\times 3}))),\\
\label{eqn:convtm2} &  \boldsymbol{\Theta}_m \rightharpoonup \boldsymbol{\Theta} \quad \text{in} \quad H^{2}(0,\hat{T};\bar{H}_{\Gamma_D}^2(\mathcal{D}_a;Skew(\mathbb{R}^{3\times 3}))),\\
\label{eqn:convtm3} & \boldsymbol{\Theta}_m \to \boldsymbol{\Theta} \quad \text{in} \quad C^{1}(\pier{[0,\hat{T}]};W^{1,p}(\mathcal{D}_a;Skew(\mathbb{R}^{3\times 3}))), \;\; p\in [1,6), \;\; \text{and} \;\; \text{a.e. in} \; \; \mathcal{D}_{a\hat{T}},\\
\label{eqn:convtm4} & \boldsymbol{\Theta}_m \to \boldsymbol{\Theta} \quad \text{in} \quad C^{1}(\pier{[0,\hat{T}]};C^0(\bar{\mathcal{D}}_a;Skew(\mathbb{R}^{3\times 3}))),\\
\label{eqn:convtm5} & \displaystyle \mathrm{e}^{\pm \boldsymbol{\Theta}_m} \to \mathrm{e}^{\pm \boldsymbol{\Theta}} \quad \text{uniformly in} \; \; \overline{\mathcal{D}}_{a\hat{T}},
\end{align}
\begin{align}
\label{eqn:convphim1} & \pier{\vec{u}}_m \overset{\ast}{\rightharpoonup} \pier{\vec{u}} \quad \text{in} \quad W^{1,\infty}(0,\hat{T};H_{\Gamma_D}^1(\mathcal{D}_a;\mathbb{R}^{3})),\\
\label{eqn:convphim2} & \pier{\vec{u}}_m \rightharpoonup \pier{\vec{u}} \quad \text{in} \quad H^{2}(0,\hat{T};H_{\Gamma_D}^{1}(\mathcal{D}_a;\mathbb{R}^{3})),\\
\label{eqn:convphim3} & \pier{\vec{u}}_m \to \pier{\vec{u}} \quad \text{in} \quad C^{1}(\pier{[0,T]};L^p(\mathcal{D}_a;\mathbb{R}^{3})), \;\; p\in [1,6), \;\; \text{and} \;\; \text{a.e. in} \; \; \mathcal{D}_{a\hat{T}},
\end{align}
\begin{align}
\label{eqn:convzm1} & \mathbf{Z}_m \overset{\ast}{\rightharpoonup} \mathbf{Z} \quad \text{in} \quad W^{1,\infty}(0,\hat{T};H_{\Gamma_N}^1(\mathcal{D}_a,\mathbb{R}^{3\times 3})\cap \mathcal{M}_{div}),\\
\label{eqn:convzm2} & \mathbf{Z}_m \rightharpoonup \mathbf{Z} \quad \text{in} \quad H^2(0,\hat{T};H_{\Gamma_N}^1(\mathcal{D}_a,\mathbb{R}^{3\times 3})\cap \mathcal{M}_{div}),\\
\label{eqn:convzm3} & \mathbf{Z}_m \to \mathbf{Z} \quad \text{in} \quad C^{1}(\pier{[0,\hat{T}]};L^{p}(\mathcal{D}_a;\mathbb{R}^{3\times 3})), \;\; p\in [1,6), \;\; \text{and} \;\; \text{a.e. in} \; \; \mathcal{D}_{a\hat{T}},
\end{align}
\begin{align}
\label{eqn:convsigmam1} & \boldsymbol{\Sigma}_m \overset{\ast}{\rightharpoonup} \boldsymbol{\Sigma} \quad \text{in} \quad L^{\infty}(0,T;L^2(\mathcal{D}_a;\mathbb{R}^{3\times 3})),\\
\label{eqn:convsigmam2} & \mathcal{G}_{L,\diver}(\curl \boldsymbol{\Sigma}_m) \overset{\ast}{\rightharpoonup} \mathcal{G}_{L,\diver}(\curl \boldsymbol{\Sigma}) \quad \text{in} \quad L^{\infty}(0,T;H^1(\mathcal{D}_a;\mathbb{R}^{3\times 3})),
\end{align}
as \pier{$m\to \infty$} and $\lambda \to 0$. We note that \eqref{eqn:convwm4} follows from \pier{\eqref{eqn:convwm3} and the compact embedding
$$ W^{1,p}(\mathcal{D}_a;Sym(\mathbb{R}^{3\times 3}))\subset C^0(\overline{\mathcal{D}}_a;Sym(\mathbb{R}^{3\times 3})),$$
holding for $p>3$.
Moreover, 
as 
$$ H^2(\mathcal{D}_a;Skew(\mathbb{R}^{3\times 3})) \, \hbox{ is compactly embedded into } \,C^0(\overline{\mathcal{D}}_a;Skew (\mathbb{R}^{3\times 3})),$$ 
the convergence \eqref{eqn:convtm2} implies a strong convergence in $C^{1}(\pier{[0,T]};C^0(\overline{\mathcal{D}}_a;Skew (\mathbb{R}^{3\times 3})))$, whence \eqref{eqn:convtm5} is easily deduced, thanks to the continuity of the exponential operator as well.}
\newline
With the convergence results \eqref{eqn:convwm1}--\eqref{eqn:convsigmam2}, we can pass to the limit in \pier{the system} \eqref{eqn:48m2} in a first step as $m\to \infty$. Let's take $\nwhat{\mathbf{W}}_m=PS_m(\nwhat{\mathbf{W}})$, $\nwhat{\boldsymbol{\Omega}}_m=PA_m(\nwhat{\boldsymbol{\Omega}})$, $\vec{v}_m=PV_m(\vec{v})$, with arbitrary $\nwhat{\mathbf{W}}\in \bar{H}_{\Gamma_D}^2(\mathcal{D}_a,Sym(\mathbb{R}^{3\times 3}))$, $\nwhat{\boldsymbol{\Omega}}\in \bar{H}_{\Gamma_D}(\mathcal{D}_a,Skew(\mathbb{R}^{3\times 3}))$, $\vec{v}\in H^1_{\Gamma_D}(\mathcal{D}_a,\mathbb{R}^3)$. Let's also rewrite \eqref{eqn:48m2}$_3$ using the convexity of $\psi_D^{\lambda}$ and the definition of the subdifferential, and moreover let us take the $L^2$ scalar product of \eqref{eqn:48m2}$_5$ with a function $\nwhat{\mathbf{Z}}\in H_{\Gamma_N}^1(\mathcal{D}_a,\mathbb{R}^{3\times 3})\cap \mathcal{M}_{div}$. We then multiply the equations by $\omega \in C_c^{\infty}([0,\hat{T}])$ and integrate over the time interval $[0,\hat{T}]$. This gives
\begin{equation}
\label{eqn:48m2cont1}
\begin{cases}
\displaystyle \int_0^{\hat{T}}\omega\int_{\mathcal{D}_a}\textrm{Sym}\left(\mathrm{e}^{-\boldsymbol{\Theta}_m}\grad \left(\mathcal{G}_L\ast \twobigdot{\vec{u}}_m\right)\right)\colon \nwhat{\mathbf{W}}_m+ \int_0^{\hat{T}}\omega\int_{\mathcal{D}_a}\twobigdot{\mathbf{W}}_m\colon \nwhat{\mathbf{W}}_m\\
\displaystyle{}+\int_0^{\hat{T}}\omega\int_{\mathcal{D}_a}\pier{\grad{}} \twobigdot{\mathbf{W}}_m::\pier{\grad{}} \nwhat{\mathbf{W}}_m+\int_0^{\hat{T}}\omega\int_{\mathcal{D}_a}\pier{\grad\grad{}} \twobigdot{\mathbf{W}}_m\colon \pier{\grad\grad{}} \nwhat{\mathbf{W}}_m\\
\displaystyle{}+\int_0^{\hat{T}}\omega\int_{\mathcal{D}_a}\textrm{Sym}\left(\mathrm{e}^{-\boldsymbol{\Theta}_m}\curl\left(
\mathcal{G}_{L,\diver}\ast \curl \boldsymbol{\Sigma}_m\right)\right)\colon \nwhat{\mathbf{W}}_m\\
\displaystyle{}+ \int_0^{\hat{T}}\omega\int_{\mathcal{D}_a}\left(\mathbf{W}_m-\mathbf{I}+\frac{\pier{{}\tilde{\psi}_{C_\alpha}}}{d \mathbf{W}_m}(\mathbf{W}_m)+\overbigdot{\mathbf{W}}_m\right)\colon \nwhat{\mathbf{W}}_m+\int_0^{\hat{T}}\omega\int_{\mathcal{D}_a}\pier{\grad{}} (\mathbf{W}_m+\overbigdot{\mathbf{W}_m})::\pier{\grad{}} \nwhat{\mathbf{W}}_m
\\
\displaystyle{}+\int_0^{\hat{T}}\omega\int_{\mathcal{D}_a}\pier{\grad\grad{}} (\mathbf{W}_m+\overbigdot{\mathbf{W}_m})\colon \pier{\grad\grad{}}\nwhat{\mathbf{W}}_m=\int_0^{\hat{T}}\omega\int_{\mathcal{D}_a}\mathbf{W}_{\text{ext}}(\mathbf{W}_m,t)\colon \nwhat{\mathbf{W}}_m\\
\displaystyle{}+\int_0^{\hat{T}}\omega\int_{\mathcal{D}_a}\textrm{Sym}\left(\mathrm{e}^{-\boldsymbol{\Theta}_m}\grad \left(\mathcal{G}_L\ast \vec{\mathcal{F}}_{\textrm ext}\right)\right)\colon \nwhat{\mathbf{W}}_m,\\ \\
\displaystyle \int_0^{\hat{T}}\omega\int_{\mathcal{D}_a}\textrm{Skew}\left(\grad \left(\mathcal{G}_L\ast \twobigdot{\vec{u}}_m\right)\mathbf{W}_m\mathrm{e}^{-\boldsymbol{\Theta}_m}\right)\colon \nwhat{\boldsymbol{\Omega}}_m+\int_0^{\hat{T}}\omega\int_{\mathcal{D}_a}\twobigdot{\boldsymbol{\Theta}}_m\colon \nwhat{\boldsymbol{\Omega}}_m\\
\displaystyle{}+\int_0^{\hat{T}}\omega\int_{\mathcal{D}_a}\grad{} \twobigdot{\boldsymbol{\Theta}}_m::\pier{\grad{}} \nwhat{\boldsymbol{\Omega}}_m+\int_0^{\hat{T}}\omega\int_{\mathcal{D}_a}\pier{\grad\grad{}} \twobigdot{\boldsymbol{\Theta}}_m\colon \pier{\grad\grad{}} \nwhat{\boldsymbol{\Omega}}_m\\
\displaystyle{}+\int_0^{\hat{T}}\omega\int_{\mathcal{D}_a}\textrm{Skew}\left(\left(\curl\left(
\mathcal{G}_{L,\diver}\ast \curl \boldsymbol{\Sigma}_m\right)\mathbf{W}_m\mathrm{e}^{-\boldsymbol{\Theta}_m}\right)\right)\colon \nwhat{\boldsymbol{\Omega}}_m \\
\displaystyle{}+\frac{1}{2}\int_0^{\hat{T}}\omega\int_{\mathcal{D}_a}\pier{\grad{}} (\boldsymbol{\Theta}_m+\overbigdot{\boldsymbol{\Theta}}_m)::\pier{\grad{}} \nwhat{\boldsymbol{\Omega}}_m+\frac{1}{2}\int_0^{\hat{T}}\omega\int_{\mathcal{D}_a}\pier{\grad\grad{}}(\boldsymbol{\Theta}+ \overbigdot{\boldsymbol{\Theta}}_m)\colon\pier{\grad\grad{}} \nwhat{\boldsymbol{\Omega}}_m\\
\displaystyle=\int_0^{\hat{T}}\omega\int_{\mathcal{D}_a}\boldsymbol{\Omega}_{\text{ext}}\abramo{(\boldsymbol{\Theta}_m,t)}\colon \nwhat{\boldsymbol{\Omega}}_m+\int_0^{\hat{T}}\omega\int_{\mathcal{D}_a}\textrm{Skew}\left(\grad \left(\mathcal{G}_L\ast \vec{\mathcal{F}}_{\textrm ext}\right)\mathbf{W}_m\mathrm{e}^{-\boldsymbol{\Theta}_m}\right)\colon \nwhat{\boldsymbol{\Omega}}_m,\\ \\
\displaystyle \int_0^{\hat{T}}\omega\int_{\mathcal{D}_a}(\curl \nwhat{\mathbf{Z}}-\curl \mathbf{Z}_m):\boldsymbol{\Sigma}_m+\int_0^{\hat{T}}\omega\int_{\mathcal{D}_a}\psi_D^{\lambda}(\curl \mathbf{Z}_m)\leq \int_0^{\hat{T}}\omega\int_{\mathcal{D}_a}\psi_D^{\lambda}(\curl \nwhat{\mathbf{Z}}),\\ \\
\displaystyle \lambda \int_0^{\hat{T}}\omega\int_{\mathcal{D}_a}\overbigdot{{\vec{u}}}_m\cdot \vec{v}_m+\int_0^{\hat{T}}\omega\int_{\mathcal{D}_a}\grad{}{{\vec{u}}}_m\colon \grad{}\vec{v}_m=\int_0^{\hat{T}}\omega\int_{\mathcal{D}_a}\left(\mathrm{e}^{\boldsymbol{\Theta}_m}\mathbf{W}_m-\mathbf{I}\right)\colon \grad{}\vec{v}_m,\\ \\
\displaystyle \int_0^{\hat{T}}\omega\int_{\mathcal{D}_a}\grad {\mathbf{Z}_m}:: \grad \nwhat{\mathbf{Z}}=\int_0^{\hat{T}}\omega\int_{\mathcal{D}_a} \left(\mathrm{e}^{\boldsymbol{\Theta_m}}\mathbf{W}_m-\mathbf{I}\right)\colon \curl \nwhat{\mathbf{Z}}.
\end{cases}
\end{equation}
We observe that 
\begin{equation}
\label{eqn:48m2cont2}
\begin{cases}
PS_m(\nwhat{\mathbf{W}})\to \nwhat{\mathbf{W}} \quad \text{in} \; \;\bar{H}_{\Gamma_D}^2(\mathcal{D}_a;Sym(\mathbb{R}^{3\times 3})),
\\[2mm]
PA_m(\nwhat{\boldsymbol{\Omega}})\to \nwhat{\boldsymbol{\Omega}}\quad \text{in} \; \;\bar{H}_{\Gamma_D}^2(\mathcal{D}_a;Skew(\mathbb{R}^{3\times 3})),\\[2mm]
PV_m(\vec{v}_m)\to \vec{v}\quad \text{in} \; \;\bar{H}_{\Gamma_D}^1(\mathcal{D}_a;\mathbb{R}^{3}),\\
\end{cases}
\end{equation}
as $m\to \infty$.
Thanks to \eqref{eqn:convwm4}, \eqref{eqn:convtm5}, \eqref{eqn:48m2cont2}$_1$ and \eqref{eqn:48m2cont2}$_2$, we have that 
\begin{align*}
& \displaystyle \omega\mathrm{e}^{\boldsymbol{\Theta}_m}\nwhat{\mathbf{W}}_m\to \omega\mathrm{e}^{\boldsymbol{\Theta}}\nwhat{\mathbf{W}}\quad \text{in} \; \; C^0([0,\hat{T}];H^{2}(\mathcal{D}_a,\mathbb{R}^{3\times 3})),\\
& \displaystyle \omega\nwhat{\boldsymbol{\Omega}}_m \mathrm{e}^{\boldsymbol{\Theta}_m}\mathbf{W}_m\to \omega\nwhat{\boldsymbol{\Omega}} \mathrm{e}^{\boldsymbol{\Theta}}\mathbf{W}\quad \text{in} \; \; C^0([0,\hat{T}];H^{2}(\mathcal{D}_a,\mathbb{R}^{3\times 3})).
\end{align*} 
Hence, using \eqref{eqn:convphim2} and \eqref{eqn:convsigmam2}, by the product of weak-strong convergence we can pass to the limit in the nonlinear coupling terms in \eqref{eqn:48m2cont1} and obtain that 
\begin{align*}
& \int_0^{\hat{T}}\omega\int_{\mathcal{D}_a}\textrm{Sym}\left(\mathrm{e}^{-\boldsymbol{\Theta}_m}\grad \left(\mathcal{G}_L\ast \twobigdot{\vec{u}}_m\right)\right)\colon \nwhat{\mathbf{W}}_m\\
& = \int_0^{\hat{T}}\omega\int_{\mathcal{D}_a}\grad \left(\mathcal{G}_L\ast \twobigdot{\vec{u}}_m\right)\colon \mathrm{e}^{\boldsymbol{\Theta}_m}\nwhat{\mathbf{W}}_m\\
& \to \int_0^{\hat{T}}\omega\int_{\mathcal{D}_a}\textrm{Sym}\left(\mathrm{e}^{-\boldsymbol{\Theta}}\grad \left(\mathcal{G}_L\ast \twobigdot{\vec{u}}\right)\right)\colon \nwhat{\mathbf{W}},
\end{align*}
as $m\to \infty$, and
\begin{align*}
 &\notag\int_0^T\omega \int_{\mathcal{D}_a}\textrm{Sym}\left(\mathrm{e}^{-\boldsymbol{\Theta}_m}\curl\left(\mathcal{G}_{L,\diver}\ast \left(
\curl \boldsymbol{\Sigma}_m\right)\right)\right)\colon \nwhat{\mathbf{W}}_m\\
&\notag
=\int_0^T\omega \int_{\mathcal{D}_a}\curl\left(\mathcal{G}_{L,\diver}\ast \left(
\curl \boldsymbol{\Sigma}_m\right)\right)\colon \mathrm{e}^{\boldsymbol{\Theta}_m}\nwhat{\mathbf{W}}_m\\
& \to \int_0^T\omega \int_{\mathcal{D}_a}\textrm{Sym}\left(\mathrm{e}^{-\boldsymbol{\Theta}}\curl\left(\mathcal{G}_{L,\diver}\ast \left(
\curl \boldsymbol{\Sigma}\right)\right)\right)\colon \nwhat{\mathbf{W}},
\end{align*}
as $m\to \infty$. Similarly, 
\begin{align*}
& \int_0^{\hat{T}}\omega\int_{\mathcal{D}_a}\textrm{Skew}\left(\grad \left(\mathcal{G}_L\ast \twobigdot{\vec{u}}_m\right)\mathbf{W}_m\mathrm{e}^{-\boldsymbol{\Theta}_m}\right)\colon \nwhat{\boldsymbol{\Omega}}_m\\
& = \int_0^{\hat{T}}\omega\int_{\mathcal{D}_a}\grad \left(\mathcal{G}_L\ast \twobigdot{\vec{u}}_m\right)\colon \nwhat{\boldsymbol{\Omega}}_m\mathrm{e}^{\boldsymbol{\Theta}_m\mathbf{W}_m}\\
& \to \int_0^{\hat{T}}\omega\int_{\mathcal{D}_a}\textrm{Skew}\left(\grad \left(\mathcal{G}_L\ast \twobigdot{\vec{u}}\right)\mathbf{W}\mathrm{e}^{-\boldsymbol{\Theta}}\right)\colon \nwhat{\boldsymbol{\Omega}}
\end{align*}
as $m\to \infty$, and
\begin{align*}
 &\notag\int_0^{\hat{T}}\omega\int_{\mathcal{D}_a}\textrm{Skew}\left(\left(\curl\left(
\mathcal{G}_{L,\diver}\ast \curl \boldsymbol{\Sigma}_m\right)\mathbf{W}_m\mathrm{e}^{-\boldsymbol{\Theta}_m}\right)\right)\colon \nwhat{\boldsymbol{\Omega}}_m\\
& \int_0^{\hat{T}}\omega\int_{\mathcal{D}_a}\left(\curl\left(
\mathcal{G}_{L,\diver}\ast \curl \boldsymbol{\Sigma}_m\right)\right)\colon \nwhat{\boldsymbol{\Omega}}_m\mathrm{e}^{\boldsymbol{\Theta}_m}\mathbf{W}_m\\
& \to \int_0^{\hat{T}}\omega\int_{\mathcal{D}_a}\textrm{Skew}\left(\left(\curl\left(
\mathcal{G}_{L,\diver}\ast \curl \boldsymbol{\Sigma}\right)\mathbf{W}\mathrm{e}^{-\boldsymbol{\Theta}}\right)\right)\colon \nwhat{\boldsymbol{\Omega}}.
\end{align*}
Similar calculations may be employed to calculate the limit of the last terms on the right hand sides of \eqref{eqn:48m2cont1}$_1$ and \eqref{eqn:48m2cont1}$_2$, those depending on the external force $\mathcal{F}_{ext}$, considering Assumption \textbf{A3}. The limit of all the other terms in \eqref{eqn:48m2cont1}$_1$ and \eqref{eqn:48m2cont1}$_2$ can be obtained straightforwardly, employing the weak and strong convergence results \eqref{eqn:convwm1}-\eqref{eqn:convtm5}, the smoothness of $\tilde{\psi}_{C_{\alpha}}$, Assumption \textbf{A3} and the Lebesgue convergence theorem to deal with the nonlinear terms. 
\newline
Thanks to \eqref{eqn:convwm4} and \eqref{eqn:convtm5} we have that
\begin{equation}
\label{strongtw}
\mathrm{e}^{\boldsymbol{\Theta}_m}{\mathbf{W}}_m\to \omega\mathrm{e}^{\boldsymbol{\Theta}}{\mathbf{W}}\quad \text{in} \; \; C^0([0,\hat{T}];C^0(\bar{\mathcal{D}}_a,\mathbb{R}^{3\times 3})),
\end{equation}
and using this result, together with \eqref{eqn:convphim1} and \eqref{eqn:convzm1}, it is straightforward to pass to the limit as $m\to \infty$ in \eqref{eqn:48m2cont1}$_4$ and \eqref{eqn:48m2cont1}$_5$.
We obtain the following limit system, as $m\to \infty$, \pier{in terms of the limit functions (restoring the index $\lambda$) $\mathbf{W}^{\lambda}, \, \boldsymbol{\Theta}^{\lambda}, \,\boldsymbol{\Sigma}^{\lambda},\, \pier{\vec{u}}^{\lambda},\, \mathbf{Z}^{\lambda} $:}
\begin{equation}
\label{eqn:48lambda}
\begin{cases}
\displaystyle \int_{\mathcal{D}_a}\textrm{Sym}\left(\mathrm{e}^{-\boldsymbol{\Theta}^{\lambda}}\grad \left(\mathcal{G}_L\ast \twobigdot{\vec{u}^{\lambda}}\right)\right)\colon \nwhat{\mathbf{W}}+ \int_{\mathcal{D}_a}\twobigdot{\mathbf{W}^{\lambda}}\colon \nwhat{\mathbf{W}}\\
\displaystyle{}+\int_{\mathcal{D}_a}\pier{\grad{}} \twobigdot{\mathbf{W}^{\lambda}}::\pier{\grad{}}  \nwhat{\mathbf{W}}+\int_{\mathcal{D}_a}\pier{\grad\grad{}} \twobigdot{\mathbf{W}^{\lambda}}\colon \pier{\grad\grad{}}  \nwhat{\mathbf{W}}\\
\displaystyle{}+\int_{\mathcal{D}_a}\textrm{Sym}\left(\mathrm{e}^{-\boldsymbol{\Theta}^{\lambda}}\curl\left(
\mathcal{G}_{L,\diver}\ast \curl \boldsymbol{\Sigma}^{\lambda}\right)\right)\colon  \nwhat{\mathbf{W}}\\
\displaystyle{}+ \int_{\mathcal{D}_a}\left(\mathbf{W}^{\lambda}-\mathbf{I}+\frac{\pier{{}\tilde{\psi}_{C_\alpha}}}{d \mathbf{W}^{\lambda}}(\mathbf{W}^{\lambda})+\overbigdot{\mathbf{W}^{\lambda}}\right)\colon \nwhat{\mathbf{W}}+\int_{\mathcal{D}_a}\pier{\grad{}} (\mathbf{W^{\lambda}}+\overbigdot{\mathbf{W}^{\lambda}})::\pier{\grad{}} \nwhat{\mathbf{W}}
\\
\displaystyle{}+\int_{\mathcal{D}_a}\pier{\grad\grad{}} (\mathbf{W}^{\lambda}+\overbigdot{\mathbf{W}^{\lambda}})\colon \pier{\grad\grad{}}\nwhat{\mathbf{W}}=\int_{\mathcal{D}_a}\mathbf{W}_{\text{ext}}(\mathbf{W}^{\lambda},t)\colon \nwhat{\mathbf{W}}\\
\displaystyle{}+\int_{\mathcal{D}_a}\textrm{Sym}\left(\mathrm{e}^{-\boldsymbol{\Theta}^{\lambda}}\grad \left(\mathcal{G}_L\ast \vec{\mathcal{F}}_{\textrm ext}\right)\right)\colon \nwhat{\mathbf{W}},\\ \\
\displaystyle \int_{\mathcal{D}_a}\textrm{Skew}\left(\grad \left(\mathcal{G}_L\ast \twobigdot{\vec{u}^{\lambda}}\right)\mathbf{W}^{\lambda}\mathrm{e}^{-\boldsymbol{\Theta}^{\lambda}}\right)\colon \nwhat{\boldsymbol{\Omega}}+\int_{\mathcal{D}_a}\twobigdot{\boldsymbol{\Theta}^{\lambda}}\colon \nwhat{\boldsymbol{\Omega}}\\
\displaystyle{}+\int_{\mathcal{D}_a}\grad{} \twobigdot{\boldsymbol{\Theta}^{\lambda}}::\pier{\grad{}} \nwhat{\boldsymbol{\Omega}}+\int_{\mathcal{D}_a}\pier{\grad\grad{}} \twobigdot{\boldsymbol{\Theta}^{\lambda}}\colon \pier{\grad\grad{}} \nwhat{\boldsymbol{\Omega}}\\
\displaystyle{}+\int_{\mathcal{D}_a}\textrm{Skew}\left(\left(\curl\left(
\mathcal{G}_{L,\diver}\ast \curl \boldsymbol{\Sigma}^{\lambda}\right)\mathbf{W}^{\lambda}\mathrm{e}^{-\boldsymbol{\Theta}^{\lambda}}\right)\right)\colon \nwhat{\boldsymbol{\Omega}} \\
\displaystyle{}+\frac{1}{2}\int_{\mathcal{D}_a}\pier{\grad{}} (\boldsymbol{\Theta}^{\lambda}+\overbigdot{\boldsymbol{\Theta}^{\lambda}})::\pier{\grad{}} \nwhat{\boldsymbol{\Omega}}+\frac{1}{2}\int_{\mathcal{D}_a}\pier{\grad\grad{}}(\boldsymbol{\Theta}^{\lambda}+ \overbigdot{\boldsymbol{\Theta}^{\lambda}})\colon\pier{\grad\grad{}} \nwhat{\boldsymbol{\Omega}}\\
\displaystyle=\int_{\mathcal{D}_a}\boldsymbol{\Omega}_{\text{ext}}\abramo{(\boldsymbol{\Theta}^{\lambda},t)}\colon \nwhat{\boldsymbol{\Omega}}+\int_{\mathcal{D}_a}\textrm{Skew}\left(\grad \left(\mathcal{G}_L\ast \vec{\mathcal{F}}_{\textrm ext}\right)\mathbf{W}^{\lambda}\mathrm{e}^{-\boldsymbol{\Theta}^{\lambda}}\right)\colon \nwhat{\boldsymbol{\Omega}},\\ \\
\displaystyle \lambda \int_{\mathcal{D}_a}\overbigdot{{\vec{u}}^{\lambda}}\cdot \vec{v}+\int_{\mathcal{D}_a}\grad{}{{\vec{u}}}^{\lambda}\colon \grad{}\vec{v}=\int_{\mathcal{D}_a}\left(\mathrm{e}^{\boldsymbol{\Theta}^{\lambda}}\mathbf{W}^{\lambda}-\mathbf{I}\right)\colon \grad{}\vec{v},\\ \\
\displaystyle \int_{\mathcal{D}_a}\grad {\mathbf{Z}^{\lambda}}:: \grad \nwhat{\mathbf{Z}}=\int_{\mathcal{D}_a} \left(\mathrm{e}^{\boldsymbol{\Theta^{\lambda}}}\mathbf{W}^{\lambda}-\mathbf{I}\right)\colon \curl \nwhat{\mathbf{Z}},
\end{cases}
\end{equation}
for a.e. $t \in [0,\hat{T}]$, for \pier{all choices of} $\nwhat{\mathbf{W}} \in \bar{H}_{\Gamma_D}^2(\mathcal{D}_a,Sym(\mathbb{R}^{3\times 3}))$, $\nwhat{\boldsymbol{\Omega}} \in \bar{H}_{\Gamma_D}(\mathcal{D}_a,Skew(\mathbb{R}^{3\times 3}))$, $\vec{v}\in H^1_{\Gamma_D}(\mathcal{D}_a,\mathbb{R}^3)$ and $\nwhat{\mathbf{Z}}\in H_{\Gamma_N}^1(\mathcal{D}_a,\mathbb{R}^{3\times 3})\cap \mathcal{M}_{div}$,
and with initial conditions \pier{(cf. the assumption \textbf{A2} and \eqref{eqn:48mic})}
\begin{equation}
    \label{eqn:48lic}
   \pier{ \mathbf{W}^\lambda(\cdot,0)=\mathbf{I}, \quad \boldsymbol{\Theta}^\lambda (\cdot,0)=\boldsymbol{0} \quad  \hbox{in }
   \, \mathcal{D}_a, \quad \vec{v}^{\lambda}(\cdot,0)=\vec{0}\quad  \hbox{a.e. in}
   \, \mathcal{D}_a.}
\end{equation}
In \pier{the system} \eqref{eqn:48lambda} we have restored the index $\lambda$, to indicate the dependence of the solutions from the regularization parameter $\lambda$. 
\newline
In order to deal with the limit of the inequality \eqref{eqn:48m2cont1}$_3$, we need to obtain a strong convergence result for $\curl \mathbf{Z}_m$. Let us take the $L^2$ scalar product of \eqref{eqn:48m2}$_5$ with $\mathbf{Z}_m$, obtaining that
\[
\pier{\|}\curl \mathbf{Z}^{\lambda}_m\pier{\|}^2=\int_{\mathcal{D}_a}\mathrm{e}^{\boldsymbol{\Theta^{\lambda}}_m}\mathbf{W}_m^{\lambda}\colon \curl \mathbf{Z}_m.
\]
Thanks to the convergence properties \eqref{strongtw} and \eqref{eqn:convzm1}, we have that
\[
\pier{\|}\curl \mathbf{Z}^{\lambda}_m\pier{\|}^2\to \int_{\mathcal{D}_a}\mathrm{e}^{\boldsymbol{\Theta^{\lambda}}}\mathbf{W}^{\lambda}\colon \curl \mathbf{Z}^{\lambda}
\]
a $m\to \infty$, and as a consequence of the weak formulation \eqref{eqn:48lambda} we conclude that
\begin{equation}
\label{strongz1}
\pier{\|}\curl \mathbf{Z}^{\lambda}_m\pier{\|}^2\to \pier{\|}\curl \mathbf{Z}^{\lambda}\pier{\|}^2.
\end{equation}
This result, together with the weak convergence \eqref{eqn:convzm1}, implies that
\begin{equation}
\label{strongz2}
\curl \mathbf{Z}^{\lambda}_m\to \curl \mathbf{Z}^{\lambda}\quad \text{in} \quad L^{\infty}(0,\hat{T};L^2(\mathcal{D}_a,\mathbb{R}^{3\times 3})\cap \mathcal{M}_{div}).
\end{equation}
Then, given the convergence results \eqref{eqn:convsigmam1} and \eqref{strongz2} we have that
\[
\int_0^{\hat{T}}\omega\int_{\mathcal{D}_a}(\curl \nwhat{\mathbf{Z}}-\curl \mathbf{Z}_m^{\lambda}):\boldsymbol{\Sigma}_m^{\lambda}\to \int_0^{\hat{T}}\omega\int_{\mathcal{D}_a}(\curl \nwhat{\mathbf{Z}}-\curl \mathbf{Z}^{\lambda}):\boldsymbol{\Sigma}^{\lambda}
\]
as $m\to \infty$. 
Moreover, using Fatou's lemma and the weak lower semicontinuity of $\psi_D^{\lambda}$, implied by its convexity and continuity, using also \eqref{eqn:convzm1}, we have that
\[
\int_0^{\hat{T}}\omega\int_{\mathcal{D}_a}\psi_D^{\lambda}(\curl \mathbf{Z}^{\lambda})\leq \int_0^{\hat{T}}\omega\int_{\mathcal{D}_a}\psi_D^{\lambda}(\curl \mathbf{Z}_m^{\lambda}).
\]
Then, in the limit as $m\to \infty$ the following inequality is valid:
\begin{equation}
 \label{lim3bis}
 \displaystyle \int_{\mathcal{D}_a}(\curl \nwhat{\mathbf{Z}}-\curl \mathbf{Z}^{\lambda}):\boldsymbol{\Sigma}^{\lambda}+\int_{\mathcal{D}_a}\psi_D^{\lambda}(\curl \mathbf{Z}^{\lambda})\leq \int_{\mathcal{D}_a}\psi_D^{\lambda}(\curl \nwhat{\mathbf{Z}}),
\end{equation}
for all $\nwhat{\mathbf{Z}}\in H_{\Gamma_N}^1(\mathcal{D}_a,\mathbb{R}^{3\times 3})\cap \mathcal{M}_{div}$ and a.e. $t\in [0,\hat{T}]$. 
\newline
We observe, without reporting all the details, that the estimates \eqref{eqn:54}, \eqref{eqn:542}, \eqref{eqn:543}, \eqref{eqn:545}, \eqref{eqn:55}, \eqref{umhigh1}, 
\eqref{umhigh2}, \eqref{zmhigh} and \eqref{strongz2} are preserved in the limit as $m\to \infty$, \pier{i.e.,} they are valid for the solutions of \pier{the system} \eqref{eqn:48lambda}. \pier{This allows us to pass to the limit as $\lambda \to 0$, up to subsequences of $\lambda$, in \pier{the system} \eqref{eqn:48lambda}, with similar calculations as the ones employed for the study of the limit problem as $m\to \infty$. On the other hand, thanks to the weak convergence
\[
\pier{\vec{u}}^{\lambda} \overset{\ast}{\rightharpoonup} \pier{\vec{u}} \quad \text{in} \quad W^{1,\infty}(0,\hat{T};H_{\Gamma_D}^1(\mathcal{D}_a;\mathbb{R}^{3})),
\] 
which is uniform in the parameter $\lambda$, we obtain that
\[
\lambda \int_{\mathcal{D}_a}\overbigdot{{\vec{u}}^{\lambda}}\cdot \vec{v}\to 0
\]
as $\lambda \to 0$. We finally obtain that the limit as $\lambda \to 0$ of system \eqref{eqn:48lambda} satisfies the system \eqref{eqn:thm1}, and the proof of Theorem \ref{thm:1} is completed.}
\section{Conclusions}
\label{sec:conclusions}
In this work we derived a model for large deformations and conditional compatibility, expressed in terms of the stretch and the rotation tensors as independent variables, which describes a viscoelastic solid subject to mixed boundary conditions, i.e. which is fixed only on a part of its boundary
and which is free to move on the other part. This model is a generalization of the model introduced in \cite{agosti4} for a viscoelastic solid subject to homogeneous Dirichlet boundary conditions and analytically studied in its quasi-stationary approximation.  
\newline
After the derivation of some technical results regarding Helmholtz-Hodge decomposition expressed in terms of elliptic problems and Green functions for elliptic operators with mixed boundary conditions, we derived the model from a generalized form of the principle of virtual powers, where the virtual velocities depend on the state variables as a consequence of internal \pierhhb{kinematic} constraints associated to the compatibility condition. The virtual velocities associated to the deformation and the defect variables were expressed through Green functions in terms of the virtual velocities associated to the stretch and the rotation tensors, thus reducing the set of independent virtual velocities and eliminating their internal constraints, obtaining a system of integro-differential coupled equations. 
The positive definiteness of the stretch matrix was imposed by adding to the free energy the indicator function of a closed and convex set whose elements are positive definite symmetric matrices with eigenvalues which are not smaller than a given positive constant at the same time. The \michhhb{internal forces} in the system were chosen in compliance with the Clausius--Duhem dissipative {\mich inequality}.
\newline
We developed the analysis of the full model with inertia. The inertia terms involve nonlinear couplings between the second order time derivative of the variables, which is the highest time derivative order in the equations, thus imposing challenges in the existence proof of a solution. Our strategy was to regularize the system, adding a time regularization in the kinematic constraints and employing the \michhhb{Moreau--Yosida} regularization of the subdifferential of the free energy for the defects associated to the threshold activation for the compatibility condition. These regularizations were both expressed in terms of a unique regularization parameter. Exploiting then a Faedo--Galerkin approximation of the regularized system and substituting the indicator function associated to the positive definiteness constraint for the stretch tensor with a smooth approximation from the interior of its proper domain, we proved the existence of a local in time weak solution in three space dimensions, studying the limit as the discretization parameter tends to zero and further as the regularization parameter tends to zero. The weak solution exists only locally in time, as long as the solution remains continuously in the interior of the proper domain of the indicator function associated to the positivity constraint, preceding the possible realization of external and internal collisions. This result is different from the one obtained in \cite{agosti4} in the quasi-stationary approximation of the model, i.e. neglecting inertia, where we obtained the global existence of strong solutions.
\newline
Further developments of the present work will be the study of the model with full incompatibility and with the possible presence of collisions in the dynamics, investigating the uniqueness and continuous dependence on data in these situations.


\section*{Acknowledgments}
\pier{A. A. acknowledges some support from the MIUR-PRIN Grant  P2022Z7ZAJ, from the GNAMPA (Gruppo Nazionale per l'Analisi Matematica,
la Probabilit\'{a} e le loro Applicazioni) of INdAM (Istituto Nazionale di Alta Matematica) through the GNAMPA project CUP E53C23001670001 and his affiliation to the GNAMPA. The authors wish to gratefully thanks Pierluigi Colli for fruitful advices and discussions.}

%
\bibliographystyle{plain}
\bibliography{biblio_CHVE} 
\end{document}